\documentclass[12pt,centertags,oneside]{amsart}
\usepackage{amsmath,amstext,amsthm,amscd}
\usepackage{amssymb}
\usepackage[mathscr]{eucal}
\usepackage{mathrsfs}
\usepackage{concmath}
\usepackage{charter}

\usepackage[a4paper,width=16.2cm,top=3cm,bottom=3cm,footskip=1cm]{geometry}

\numberwithin{equation}{section}

\newcommand{\field}[1]{\mathbb{#1}}

\newcommand{\R}{\field{R}}
\newcommand{\C}{\field{C}}
\newcommand{\N}{\field{N}}

 \def\cC{\mathscr{C}}
 \def\cF{\mathscr{F}}
 \def\cK{\mathscr{K}}
\def\cL{\mathscr{L}}
\def\cO{\mathscr{O}}
\def\cP{\mathscr{P}}

\def\mB{\mathcal{B}}

\def\mO{\mathcal{O}}

\def\mR{\mathcal{R}}

\newcommand{\boldsym}[1]{\boldsymbol{#1}}
\newcommand\bb{\boldsym{b}}

\newcommand\bI{\boldsym{I}}

\newcommand\br{\boldsym{r}}

\newcommand{\imat}{\sqrt{-1}}
\DeclareMathOperator{\End}{End}

\DeclareMathOperator{\Ker}{Ker}

\DeclareMathOperator{\rank}{rank}
\DeclareMathOperator{\Id}{Id}

\DeclareMathOperator{\tr}{Tr}

\DeclareMathOperator{\td}{Td}

\DeclareMathOperator{\ch}{ch}
\DeclareMathOperator{\ric}{Ric}
\DeclareMathOperator{\spec}{Spec}

\newcommand{\spin}{$\text{spin}^c$ }
\newcommand{\norm}[1]{\lVert#1\rVert}
\newcommand{\abs}[1]{\lvert#1\rvert}
\newcommand{\om}{\omega}
\newtheorem{thm}{Theorem}[section]
\newtheorem{lemma}[thm]{Lemma}

\newtheorem{cor}[thm]{Corollary}
\theoremstyle{definition}
\newtheorem{rem}[thm]{Remark}
\theoremstyle{definition}

\newtheorem{notation}[thm]{Notation}

\newcommand{\be}{\begin{eqnarray}}
\newcommand{\ee}{\end{eqnarray}}
\newcommand{\ov}{\overline}
\newcommand{\ovz}{\overline{z}}
\newcommand{\wi}{\widetilde}
\newcommand{\var}{\varepsilon}

\newcommand{\comment}[1]{}
\begin{document}

\title{Berezin-Toeplitz quantization on K{\"a}hler manifolds}
%\title{Berezin-Toeplitz quantization revisited}
\date{\today}
\author{Xiaonan Ma}
\address{Universit{\'e} Paris Diderot - Paris 7,
UFR de Math{\'e}matiques, Case 7012,
Site Chevaleret,
75205 Paris Cedex 13, France}
\email{ma@math.jussieu.fr}

\author{George Marinescu}
\address{Universit{\"a}t zu K{\"o}ln,  Mathematisches Institut,
    Weyertal 86-90,   50931 K{\"o}ln, Germany\\
    \& Institute of Mathematics `Simion Stoilow', Romanian Academy,
Bucharest, Romania}
%\thanks{}
\email{gmarines@math.uni-koeln.de}


\begin{abstract} We study Berezin-Toeplitz quantization
on K{\"a}hler manifolds. We explain first how to compute
various associated asymptotic expansions, then
we compute explicitly the first terms of the expansion
of the kernel of the  Berezin-Toeplitz operators,
and of the composition of two  Berezin-Toeplitz operators. 
As an application, we estimate the norm of Donaldson's $Q$-operator.
\end{abstract}

\dedicatory{Dedicated to Professor Jochen  Br\"uning on the occasion of his 65th birthday}

\maketitle
\setcounter{section}{-1}
%%%%%%%%%%%%%%%%%%%%%%%%%%%%%%%%%%%%%%%%%%%
%%%%%%%%%%%%%%%%%%%%%%%%%%%%%%%%%%%%%%%%%%%
\section{Introduction} \label{s1}
%%%%%%%%%%%%%%%%%%%%%%%%%%%%%%%%%%%%%%%%%%%

Berezin-Toeplitz operators are important in geometric quantization and
the properties of their kernels turn out to be 
deeply related to various problems in K{\"a}hler geometry 
(see e.g.\ \cite{Fine08,Fine10}). In this paper, 
we will study the precise asymptotic expansion of these kernels.
We refer the reader to the book \cite{MM07}
for a comprehensive study of the Bergman kernel,
Berezin-Toeplitz quantization and its applications. 
See also the survey \cite{Ma10}.

The setting of Berezin-Toeplitz quantization
on  K{\"a}hler manifolds is the following. 
Let $(X,\omega, J)$ be a compact K{\"a}hler manifold of
$\dim_{\C}X=n$ with K{\"a}hler form $\omega$ and complex structure $J$.
% Consider further 
Let $(L, h^L)$ be a holomorphic Hermitian line bundle on $X$, 
and let $(E, h^E)$ be a holomorphic Hermitian vector bundle on $X$.
Let $\nabla ^L, \nabla^E$ be the holomorphic Hermitian
connections on $(L,h^L)$, $(E, h^E)$  with curvatures
 $R^L=(\nabla^L)^2$, $R^E=(\nabla^E)^2$, respectively.
 We assume that $(L,h^L,\nabla^L)$ is a prequantum line bundle, i.e.,
%satisfies the {\em prequantization condition}, that is
 \begin{align} \label{toe2.1}
     \om= \frac{\sqrt{-1}}{2 \pi} R^L.
\end{align}

Let $g^{TX}(\cdot,\cdot):= \om(\cdot,J\cdot)$ be the Riemannian 
metric on $TX$
induced by $\om$ and $J$. The Riemannian volume
form $dv_X$ of $(X, g^{TX})$ has the form $dv_X= \om^n/n!$.
The $L^2$--Hermitian product on the space
$\cC^\infty(X,L^p\otimes E)$
of smooth sections of $L^p\otimes E$ on $X$, 
with $L^{p}:=L^{\otimes p}$, is given by
\begin{equation}\label{toe2.2}
\langle s_1,s_2\rangle=\int_X\langle s_1,s_2\rangle(x)\,dv_X(x)\,.
\end{equation}
We denote the corresponding norm by $\norm{\cdot}_{L^2}$ and
by $L^2(X,L^p\otimes E)$ the completion of
$\cC^\infty(X,L^p\otimes E)$ with respect to this norm.

Given a continuous smoothing linear operator 
$K:L^2(X,L^p\otimes E)\longrightarrow L^2(X,L^p\otimes E)$, 
the Schwartz kernel theorem \cite[Th.\,B.2.7]{MM07} 
guarantees the existence of an integral kernel with respect to $dv_X$, 
denoted by $K(x,x')\in(L^p\otimes E)_x\otimes(L^p\otimes E)^*_{x'}$, 
for $x,x'\in X$, i.e.,
\begin{equation}\label{toe2.21}
(KS)(x)=\int_X K(x,x')S(x')\,dv_X(x')\,,\quad S\in L^2(X,L^p\otimes E)\,.
\end{equation}
Consider now the space $H^0(X,L^p\otimes E)$ of holomorphic sections 
of $L^p\otimes E$ on $X$ and let
$P_p:L^2(X,L^p\otimes E)\to H^0(X,L^p\otimes E)$ 
be the orthogonal (Bergman) projection.
Its kernel $P_{p}(x,x')$ with respect to $dv_X(x')$ is smooth;
it is called the {\em Bergman kernel}.
%the orthogonal projection  $P_{p}$ from 
%$(\cC ^\infty(X, L^p\otimes E),\langle\cdot\,,\cdot \rangle)$ onto ,
The {\em Berezin-Toeplitz quantization}
of a section $f\in \cC^\infty(X, \End(E))$
is the {\em Berezin-Toeplitz operator} $\{T_{f,p}\}_{p\in \N}$ 
which is a sequence of linear operators $T_{f,\,p}$ defined by
 \begin{equation}\label{toe2.4}
 T_{f,\,p}:L^2(X,L^p\otimes E)\longrightarrow L^2(X,L^p\otimes E)\,,
 \quad T_{f,\,p}=P_p\,f\,P_p\,.
 \end{equation}
The kernel $T_{f,\,p}(x,x')$ of $T_{f,\,p}$ with respect to 
$dv_X(x')$ is also smooth. Since $\End(L)=\C$, 
we have $T_{f,\,p}(x,x)\in \End(E)_{x}$ for $x\in X$.
%kernel of the  Berezin-Toeplitz (BT) operator $T_{f,\,p}$  .

We introduce now the relevant geometric objects used in 
Theorems \ref{toet4.1}, \ref{toet4.6}  and \ref{toet4.5}.
%the study of the asymptotic behavior of the Berezin-Toeplitz operators. 
Let $T^{(1,0)}X$ be the holomorphic tangent bundle on $X$,
and $T^{*(1,0)}X$ its dual bundle.
Let $\nabla^{TX}$ be the Levi-Civita connection on $(X, g^{TX})$. We denote by
$R^{TX}=(\nabla^{TX})^2$ the curvature, by $\ric$ the Ricci curvature and 
by $\br$ the scalar curvature of $\nabla^{TX}$ (cf.\,\eqref{alm01.1}). 
%Let $\ric_\om= \ric(J\cdot,\cdot)$
%be the $(1,1)$-form associated to $\ric$.

We still denote by $\nabla ^{E}$  the connection on
 $\End (E)$ induced by $\nabla ^E$.
Consider %further 
the (positive) Laplacian $\Delta$ %on $(X, g^{TX})$
acting on the functions on $(X, g^{TX})$ and the %(positive) 
Bochner Laplacian $\Delta^{E}$
on $\cC^{\infty}(X, E)$ and on $\cC^{\infty}(X, \End(E))$. 
Let $\{e_k\}$ be a (local) orthonormal frame of $(TX, g^{TX})$. Then
\begin{align} \label{toe2.3}
    \Delta^{E} = - \sum_k (\nabla^{E}_{e_{k}}\nabla^{E}_{e_{k}}
- \nabla^{E}_{\nabla^{TX}_{e_{k}}e_{k}}).
\end{align}
%where $\{e_k\}$ is a (local) orthonormal frame of $(TX, g^{TX})$.

Let $\Omega^{q,\,r}(X, \End(E))$ be the space of $(q,r)$-forms on $X$ 
with values in $\End(E)$, and let
\begin{align} \label{abk2.5}
\nabla^{1,0}:\Omega^{q,*}(X, \End(E))\to \Omega^{q+1,*}(X,\End(E))
\end{align}
 be the $(1,0)$-component of the connection $\nabla^E$.
 Let $(\nabla^{E})^*$, $\nabla^{1,0*}, \ov{\partial}^{E*}$
  be the adjoints of $\nabla^{E}$,  $\nabla^{1,0}, \ov{\partial}^{E}$,
 respectively.
Let $D^{1,0}, D^{0,1}$ be the $(1,0)$ and $(0,1)$ components of the
 connection
$\nabla^{T^{*}X}: \cC^\infty(X,T^* X)\to \cC^\infty(X,T^* X\otimes T^*X)$
induced by $\nabla^{TX}$.

In the following, we denote by 
\[
\langle\cdot\,,\cdot \rangle_{\om}:
\Omega^{*,\,*}(X,\End(E))\times\Omega^{*,\,*}(X,\End(E))
\to\cC^\infty(X,\End(E))
\]
the $\C$-bilinear pairing
$\langle\alpha\otimes f,\beta\otimes g \rangle_{\om}
=\langle\alpha,\beta\rangle f\cdot g$, 
for forms $\alpha,\beta\in\Omega^{*,\,*}(X)$ and sections 
$f,g\in\cC^\infty(X,\End(E))$ (cf.  \eqref{abk4.4},
\eqref{toe4.31},  \eqref{toe4.32}). Put 
%In this sense, let $R^E_{\Lambda} $ be the contraction of $R^E$
 %with respect to $\om$, i.e.,
\begin{align} \label{bk2.5}
\begin{split}
&R^E_{\Lambda}
=\left \langle R^E, \om\right \rangle_{\om}\, .
\end{split}\end{align}
%We introduce now some coefficients appearing in our expansions. 
Let $\ric_\om= \ric(J\cdot,\cdot)$
be the $(1,1)$-form associated to $\ric$.
Set 
\[
 |\ric_\om|^2 = \sum_{i<j}\ric_\om(e_{i},e_{j})^{2}\,,
 \quad |R^{TX}|^2 = \sum_{i<j}\sum_{k<l} \langle 
R^{TX}(e_{i},e_{j})e_{k},e_{l}\rangle ^{2},
\] 
and let
\begin{equation}\label{abk2.6}\begin{split}
\bb_{2\C}= & - \frac{\Delta \br}{48}+ \frac{1}{96}|R^{TX}|^2 
- \frac{1}{24} |\ric_\om|^2 +  \frac{1}{128} \br^2,\\
\bb_{2E}=&\frac{\sqrt{-1}}{32} \Big(2 \br R^E_{\Lambda} 
- 4 \langle\ric_\om, R^ E\rangle_{\om} + \Delta^E  R^E_{\Lambda}\Big) \\
&- \frac{1}{8} (R^E_{\Lambda})^2 + 
 \frac{1}{8}\langle R^ E, R^ E\rangle_{\om}
+ \frac{3}{16} \ov{\partial}^{E*}\nabla^{1,0*} R^ E,\\
\bb_{1} =& \frac{\br}{8\pi}
+ \frac{\sqrt{-1}}{2\pi}R^E_{\Lambda},\quad  \,
\bb_{2} = \frac{1}{\pi^2}(\bb_{2\C} + \bb_{2E}).
\end{split}\end{equation}
We use now the notation from \eqref{lm01.4}.
By our convention (cf.\;\eqref{lm01.3}), we have at $x_{0}\in X$,
\[\big\langle\alpha_{\ell\ov{m}}\, dz_\ell\wedge  d\ov{z}_m,\beta_{k\ov{q}}\, 
dz_k\wedge  d\ov{z}_q\big\rangle
= -4 \alpha_{\ell\ov{m}}\,\beta_{m\ov{\ell}}\;\;,\quad
\big\langle \alpha_{\ov{m}\,\ov{q}}\, d\ov{z}_m\otimes  d\ov{z}_q,
\beta_{k\ell}\, dz_k\otimes dz_\ell\big\rangle
= 4  \alpha_{\ov{m}\,\ov{q}}\, \beta_{mq}
\]
(note that $|dz_q|^2=2$).
Then by Lemma \ref{lmt1.6}, \eqref{abk4.3} and \eqref{bk4.2a},
we have  at $x_{0}\in X$,\footnote{
The Bianchi identity reads $[\nabla^E, R^E]=0$. Take the derivative of 
$[\nabla^E, R^E]( \widetilde{\tfrac{\partial}{\partial z_i}}, 
\widetilde{\tfrac{\partial}{\partial z_k}}, 
\widetilde{\tfrac{\partial}{\partial \overline{z}_k}})=0,$
and use \eqref{lm01.5}, \eqref{lm01.27}, \eqref{bk2.85}
to obtain $[\widetilde{\tfrac{\partial}{\partial z_i}}, 
\widetilde{\tfrac{\partial}{\partial \ov{z}_k}}]=
-\frac{1}{2}R^{TX}_{x_{0}}(\tfrac{\partial}{\partial z_i}, 
\tfrac{\partial}{\partial \ov{z}_k}) \mR +  \cO(|Z|^{2})$,
$[\widetilde{\tfrac{\partial}{\partial z_i}}, 
\widetilde{\tfrac{\partial}{\partial z_k}}]=\cO(|Z|^{2})$.
We conclude that
$$\quad R^E_{m\overline{k}\,;\,k\overline{m}} 
= R^E_{k\overline{k}\,;\,m\overline{m}}\,.$$
The term $\frac{1}{4}\left(-R^E_{k\overline{k}\,;\,m\overline{m}}
+3 R^E_{m\overline{k}\,;\,k\overline{m}}\right)$ in (\ref{bk2.6})
can be thus replaced by
$\frac{1}{2}R^E_{k\overline{k}\,;\,m\overline{m}}$\,. Equivalently,
by using (\ref{bk4.3a}),
one can replace
$+\frac{\sqrt{-1}}{32}\Delta^E  R^E_{\Lambda}
+ \frac{3}{16} \overline{\partial}^{E*}\nabla^{1,0*} R^ E$
in (\ref{abk2.6}) by
$- \frac{\sqrt{-1}}{16}\Delta^E  R^E_{\Lambda}.$}
\begin{equation}\label{bk2.6}\begin{split}
&\ric_{\ell\ov{k}}=2R_{\ell\ov{k}q\ov{q}}=2R_{\ell\ov{q}q\ov{k}}\,, 
\quad \br=8R_{\ell\ov{\ell}q\ov{q}},\quad 
\ric_\om= \sqrt{-1} \ric_{\ell\ov{k}}\, dz_\ell\wedge d\ov{z}_k\,, \\ 
&\sqrt{-1}R^E_{\Lambda}= 2  R^E_{k\ov{k}},  \quad
\bb_{1} = \frac{1}{\pi}\left( R_{k\ov{k} m\ov{m}}+ R^E_{m\ov{m}}\right)\,,\\
&\bb_{2\C}= - \frac{\Delta \br}{48}
+ \frac{1}{6} R_{k\ov{\ell} m\ov{q}}  R_{\ell\ov{k}q\ov{m}}
- \frac{2}{3} R_{\ell\ov{\ell} m\ov{q}}  R_{k\ov{k}q\ov{m}}
+  \frac{1}{2} R_{\ell\ov{\ell}q\ov{q}} R_{k\ov{k}m\ov{m}}  \,,\\
&\bb_{2E}= R^E_{q\ov{q}}  R_{k\ov{k}m\ov{m}}
- R^E_{m\ov{q}}  R_{k\ov{k}q\ov{m}}
+ \frac{1}{2}\left(R^E_{q\ov{q}} R^E_{m\ov{m}} 
- R^E_{m\ov{q}} R^E_{q\ov{m}}\right)\\
&\hspace{15mm}+\frac{1}{4}\left(-R^E_{k\ov{k}\,;\,m\ov{m}}
+3 R^E_{m\ov{k}\,;\,k\ov{m}}\right)\,.
\end{split}\end{equation}
%
%\begin{equation}\label{bk2.6}\begin{split}
%&\ric_\om= \sqrt{-1} \ric_{\ell\ov{k}}\, dz_\ell\wedge d\ov{z}_k,  \quad
%\sqrt{-1}R^E_{\Lambda}= 2  R^E_{k\ov{k}},  \quad
%\bb_{1} = \frac{1}{\pi} ( R_{k\ov{k} m\ov{m}}+ R^E_{m\ov{m}}),\\
%&\bb_{2\C}= - \frac{\Delta \br}{48}
%+ \frac{1}{6} R_{k\ov{\ell} m\ov{q}}  R_{\ell\ov{k}q\ov{m}}
%- \frac{2}{3} R_{\ell\ov{\ell} m\ov{q}}  R_{k\ov{k}q\ov{m}}
%+  \frac{1}{2} R_{\ell\ov{\ell}q\ov{q}} R_{k\ov{k}m\ov{m}}  ,\\
%&\bb_{2E}= R^E_{q\ov{q}}  R_{k\ov{k}m\ov{m}}
%- R^E_{m\ov{q}}  R_{k\ov{k}q\ov{m}}
%+ \frac{1}{2} (R^E_{q\ov{q}} R^E_{m\ov{m}} - R^E_{m\ov{q}} R^E_{q\ov{m}})
%+\frac{1}{4}(-R^E_{k\ov{k}\,;\,m\ov{m}}+3 R^E_{m\ov{k}\,;\,k\ov{m}}).
%\end{split}\end{equation}
We say that a sequence $\Theta_{p}\in\cC^\infty(X,\End(E))$ 
has an asymptotic expansion of the form
\begin{equation}\label{bk4.2o}
\Theta_{p}(x)= \sum_{r=0}^{\infty} \boldsym{A}_{r}(x) p^{n-r}
+\cO(p^{-\infty})\,,
\quad \boldsym{A}_{r}\in\cC^\infty(X,\End(E))\,,
\end{equation}
if for any $k,l\in \N$, there exists
$C_{k,l}>0$ such that for any $p\in \N^*$, we have
\begin{equation}\label{bk4.23}
\Big |\Theta_{p}(x)- \sum_{r=0}^{k} \boldsym{A}_{r}(x) p^{n-r} \Big |_{\cC^l(X)}
\leqslant C_{k,l} \,p^{n-k-1},
\end{equation}
where $|\cdot|_{\cC^{l}(X)}$ is the $\cC^{l}$-norm on $X$.
%--------------------------
\begin{thm} \label{toet4.1}
For any $f\in\cC^\infty(X,\End(E))$, we have
%the kernel of the Toeplitz operator $T_{f,\,p}$ 
%has an asymptotic expansion on the diagonal
\begin{equation}\label{bk4.2}
T_{f,\,p}(x,x)= \sum_{r=0}^{\infty} \bb_{r,f}(x) p^{n-r}+\cO(p^{-\infty})\,,
\quad \bb_{r,f}\in\cC^\infty(X,\End(E))\,.
\end{equation}
%in the sense that
%for any $k,l\in \N$, there exists
%$C_{k,l}>0$ such that for any $p\in \N^*$, we have
%\begin{equation}\label{bk4.23}
%\Big |T_{f,\,p}(x,x)- \sum_{r=0}^{k} \bb_{r,f}(x) p^{n-r} \Big |_{\cC^l(X)}
%\leqslant C_{k,l} \,p^{n-k-1},
%\end{equation}
%where $|\cdot|_{\cC^{l}(X)}$ is the $\cC^{l}$-norm on $X$. 
Moreover,
 \begin{align}\label{bk4.3}
\bb_{0,f}=&f, \quad \bb_{1,f} = \frac{\br}{8\pi}  f
+ \frac{\sqrt{-1}}{4\pi}
\left(R^E_{\Lambda}  f + fR^E_{\Lambda} \right)
-  \frac{1}{4\pi} \Delta^{E} f.
\end{align}
If $f\in \cC^\infty(X)$, then
%\footnote{For $\alpha,\beta$
%  two $(1,1)$-forms, under the notation of \eqref{lm01.3}, if
%$\alpha=\alpha_{l\ov{j}} dz_l\wedge  d\ov{z}_j$,
%$\beta=\beta_{k\ov{i}} dz_k\wedge  d\ov{z}_i$, then
%$\langle\alpha,\beta\rangle= -4 \alpha_{l\ov{j}}\beta_{j\ov{l}}$.}
\begin{equation}\label{abk4.4}
\begin{split}
\pi^2 \bb_{2,f} & = \pi^2  \bb_{2}  f
+ \frac{1}{32}\Delta^2 f
- \frac{1}{32}  \br \Delta f
- \frac{\sqrt{-1}}{8} \big\langle \ric_\om, 
\partial\ov{\partial}f\big\rangle
+ \frac{\sqrt{-1}}{24} \big\langle df, 
\nabla^E R^E_{\Lambda}\big\rangle_{\om}\\
 &+ \frac{1}{24} \big\langle \partial f, \nabla^{1,0*} R^E\big\rangle_{\om}
 - \frac{1}{24} \big\langle \ov{\partial} f, 
\ov{\partial}^{E*} R^E\big\rangle_{\om}
-  \frac{\sqrt{-1}}{8} (\Delta f)R^E_{\Lambda}
+ \frac{1}{4}\big\langle \partial \ov{\partial} f, R^E\big\rangle_{\om}.
\end{split}
\end{equation}
\end{thm}
%\begin{rem}\label{toet4.2}

%For the kernel of the composition of the Berezin-Toeplitz operators,
%we have the following results.

\begin{thm} \label{toet4.6}
For any $f,g\in\cC^\infty(X,\End(E))$, the kernel of 
the composition $T_{f,\,p}\circ T_{g,\,p}$
has an asymptotic expansion on the diagonal
\begin{equation}\label{toe4.30}
(T_{f,\,p}\,\circ T_{g,\,p})(x,x)
=\sum_{r=0}^{\infty} \bb_{r,\,f,\,g}(x) p^{n-r}+\cO(p^{-\infty})\,,
\quad  \bb_{r,\,f,\,g}\in \cC^\infty(X,\End(E))\,,
\end{equation}
in the sense of \eqref{bk4.23}.
\comment{that for any $k,l\in \N$, there exists
$C_{k,l}>0$ such that the following estimate
\begin{equation}\label{toe4.301}
\Big | (T_{f,\,p}\, \cdot T_{g,\,p})(x,x)
- \sum_{r=0}^{k} \bb_{r,f,g}(x) p^{n-r} \Big |_{\cC^l(X)}
\leqslant C_{k,l} \, p^{n-k-1}
\end{equation}
holds for any $p\in \N$.}
Moreover, $\bb_{0,\,f,\,g}=fg$ and 
 \begin{equation}\label{toe4.31} \begin{split}
%\bb_{0,f,g}=&fg, \\
\bb_{1,\,f,\,g}= &\frac{1}{8\pi} \br fg
+  \frac{\sqrt{-1}}{4\pi}\big(R^E_{\Lambda}  fg + fgR^E_{\Lambda} \big)\\
&-\frac{1}{4\pi}\big(f\Delta^{E} g + (\Delta^{E} f)g\big) 
+ \frac{1}{2\pi}\big\langle\, \ov{\partial}^{E} f , \nabla^{1,0} g
\big\rangle_{\om}\; . %\in \cC^\infty(X,\End(E)).
\end{split}\end{equation}
If $f,g\in \cC^\infty(X)$, then
\begin{equation}\label{toe4.32}
\begin{split}
 \bb_{2,\,f,\,g} &= f\, \bb_{2,g} + g \,\bb_{2,f} - fg \,\bb_{2}
+ \frac{1}{\pi^2}\Big\{
-\frac{1}{8}\big\langle \,\ov{\partial} f, \partial \Delta g\big\rangle
 - \frac{1}{8}\big\langle \,\ov{\partial} \Delta f, \partial g\big\rangle\\
&+  \frac{1}{2}\big\langle \,\ov{\partial} f, \partial g\big\rangle \pi \bb_{1}
- \frac{1}{4}\big\langle \,\ov{\partial} f\wedge\partial g, R^E\big\rangle_{\om}
+ \frac{1}{16} \Delta f \cdot\Delta g
+ \frac{1}{8} \big\langle D^{0,1}\ov{\partial} f, 
D^{1,0}\partial g\big\rangle\Big\}.
\end{split}\end{equation}
\end{thm}
The existence of the expansions \eqref{bk4.2} and \eqref{toe4.30} 
and the formulas for
the leading terms hold in fact in the symplectic setting 
and are consequences of \cite[Lemma\,4.6 and (4.79)]{MM08b}
or \cite[Lemma\,7.2.4 and (7.4.6)]{MM07} (cf.\ Lemma \ref{toet2.3}).
The novel point of Theorems \ref{toet4.1}, \ref{toet4.6} is the calculation of 
the coefficients $\bb_{1,f}$, $\bb_{2,f}$,
 $\bb_{1,f,g}$ and $\bb_{2,f,g}$\,.
Note that the precise formula $\bb_{1,f}$ for a function $f\in \cC^\infty(X)$ 
was already given in \cite[Problem\,7.2]{MM07}. 
%The formula \eqref{abk4.4} is new even for $f\in\cC^\infty(X)$. 
In Theorem \ref{toet4.1a}, we find a general formula
of $ \bb_{2,f}$ for any $f\in \cC^\infty(X,\End(E))$.

If $f=1$, then $T_{f,\,p}=P_{p}$, and
the existence of the expansion (\ref{bk4.2}) and the form of the leading term
was proved by \cite{T90}, \cite{Catlin99}, \cite{Ze98}. 
The terms $\bb_{1}$, $\bb_{2}$ were computed by
Lu \cite{Lu00} (for $E=\C$, the trivial line bundle with trivial 
metric),  X.~Wang \cite{Wang05},
L.~Wang \cite{Wangl03}, in various degree of generality. 
The method of these authors is to construct appropriate peak sections 
as in \cite{T90}, using 
H{\"o}rmander's $L^2$ $\overline\partial$-method.
In \cite[\S 5.1]{DLM04a}, Dai-Liu-Ma computed  $\bb_{1}$
by using the heat  kernel, and in \cite[\S 2]{MM08a}, \cite[\S 2]{MM06}
(cf.\ also \cite[\S 4.1.8, \S 8.3.4]{MM07}), we 
%worked out the computation of 
computed $\bb_{1}$ in the symplectic case.
%The term $\bb_{2,f}$ for $f\in\cC^\infty(X,\End(E))$
%is more complicate, as we will only use \eqref{bk4.4} in
%the applications, thus we only compute precisely \eqref{bk4.4}.
%\end{rem}

%Note that 
%The expansion of the Bergman kernel was rederived
The expansion of the Bergman kernel $P_{p}(x,x)$ on the diagonal, 
for $E=\C$, was rederived
by Douglas and Klevtsov \cite{DouKle10} by using path integral 
and perturbation theory. They give physics interpretations
in terms of supersymmetric quantum mechanics, 
the quantum Hall effect and black holes (cf.\;also \cite{DouKle08}).

An interesting consequence of Theorem \ref{toet4.6}
is the following precise computation of the expansion of the composition
of two Berezin-Toeplitz operators.

\begin{thm} \label{toet4.5}
Let $f,g\in\cC^\infty(X,\End(E))$. The product of the Toeplitz operators
$T_{f,\,p}$ and  $T_{g,p}$ is a Toeplitz operator,
more precisely, it admits the asymptotic expansion
\begin{equation}\label{toe4.2}
T_{f,\,p}\circ T_{g,p}=\sum^\infty_{r=0}p^{-r}T_{C_r(f,g),p}
+\mO(p^{-\infty}),
\end{equation}
where $C_r$ are bidifferential operators, in the sense
that for any $k\geqslant0$, there exists
$c_k>0$ with
%----------------------------------------------------------------------------
\begin{equation}\label{toe4.1}
\Big\|T_{f,\,p}\circ  T_{g,\,p}- \sum_{l=0}^k p^{-l}T_{C_r(f,g),p}\Big\|
\leqslant c_k\,p^{-k-1},
\end{equation}
where $\norm{\cdot}$ denotes the operator norm on the space of
bounded operators.
We have
 \begin{equation}\label{toe4.3} \begin{split}
 C_0(f,g)=&fg, \\
C_1(f,g)=&-  \frac{1}{2\pi}
\langle \nabla^{1,0} f, \ov{\partial}^{E} g\rangle_{\om}\in 
\cC^{\infty}(X,\End(E)),\\
C_2(f,g)=& \, \bb_{2,f,g} - \bb_{2,fg}- \bb_{1,C_1(f,g)}.
\end{split}\end{equation}
%Here $\langle\cdot\,,\cdot \rangle$ acts $\C$-bilinearly
%(and pointwisely) on $TX$ and not on $E$.
If $f,g\in \cC^\infty(X)$, then
 \begin{align}\label{toe4.3a} \begin{split}
C_2(f,g)= &
\frac{1}{8\pi^2 } \langle  D^{1,0}\partial f, D^{0,1}\ov{\partial} g\rangle
+ \frac{\sqrt{-1}}{4\pi^2 } \langle  \ric_\om, 
\partial f \wedge\ov{\partial} g\rangle
 %+ \frac{\sqrt{-1}}{2\pi^2 }\langle \partial f, 
%\ov{\partial} g\rangle R^E_{\Lambda} 
 -\frac{1}{4\pi^2} \langle\partial f\wedge \ov{\partial} g, R^E\rangle_{\om}\,.
\end{split}\end{align}
\end{thm}
The existence of the expansion \eqref{toe4.2} is
a special case of \cite[Th.\;1.1]{MM08b} 
(cf.\;also \cite[Th.\;7.4.1,\;8.1.10]{MM07}),
where we found a symplectic version in which the Toeplitz operators 
\eqref{toe2.4} are constructed by using the projection to
the kernel of the Dirac operator.
Note that the precise values of $\bb_{2}$ are not used to derive 
\eqref{toe4.3a} (cf.\ Section \ref{toes4.3}).
%instead of $H^{0}(X,L^p\otimes E)$ considered here.

The existence of the expansion \eqref{toe4.2} 
for $E=\C$ 
was first established by
Bordemann, Meinrenken and Schlichenmaier \cite{BMS94},
Schlichenmaier \cite{Schlich:00} (cf.\ also  \cite{KS01})
using the theory of Toeplitz structures %(generalized Szeg{\"o} operators)
by Boutet de Monvel and Guillemin \cite{BouGu81}. 
Charles \cite{Charles03} calculated  $C_1(f,g)$ for $E=\C$. 
The asymptotic expansion \eqref{toe4.2} with a twisting bundle $E$ 
was derived by Hawkins \cite[Lemma 4.1]{Haw00} up to order one
(i.e., \eqref{toe4.1} for $k=0$).

Also, there is related work of Engli\v{s} \cite{Englis00,Englis02} dealing with 
expressing asymptotic expansions of Bergman kernel and coefficients of 
the Berezin-Toeplitz expansion \eqref{toe4.2} in terms of the metric. 
Engli\v{s} \cite[Cor.\,15]{Englis02} computed $C_{1}(f,g)$ and $C_{2}(f,g)$
for a smoothly bounded pseudoconvex domain 
$X=\{z\in\C^{n}:\varphi(z)>0\}$, where 
$\varphi$ is a defining function such that $-\log\varphi$ 
is strictly plurisubharmonic, and for the trivial line bundle $L=\C$ over $X$, 
equipped with the nontrivial metric $h^{L}=\varphi$ of positive curvature.

%The coefficient $C_1(f,g)$ was already calculated by Charles 
%\cite{Charles2} for the case of a trivial bundle $E$. The expression 
%for $C_2(f,g)$ given in \eqref{toe4.3} is new even in that case.

Note that we work throughout the paper with a non-trivial twisting bundle $E$.
%which is a generalization to manifolds of 
%the matrix-valued Berezin-Toeplitz quantization.
Moreover, we have shown in \cite[\S5-6]{MM08b} (cf.\ also \cite[\S 7.5]{MM07}) 
that Berezin-Toeplitz quantization
holds for complete K{\"a}hler manifolds and orbifolds endowed 
with a prequantum line bundle.
The calculations of the coefficients in the present paper 
being local in nature, they hold also for the above cases.

For some applications of the results of this paper 
to K\"ahler geometry see the paper \cite{Fine10} by Fine.

We close the introduction with some remarks about the Berezin-Toeplitz 
star-product. Following the ground-breaking work of Berezin \cite{Berez:74}, 
one can define a star-product by using Toeplitz operators. 
Note that formal star-products are known to exist on symplectic manifolds by
\cite{DeLe83, Fedo:96}. %\cite{Fedo:96,Guill:95}. 
The Berezin-Toeplitz star-product gives 
a very concrete and geometric realization of such product. 
For general symplectic manifolds this was realized in \cite{MM07,MM08b}
by using Toeplitz operators obtained by projecting on the kernel of 
the Dirac operator.

Consider now a compact K\"ahler manifold $(X,\omega)$ and 
a prequantum line bundle $L$. For every $f,g\in\cC^{\infty}(X)$
one defines the Berezin-Toeplitz star-product (cf.\ \cite{KS01,Schlich:00} 
and \cite{MM07,MM08b} for the symplectic case) by
\begin{equation}\label{toe4.4}
f*g:=\sum_{k=0}^\infty C_k(f,g) \hbar^{k}\in\cC^\infty(X)[[\hbar]].
\end{equation}
This star-product is associative. Moreover, for $f,g\in\cC^{\infty}(X)$ 
we have (cf.\ \cite[(7.4.3)]{MM07}, \cite[(4.89)]{MM08b})
 \begin{equation}\label{toe4.5}
 C_{0}(f,g)=fg=C_{0}(g,f)\,,\quad C_1(f,g)-C_1(g,f)=\imat\{f,g\}\,,
 \end{equation}
 where $\{f,g\}$ is the Poisson bracket associated to $2\pi\omega$. Therefore
 \begin{equation}\label{toe4.4a}
\big[T_{f,\,p}\,,T_{g,\,p}\big]=\tfrac{\sqrt{-1}}{\, p}\,T_{\{f,g\},\,p}
+\mO\big(p^{-2}\big)\,,\quad p\to\infty.
\end{equation}
Consider a twisting holomorphic Hermitian vector bundle $E$ 
and $f,g\in\cC^\infty(X,\End(E))$ as in Theorem \ref{toet4.5}. 
This corresponds to matrix-valued Berezin-Toeplitz quantization, 
which models a quantum system with $r=\rank E$ degrees of freedom.
By \eqref{toe4.2}, this Berezin-Toeplitz quantization has 
the expected semi-classical behaviour. Moreover, 
by \cite[Th.\,7.4.2]{MM07}, \cite[Th.\,4.19]{MM08b} we have
\begin{equation}\label{toe4.17}
\lim_{p\to\infty}\norm{T_{f,\,p}}={\norm f}_\infty
:=\sup_{0\neq u\in E_x, x\in X} |f(x)(u)|_{h^{E}}/ |u|_{h^{E}}\,.
\end{equation}
%We can thus define a star-product as in \eqref{toe4.4} for 
%$f,g\in\cC^\infty(X,\End(E))$. If the symbols $f,g$ of the Toeplitz operators 
%$T_{f,p}$, $T_{g,p}$ commute we obtain the following commutation relation.
\begin{cor}
Let $f,g\in\cC^\infty(X,\End(E))$. Set
\begin{equation}\label{toe4.4c}
f*g:=\sum_{k=0}^\infty C_k(f,g) \hbar^{k}\in\cC^\infty(X,\End(E))[[\hbar]].
\end{equation}
where $C_{r}(f,g)$ are determined by \eqref{toe4.2}. Then \eqref{toe4.4c} defines an associative star-product on $\cC^\infty(X,\End(E))$. Set moreover
\begin{equation}\label{toe4.5a}
\{\!\{f,g\}\!\}:=%\frac{1}{\imat}(C_1(f,g)-C_1(g,f))
\frac{1}{2\pi\imat}
\big(\langle \nabla^{1,0} g, \ov{\partial}^{E} f\rangle_{\om}-\langle \nabla^{1,0} f, \ov{\partial}^{E} g\rangle_{\om}\big)\,.
\end{equation}
If $fg=gf$ on $X$ we have 
\begin{equation}\label{toe4.4b}
\big[T_{f,\,p}\,,T_{g,\,p}\big]=\tfrac{\sqrt{-1}}{\, p}\,T_{\{\!\{f,g\}\!\},\,p}
+\mO\big(p^{-2}\big)\,,\quad p\to\infty.
\end{equation}
\end{cor}
The associativity of the star-product \eqref{toe4.4c} follows immediately 
from the associativity rule for the composition of Toeplitz operators, 
$(T_{f,\,p}\circ T_{g,\,p})\circ T_{k,\,p}=
T_{f,\,p}\circ (T_{g,\,p}\circ T_{k,\,p})$ for any  
$f,g,k\in\cC^\infty(X,\End(E))$, and from the asymptotic expansion 
\eqref{toe4.2} applied to both sides of the latter equality.

Due to the fact that $\{\!\{f,g\}\!\}=\{f,g\}$ if $E$ is trivial and 
comparing \eqref{toe4.4a} to \eqref{toe4.4b},  one can regard 
$\{\!\{f,g\}\!\}$ defined in \eqref{toe4.5a} 
as a non-commutative Poisson bracket.% by  \eqref{toe4.3}.
 
\begin{rem} \label{toet4.25}
Throughout the paper we suppose that 
%the Riemannian metric $g^{TX}$ was the metric associated to $\omega$, 
%that is, 
$g^{TX}(u,v)=\omega(u,Jv)$. 
The results presented so far still hold for a general non-K\"ahler 
$J$-invariant Riemannian 
metric $g^{TX}$. To explain this point we follow \cite[\S 4.1.9]{MM07}.

Let us denote the metric associated to $\omega$ by 
$g^{TX}_{\om}:=\om(\cdot,J\cdot)$.
We identify the $2$-form $R^L$ with the Hermitian matrix
 $\dot{R}^L \in \End(T^{(1,0)}X)$ via $g^{TX}$.
Then the Riemannian volume form of $g^{TX}_{\om}$ 
is given by $dv_{X,\,\om} =(2\pi)^{-n}{\det}(\dot{R}^L) dv_{X}$ 
(where $dv_{X}$ is the Riemannian volume form of $g^{TX}$).
Moreover, $h^E_\om:={\det} (\frac{\dot{R}^L}{2\pi})^{-1} h^E$ 
defines a metric on $E$. 
We add a subscript $\om$ to indicate the objects associated to $g^{TX}_{\om}$, $h^L$ and $h^E_{\om}$.
%Let $R^E_\om$ be the curvature associated 
%to the holomorphic Hermitian connection of $(E, h^E_\om)$, then 
%\begin{equation}\label{bk2.92}
%R^E_\om =R^E -\ov{\partial}\partial \log {\det}(\dot{R}^L).
%\end{equation}
Hence $\left\langle\,\cdot,\cdot \right \rangle_\om$ denotes 
the $L^2$ Hermitian product on 
$\cC^\infty (X, L^p\otimes E)$
induced by $g^{TX}_\om$, $h^L$, $h^E_\om$. 
This product is equivalent to the product  
$\left\langle\,\cdot,\cdot \right \rangle$ induced by $g^{TX}$, $h^L$, $h^E$.

Moreover, $H^{0}(X, L^p\otimes E)$ does not depend on the Riemannian metric 
on $X$ or on the Hermitian metrics on $L$, $E$. Therefore, the orthogonal 
projections from 
$(\cC^\infty (X, L^p\otimes E), \left\langle\,\cdot,\cdot \right \rangle_\om)$ and 
$(\cC^\infty (X, L^p\otimes E), \left\langle\,\cdot,\cdot \right \rangle)$
onto $H^{0}(X, L^p\otimes E)$ are the same. Hence $P_p=P_{p,\,\om}$ 
and therefore $T_{f,\,p}=T_{f,\,p,\,\om}$ as operators. 
However, their kernels are different.
If %$P_{p,\,\om}(x,x')$, 
$T_{f,\,p,\,\om}(x,x')$, ($x,x'\in X$), denotes the smooth 
kernels of %$P_{p,\,\om}$, 
$T_{f,\,p,\,\om}$ with respect to $dv_{X,\om}(x')$,
we have
\begin{equation}\label{bk2.95}
\begin{split}
%&P_{p}(x,x')=(2\pi)^{-n} {\det}(\dot{R}^L)(x') P_{p,\,\om}(x,x')\,,\\
&T_{f,\,p}(x,x')=(2\pi)^{-n} {\det}(\dot{R}^L)(x') T_{f,\,p,\,\om}(x,x')\,.
\end{split}
\end{equation}
For the kernel $T_{f,\,p,\,\om}(x,x')$, we can apply Theorem \ref{toet4.1} 
since $g^{TX}_\om(\cdot, \cdot)= \om(\cdot, J\cdot)$ 
is a K\"ahler metric on $TX$.
We obtain in this way the expansion of $T_{f,\,p}(x,x)$ for 
a non-K\"ahler metric $g^{TX}$ on $X$.
%see \cite[Th.\,4.1.1,\,4.1.3]{MM07}. 
%Using the expansion of the Bergman kernel $P_{p,\,\om}(\cdot,\cdot)$ 
%we can deduce the expansion of the Toeplitz operators $T_{f,\,p,\,\om}$ 
%and their kernels, analogous to Theorem \ref{toet2.3}, Corollary \ref{toec2.1} 
%and Theorem \ref{toet4.1b}. 
By \eqref{bk2.95}, the coefficients of these 
expansions (\ref{bk4.2}), (\ref{toe4.2}) satisfy 
\begin{equation}\label{bk2.951}
\begin{split}
&\bb_{r,\,f}=(2\pi)^{-n} {\det}(\dot{R}^L)\bb_{r,\, f,\,\om}\,,\\
&C_r(f,g)=C_{r,\,\om}(f,g)\,.
\end{split}
\end{equation}
\end{rem}

This paper is organized as follows. In Section~\ref{toes1}, we recall
the formal calculus on $\C^n$ for the model operator $\cL$, which is the
main ingredient of our approach. In Section~\ref{toes2}, we review the
asymptotic expansion of the kernel of Berezin-Toeplitz operators and
explain the strategy of our computation.
In Section \ref{toes3}, we obtain explicitly the first terms of the
Taylor expansion of our rescaled operator $\cL_t$\,.
In Section \ref{bks3}, we study in detail the contribution of
$\mO_2,\mO_4$ to the term $\cF_4$ from \eqref{bk2.32}.
In Section \ref{toes4},  by applying the formal calculi on $\C^n$ and
the results from Section \ref{bks3}, we establish Theorems
\ref{toet4.1}, \ref{toet4.6} and \ref{toet4.5}. We also verify that
our calculations are compatible with Riemann-Roch-Hirzebruch Theorem.
  In Section \ref{toes5}, we estimate the $\cC^m$-norm of
Donaldson's $Q$-operator, thus continuing
 \cite{LM07,LM09}.

%\smallskip

We shall use the following notations. For 
$\alpha=(\alpha_1, \ldots, \alpha_n)\in \N^{m}$,
$Z\in \C^{m}$, we set $|\alpha|:=\sum_{j=1}^{m}
\alpha_j$ and $Z^\alpha:= Z_1^{\alpha_1}\cdots Z_{m}^{\alpha_{m}}$.
Moreover, when an index variable appears twice in a single term,
it means that we are summing over all its possible values.

\smallskip
{\small\emph{
\textbf{Acknowledgments.} We would like to thank Joel Fine for motivating 
and helpful discussions which led to the writing of this paper.
We are grateful to the referee for pointing out some interesting references.
%The first-named author 
X.~M. thanks Institut Universitaire de France for support.
%The second-named author 
G.~M. was partially supported by SFB/TR 12 and 
Fondation Sciences Math{\'e}matiques de Paris. 
We were also supported by the DAAD Procope Program.
 }}

%%%%%%%%%%%%%%%%%%%%%%%%%%%%%%%%%%%%%%%%%%
%%%%%%%%%%%%%%%%%%%%%%%%%%%%%%%%%%%%%%%%%%
 \section{Kernel calculus on $\C^n$}\label{toes1}%\label{kernel calculus}
%%%%%%%%%%%%%%%%%%%%%%%%%%%%%%%%%%%%%%%%%%

In this section we recall
the formal calculus on $\C^n$ for the model operator $\cL$
introduced in  \cite[\S 2]{MM08b}, \cite[\S  7.1]{MM07}
(with $a_{j}=2\pi$ therein),
and we derive the properties of the calculus
of the kernels $(F\cP)(Z,Z^\prime)$, where $F\in\C[Z,Z^\prime]$ and
$\cP(Z,Z^\prime)$ is the kernel of the projection on the null space
of the model operator $\cL$.
This calculus is the main ingredient of our approach.

Let us consider the canonical coordinates $(Z_1,\dotsc,Z_{2n})$
on the real vector space $\R^{2n}$. On the complex vector space $\C^n$
we consider the complex coordinates $(z_1,\dotsc,z_n)$.
The two sets of coordinates are linked by the relation
$z_j=Z_{2j-1}+\imat Z_{2j}$, $j=1,\dotsc,n$.

We consider the $L^2$-norm
$\norm{\,\cdot\,}_{L^2}=\big(\int_{\R^{2n}}\abs{{\,\cdot\,}}^2\,dZ\big)^{1/2}$
on $\R^{2n}$, where $dZ=dZ_1\cdots dZ_{2n}$
is the standard Euclidean volume form.
We define the differential operators:
\begin{equation}\label{toe1.1}
\begin{split}
&b_i=-2{\frac{\partial}{\partial z
_i}}+\pi \overline{z}_i\,,\quad
b^{+}_i=2{\frac{\partial}{\partial\overline{z}_i}}+\pi
z_i\,,\quad b=(b_1,\ldots,b_n)\,,\quad
\cL=\sum_i b_i\, b^{+}_i\,,
\end{split}
\end{equation}
which extend to closed densely defined operators on 
$(L^2(\R^{2n}), \norm{\,\cdot\,}_{L^2})$. As such,
$b^{+}_i$ is the adjoint of $b_i$
and $\cL$ defines as a densely defined self-adjoint operator on
$(L^2(\R^{2n}),\norm{\,\cdot\,}_{L^2})$.

%Set
%\begin{equation}\label{toe1.1a}
%\cL=\sum_i b_i\, b^{+}_i\,.
%\end{equation}
%Then $\cL$ acts as a densely defined self-adjoint operator on
%$(L^2(\R^{2n}),\norm{\,\cdot\,}_{L^2})$.
%----------------------------------------------------------------------------
\noindent
The following result was established in \cite[Th.\,1.15]{MM08a}
(cf.\,also \cite[Th.\,4.1.20]{MM07}).
\begin{thm}\label{bkt2.17}
The spectrum of $\cL$ on $L^2(\R^{2n})$ is given by
\begin{equation}\label{bk2.68}
{\spec}(\cL)=\Big\{ 4\pi |\alpha| \,:\, \alpha\in\N ^n\Big\}\,.
\end{equation}
Each $\lambda\in\spec(\cL)$ is an eigenvalue of infinite multiplicity and 
an orthogonal basis of the corresponding eigenspace is given by
\begin{equation}\label{bk2.69}
B_\lambda=\Big\{b^{\alpha}\big(z^{\beta} e^{-\pi\sum_i |z_i|^2/2}\big):
\text{$\alpha\in\N^n$ with $4\pi|\alpha|=\lambda$, $\beta\in\N^n$}\Big\}
\end{equation}
and $\bigcup\{B_\lambda:\lambda\in\spec(\cL)\}$
forms a complete orthogonal basis of $L^2(\R^{2n})$.
In particular, an orthonormal basis of
$\Ker (\cL)$ is
\begin{equation}\label{bk2.70}
\Big\{\varphi_\beta(z)=\big(\tfrac{\pi ^{|\beta|}}{\beta!}\big)^{1/2}z^\beta
e^{-\pi\sum_i |z_i|^2/2}\,:\beta\in\N^n\Big\}\,.
\end{equation}

\end{thm}
%---------------------------------------------------------------------------
\noindent

Let $\cP(Z,Z')$ denote the kernel of the orthogonal projection
$\cP:L^2 (\R^{2n})\longrightarrow\Ker(\cL)$ with respect to $dZ'$.
Let $\cP^\bot = \Id -\cP$.
We call $\cP(\cdot,\cdot)$ the Bergman kernel of $\cL$.

Obviously
$\cP(Z,Z')=\sum_{\beta}\varphi_\beta(z)\,\overline{\varphi_\beta(z')}$ so
we infer from \eqref{bk2.70} that
\begin{equation}\label{toe1.3}
\cP(Z,Z') = \exp\big(-\tfrac{\pi}{2}\textstyle\sum_{i=1}^n
\big(|z_i|^2+|z^{\prime}_i|^2 -2z_i\overline{z}_i'\big)\big)\,.
\end{equation}

In the calculations involving the kernel $\cP(\cdot,\cdot)$,
 we prefer however to use the orthogonal decomposition of $L^2(\R^{2n})$
given in Theorem \ref{bkt2.17} and the fact that $\cP$ is
an orthogonal projection, rather than
integrating against the expression \eqref{toe1.3} of $\cP(\cdot,\cdot)$.
This point of view helps simplify a lot the computations
and understand better the operations.
As an example, if $\varphi(Z)=b^{\alpha}\big(z^{\beta} 
e^{-\pi\sum_i |z_i|^2/2}\big)$
%$\varphi(Z)= b^\alpha z^\beta\exp\Big (-\frac{\pi}{2} 
%\sum_{j=1}^n |z_j|^2\Big )$
with $\alpha,\beta\in \N^n$, then Theorem \ref{bkt2.17}
implies immediately that
 \begin{align}\label{toe1.4}
(\cP \varphi)(Z)= \left \{  \begin{array}{ll}
\displaystyle{z^\beta e^{-\pi\sum_i |z_i|^2/2}\,, }
& \mbox{if} \,\, |\alpha|=0,  \\
 0\,,& \mbox{if} \,\, |\alpha|>0.\end{array}\right.
\end{align}
\noindent
The following commutation relations are very useful in the computations.
Namely, for any polynomial  $g(z,\ov{z})$ in $z$ and $\ov{z}$, we have
\begin{align}\label{bk2.66}
\begin{split}
&[b_i,b^{+}_j]=b_i b^{+}_j-b^{+}_j b_i =-4\pi \delta_{i\,j},\\
&[b_i,b_j]=[b^{+}_i,b^{+}_j]=0\, ,\\
& [g(z,\ov{z}),b_j]=  2 \frac{\partial}{\partial z_j}g(z,\ov{z}), \\
&  [g(z,\ov{z}),b_j^+]
= - 2\frac{\partial}{\partial \overline{z}_j}g(z,\ov{z})\,.
\end{split}\end{align}
\noindent
For a polynomial $F$ in $Z,Z^\prime$, we denote by $F\cP$ 
the operator on $L^2(\R^{2n})$ defined by the kernel
$F(Z,Z^\prime)\cP(Z,Z^\prime)$ and the volume form $dZ$ 
according to \eqref{toe2.21}.

The following very useful Lemma \cite[Lemma 7.1.1]{MM07} 
describes the calculus of the kernels
$(F\cP)(Z,Z^\prime):=F(Z,Z^\prime)\cP(Z,Z^\prime)$.

%We will add a subscript $z$ or
%$z^{\prime}$ when we need to specify the operator is acting on the
%variables $Z$ or $Z^{\prime}$.
%---------------------------------------------------------------------
\begin{lemma}\label{toet1.1}
For any
$F,G\in\C[Z, Z^{\prime}]$ there exists a polynomial
$\cK[F,G]\in\C[Z, Z^{\prime}]$
with degree $\deg\cK[F,G]$ of the same parity as
$\deg F+\deg G$, such that
\begin{equation}\label{toe1.6}
((F\cP) \circ (G\cP))(Z, Z^{\prime}) =
\cK[F,G](Z, Z^{\prime}) \cP( Z, Z^{\prime}).
\end{equation}
\end{lemma}

%%%%%%%%%%%%%%%%%%%%%%%%%%%%%%%%%%%%%%%%%
%%%%%%%%%%%%%%%%%%%%%%%%%%%%%%%%%%%%%%%%%
\section{Expansion of the kernel of
Berezin-Toeplitz operators}\label{toes2}
%%%%%%%%%%%%%%%%%%%%%%%%%%%%%%%%%%%

In this section, we review some results from \cite{MM08a}, \cite{MM08b}
(cf.\ also \cite[\S\,7.2]{MM07}).  We explain then how to compute
the coefficients of various expansions considered in this paper.
We keep the notations and assumptions from the Introduction.

%Let's recall the basic setting and notation in \cite{LM07}.
%and we will use freely the notation in  \cite{LM07}.

\comment{
Let $(X,\omega, J)$ be a compact K{\"a}hler manifold of
$\dim_{\C}X=n$ with complex structure $J$,
and let $(L, h^L)$ be a holomorphic Hermitian line
bundle on $X$. Let $\nabla ^L$ be the holomorphic Hermitian
connection on $(L,h^L)$ with curvature $R^L$. We assume that
 \begin{align} \label{toe2.1}
\frac{\sqrt{-1}}{2 \pi} R^L=\om.
\end{align}

Let $(E, h^E)$ be a holomorphic Hermitian vector bundle on $X$. Let
$\nabla^E$ be the holomorphic Hermitian connection on $(E, h^E)$
with curvature $R^E$.
   Let $\nabla ^{\End (E)}$ be the connection on
 $\End (E)$ induced by $\nabla ^E$.

Let $g^{TX}(\cdot,\cdot):= \om(\cdot,J\cdot)$ be 
the Riemannian metric on $TX$ induced by $\om, J$. 
Let $dv_X$ be the Riemannian volume
form of $(TX, g^{TX})$, then $dv_X= \om^n/n!$.
Let $\Delta$ be the (positive) Laplace operator on $(X, g^{TX})$
acting on the functions on $X$.

The $L^2$--Hermitian product on the space
$\cC^\infty(X,L^p\otimes E)$
of smooth sections of $L^p\otimes E$ is given by
\begin{equation}\label{toe2.2}
\langle s_1,s_2\rangle=\int_X\langle s_1(x),s_2(x)\rangle\,dv_X(x)\,.
\end{equation}
We denote the corresponding norm with $\norm{\cdot}_{L^2}$ and
with $L^2(X,L^p\otimes E)$ the completion of
$\cC^\infty(X,L^p\otimes E)$ with respect to this norm.

Let $\nabla^{TX}$ be the Levi-Civita connection on $(X, g^{TX})$.
}

\noindent
\textbf{\emph{Kodaira-Laplace operator.\/}}
Let $\partial ^{L^{p}\otimes E,*}$ be the adjoint of the Dolbeault
operator $\partial ^{L^{p}\otimes E}$.
Let $\Box_{p}= \partial ^{L^{p}\otimes E,*}\partial ^{L^{p}\otimes E}$
be the restriction of the Kodaira Laplacian to
$\cC^{\infty}(X, L^p\otimes E)$.
Let $\Delta^{L^{p}\otimes E}$ be the Bochner Laplacian
on $\cC^{\infty}(X, L^p\otimes E)$ associated to $\nabla^{L}$, 
$\nabla^{E}$, $g^{TX}$, defined as in \eqref{toe2.3}.
\comment{Let
\begin{align} \label{toe2.3}
    \Delta^{L^{p}\otimes E}
= - \nabla^{L^{p}\otimes E}_{e_{i}}\nabla^{L^{p}\otimes E}_{e_{i}}
- \nabla^{L^{p}\otimes E}_{\nabla^{TX}_{e_{i}}e_{i}}
\end{align}
be the Bochner Laplacian
on $\cC^{\infty}(X, L^p\otimes E)$. }
Then we have (cf. \cite[Remark\,1.4.8]{MM07})
\begin{align}\label{toe2.5}
2 \, \Box_{p} = \Delta^{L^{p}\otimes E} 
- \tfrac{\sqrt{-1}}{2}R^{E}(e_{i}, Je_i) - 2n p.
\end{align}
Moreover, by Hodge theory (cf.\ \cite[Th.\,1.4.1]{MM07}) we have
\begin{align}\label{toe2.6}
\Ker (\Box_{p}) =H^0(X,L^p\otimes E).
\end{align}
This identification is important since the computations 
are performed by rescaling $\Box_{p}$
and expanding the rescaled operator.
%Our starting point is \eqref{toe2.6}, especially, 
%we get various computations under the help of
%the asymptotic expansion of the rescaled operator from $\Box_{p}$.

\comment{
Let $P_{p}(x,x')$ $(x,x'\in X)$ be the smooth kernel of
the orthogonal projection  $P_{p}$ from
$(\cC ^\infty(X, L^p\otimes E),\langle\cdot\,,\cdot \rangle )$
onto $H^0(X,L^p\otimes E)$,
the space of the holomorphic sections of $L^p\otimes E$ on $X$,
 with respect to $dv_X$. For $f\in \cC^\infty(X, \End(E))$,
its Berezin-Toeplitz quantization $T_{f,\,p}$ is defined by
 \begin{equation}\label{toe2.4}
 T_{f,\,p}:L^2(X,L^p\otimes E)\longrightarrow L^2(X,L^p\otimes E)\,,
 \quad T_{f,\,p}=P_p\,f\,P_p\,.
 \end{equation}
Let $T_{f,\,p}(x,x')$ $(x,x'\in X)$ be the smooth kernel of
$T_{f,\,p}$  with respect to $dv_X$.
 }

\noindent
\textbf{\emph{Normal coordinates.\/}}
 Let $a^X$ be the injectivity radius of $(X, g^{TX})$.
 We denote by $B^{X}(x,\var)$ and  $B^{T_xX}(0,\var)$ the open balls
 in $X$ and $T_x X$ with center $x$ and radius $\var$, respectively.
 Then the exponential map $ T_x X\ni Z \to \exp^X_x(Z)\in X$ is a
 diffeomorphism from $B^{T_xX} (0,\var)$ onto $B^{X} (x,\var)$ for
 $\var\leqslant a^X$.  {}From now on, we identify $B^{T_xX}(0,\var)$
 with $B^{X}(x,\var)$ via the exponential map for $\var \leqslant a^X$.
 Throughout what follows, $\varepsilon$ runs in the 
fixed interval $]0, a^X/4[$.

%   Let $\nabla ^{\End (\bE)}$ be the connection on
% $\End (\Lambda (T^{*(0,1)}X)\otimes E)$ induced by
 %$\nabla ^{\text{Cliff}}$ and $\nabla ^E$.

\noindent
\textbf{\emph{Basic trivialization.\/}}
 We fix $x_0\in X$.
 For $Z\in B^{T_{x_0}X}(0,\var)$ we identify
 $(L_Z, h^L_Z)$, $(E_Z, h^E_Z)$ and $(L^p\otimes E)_Z$
 to $(L_{x_0},h^L_{x_0})$, $(E_{x_0},h^E_{x_0})$ and $(L^p\otimes E)_{x_0}$
 by parallel transport with respect to the connections
 $\nabla ^L$, $\nabla ^E$ and $\nabla^{L^p\otimes E}$ along the curve
 $\gamma_Z :[0,1]\ni u \to \exp^X_{x_0} (uZ)$.
This is the basic trivialization we use in this paper.

Using this trivialization we identify $f\in \cC^\infty(X,\End(E))$ to a family
$\{f_{x_0}\}_{x_0\in X}$ where $f_{x_0}$ is the function $f$ in 
normal coordinates near $x_0$, i.e.,
$f_{x_0}:B^{T_{x_0}X}(0,\var)\to\End(E_{x_0})$, 
$f_{x_0}(Z)=f\circ\exp^X_{x_0}(Z)$.
%$f_{x_0}(Z)\in \End(E_{x_0})$ (with parameter $x_0\in X$)
%of functions in $Z$ in normal coordinates near $x_0$.
In general, for functions in the normal coordinates,
we will add a subscript $x_0$ to indicate the base point $x_0\in X$.
Similarly, 
$P_p(x,x')$ induces in terms of the basic trivialization a smooth section
$(Z,Z')\mapsto P_{p,\,x_0}(Z,Z')$
of $\pi ^* \End(E)$ over $\{(Z,Z')\in TX\times_{X} TX:|Z|,|Z'|<\var\}$, 
which depends smoothly on $x_0$. Here we identify a section 
$S\in \cC^\infty \big(TX\times_{X}TX,\pi ^* \End (E)\big)$
 with the family $(S_x)_{x\in X}$, where
 $S_x=S|_{\pi^{-1}(x)}$. %, $\End(E)=\C$

 \comment{
Then under our identification,
 we can view  $P_{p,x_0}(Z,Z')$ as a smooth section of $\pi^*(\End( E))$ on
  $TX\times_{X} TX$ (which is defined for $|Z|,|Z^\prime|\leq \var$)
 by identifying a section $S\in \cC^\infty (TX\times_{X}TX,\pi ^* \End (E))$
 with the family $(S_x)_{x\in X}$, where
 $S_x=S|_{\pi^{-1}(x)}$, $\End(E)=\C$.
}

\noindent
\textbf{\emph{Coordinates on $T_{x_0}X$.\/}}
Let us choose
an orthonormal basis $\{ w_i\}_{i=1}^n$ of $T^{(1,0)}_{x_0} X$.
Then $e_{2j-1}=\tfrac{1}{\sqrt{2}}(w_j+\overline{w}_j)$ and
$e_{2j}=\tfrac{\sqrt{-1}}{\sqrt{2}}(w_j-\overline{w}_j)$, $j=1,\dotsc,n\, $
form an orthonormal basis of $T_{x_0}X$.
We use coordinates on $T_{x_0}X\simeq\R^{2n}$
given by the identification 
\begin{equation}\label{n11}
\R^{2n}\ni (Z_1,\ldots, Z_{2n}) \longmapsto \sum_i
Z_i e_i\in T_{x_0}X.
\end{equation}
In what follows we also use complex coordinates $z=(z_1,\ldots,z_n)$
on $\C^n\simeq\R^{2n}$.

\noindent
\textbf{\emph{Volume form on $T_{x_0}X$.\/}}
If $dv_{TX}$ is the Riemannian volume form
 on $(T_{x_0}X, g^{T_{x_0}X})$, there exists
 a smooth positive function $\kappa_{x_0}:T_{x_0}X\to\R$, 
$Z\mapsto\kappa_{x_0}(Z)$ defined by
 \begin{equation} \label{atoe2.7}
 dv_X(Z)= \kappa_{x_0}(Z) dv_{TX}(Z),\quad \kappa_{x_0}(0)=1,
 \end{equation}
 where the subscript $x_0$ of
 $\kappa_{x_0}(Z)$ indicates the base point $x_0\in X$.

\noindent
\textbf{\emph{Sequences of operators.\/}}
Let $\Theta_p: L^2(X,L^p\otimes E)\longrightarrow L^2(X,L^p\otimes E)$ 
be a sequence of continuous linear  operators with smooth kernel 
$\Theta_p(\cdot,\cdot)$ with respect to $dv_X$ (e.g.\,$\Theta_p=T_{f,p}$).
 Let $\pi : TX\times_{X} TX \to X$ be the natural projection from the
 fiberwise product of $TX$ on $X$.
In terms of our basic trivialization, $\Theta_p(x,y)$ induces 
a family of smooth sections
$Z,Z'\mapsto \Theta_{p,\,x_0}(Z,Z^\prime)$
of $\pi^*\End(E)$ over $\{(Z,Z')\in TX\times_{X} TX:|Z|,|Z'|<\var\}$, 
which depends smoothly on $x_0$. %Here we identify a section 
 %$S\in \cC^\infty \big(TX\times_{X}TX,\pi ^* \End (E)\big)$
 % with the family $(S_x)_{x\in X}$, where
 % $S_x=S|_{\pi^{-1}(x)}$.

We denote by $\abs{\Theta_{p,\,x_0}(Z,Z^\prime)}_{\cC^l(X)}$ 
the $\cC^{l}$ norm with respect to the parameter $x_0\in X$. We say that
$\Theta_{p,\,x_0}(Z,Z^\prime)=\mO(p^{-\infty})$ if
for any $l,m\in \N$, there exists $C_{l,m}>0$ such that 
$\abs{\Theta_{p,\,x_0}(Z,Z^\prime)}_{\cC^{m}(X)}\leqslant C_{l,m}\, p^{-l}$.

\begin{notation}\label{noe2.7}
Recall that $\cP_{x_0}=\cP$ was defined in \eqref{toe1.3}. 
Fix $k\in\N$ and $\var^\prime\in\,]0,a^X[$\,. 
Let $\{Q_{r,\,x_0}\}_{0\leqslant r\leqslant k,x_0\in X}$  be a family
of polynomials $Q_{r,\,x_0}\in \End(E)_{x_0}[Z,Z^{\prime}]$
in $Z,Z^ \prime$, which is smooth with respect to
the parameter $x_0\in X$. We say that
\begin{equation} \label{toe2.7}
p^{-n} \Theta_{p,x_0}(Z,Z^\prime)\cong \sum_{r=0}^k
(Q_{r,\,x_0} \cP_{x_0})(\sqrt{p}Z,\sqrt{p}Z^{\prime})p^{-r/2}
+\mO(p^{-(k+1)/2})\,,
\end{equation}
on $\{(Z,Z^\prime)\in TX\times_X TX:\abs{Z},\abs{Z^{\prime}}<\var^\prime\}$ 
if there exist $C_0>0$ and a decomposition
\begin{equation} \label{toe2.71}
\begin{split}
p^{-n} \Theta_{p,x_0}(Z,Z^\prime)\kappa^{1/2}_{x_0}(Z)\kappa^{1/2}_{x_0}(Z')-\sum_{r=0}^k
(Q_{r,\,x_0} \cP_{x_0})(\sqrt{p}Z,\sqrt{p}Z^{\prime})p^{-r/2}\\
=\Psi_{p,k,x_0}(Z,Z^\prime)+\mO(p^{-\infty})\,,
\end{split}
\end{equation}
where $\Psi_{p,k,x_0}$ satisfies the following estimate on 
$\{(Z,Z^\prime)\in TX\times_X TX:\abs{Z},\abs{Z^{\prime}}<\var^\prime\}$:
for every $l\in\N$ there exist $C_{k,\,l}>0$, $M>0$ such that for all 
$p\in\N^{*}$
\begin{equation}
|\Psi_{p,k,x_0}(Z,Z^\prime)|_{\cC^l(X)}\leqslant  \,C_{k,\,l}\,p^{-(k+1)/2}
(1+\sqrt{p}\,|Z|+\sqrt{p}\,|Z^{\prime}|)^M \,
e^{-C_0\,\sqrt{ p}\,|Z-Z^{\prime}|}\,.
\end{equation}
 \end{notation}

\noindent
\textbf{\emph{The sequence $P_p$\,.\/}}
By \cite[Proposition\,4.1]{DLM04a} we know that the Bergman kernel 
decays very fast outside the diagonal of $X\times X$. 
Namely, for any $l,m\in \N$, $\var>0$, there exists $C_{l,m,\var}>0$
 such that for all $p\geqslant 1$ we have
 \begin{align}\label{1n13}
 &|P_{p}(x,x')|_{\cC^m}
 \leqslant C_{l,m,\var}\, p^{-l} \quad \text{on $\{(x,x')
 \in X\times X: d(x,x') \geqslant \var\}$} \,.
 \end{align}
 Here the $\cC^m$-norm is induced by $\nabla ^L$, $\nabla ^E$, 
 $\nabla ^{TX}$ and $h^L, h^E, g^{TX}$.

By \cite[Th.\,\,4.18$^\prime$]{DLM04a},
there exist polynomials $J_{r,\,x_{0}}(Z,Z')\in \End(E)_{x_0}$
in $Z,Z'$ with the same parity as $r$, such that
for any $k\in \N$, $\varepsilon\in]0, a^X/4[$\,,
%$Z,Z^\prime \in T_{x_0}X$, $\abs{Z},\abs{Z^{\prime}}<2 \var$,
we have
\begin{equation} \label{toe2.9}
p^{-n} P_{p,\,x_0}(Z,Z^\prime)\cong \sum_{r=0}^k
(J_{r,\,x_0} \cP_{x_0})(\sqrt{p}Z,\sqrt{p}Z^{\prime})p^{-\frac{r}{2}}
+\mO(p^{-\frac{k+1}{2}})\,,
\end{equation}
on the set $\{(Z,Z^\prime)\in TX\times_X TX:\abs{Z},\abs{Z^{\prime}}<2\var\}$,
in the sense of Notation \ref{noe2.7}.

\noindent
\textbf{\emph{The sequence $T_{f,\,p}$\,.\/}}
From (\ref{toe2.9}), we get the following result
(cf.\,\cite[Lemma\,4.6]{MM08b}, \cite[Lemma\,7.2.4]{MM07}).
%--------------------------------------------------------------------------
\begin{lemma} \label{toet2.3}
Let $f\in\cC^\infty(X,\End(E))$.
There exists a family $\{Q_{r,\,x_0}(f)\}_{r\in\N,\,x_0\in X}$, 
depending smoothly on the parameter $x_0\in X$, where
$Q_{r,\,x_0}(f)\in\End(E)_{x_0}[Z,Z^{\prime}]$
are polynomials with the same parity as $r$ and such that for every $k\in \N$,
$\varepsilon\in]0, a^X/4[$\,,
%$x_0\in X$, $Z,Z^\prime \in T_{x_0}X$, 
%$\abs{Z},\abs{Z^{\prime}}<\var/2$ we have
%--------------------------------------------------------------------------
\begin{equation} \label{toe2.13}
p^{-n}T_{f,\,p,\,x_0}(Z,Z^{\prime})
\cong \sum^k_{r=0}(Q_{r,\,x_0}(f)\cP_{x_0})(\sqrt{p}Z,\sqrt{p}Z^{\prime})
p^{-r/2} + \mO(p^{-(k+1)/2})\,,
\end{equation}
%--------------------------------------------------------------------------
on the set $\{(Z,Z^\prime)\in TX\times_X TX:\abs{Z},\abs{Z^{\prime}}<2\var\}$, 
in the sense of Notation \ref{noe2.7}. Moreover,
  $Q_{r,\,x_0}(f)$ are expressed by
  \begin{equation} \label{toe2.14}
  Q_{r,\,x_0}(f) = \sum_{r_1+r_2+|\alpha|=r}
    \cK\Big[J_{r_1,\,x_0}\;,\;
  \frac{\partial ^\alpha f_{\,x_0}}{\partial Z^\alpha}(0)
  \frac{Z^\alpha}{\alpha !} J_{r_2,\,x_0}\Big]\,.
  \end{equation}
  %--------------------------------------------------------------------------
  Especially,
  \begin{align} \label{toe2.15}
  Q_{0,\,x_0}(f)= f(x_0) .
  \end{align}
  \end{lemma}

Our goal is of course to compute the coefficients $Q_{r,\,x_0}(f)$. 
For this we need $J_{r,\,x_{0}}$, which are obtained by computing
the operators $\cF_{r,\,x_{0}}$ defined by the smooth kernels
\begin{align}\label{bk2.24}
    \cF_{r,\,x_{0}}(Z,Z')= J_{r,\,x_{0}}(Z,Z')\cP(Z,Z')
\end{align}
with respect to $dZ'$. Our strategy (already used in \cite{MM07,MM08a}) 
is to rescale the Kodaira-Laplace operator, take the Taylor expansion 
of the rescaled operator and apply resolvent analysis. 
In the remaining of this section we outline the main steps 
and continue the calculation in Section \ref{toes3}.

\noindent
\textbf{\emph{Rescaling $\Box_p$ and Taylor expansion.\/}}
For  $s \in \cC^{\infty}(\R^{2n}, E_{x_0})$, $Z\in \R^{2n}$,
 $|Z|\leq 2\var$, and
for $t=\frac{1}{\sqrt{p}}$, set
\begin{align}\label{bk2.21}
\begin{split}
&(S_{t} s ) (Z) :=s (Z/t),  \\
& \nabla'_{t}:=  S_t^{-1}
t \, \nabla ^{L^p\otimes E}  S_t,\\
& \nabla_{t}:=  S_t^{-1}
t \, \kappa^{1/2}\nabla ^{L^p\otimes E} \kappa^{-1/2} S_t
= \kappa^{1/2}(tZ) \nabla'_{t}\,  \kappa^{-1/2}(tZ), \\
&  \cL_{t}:= S_t^{-1} \kappa^{1/2}\, t^2 (2\, \Box_{p})\kappa^{-1/2} S_t.
\end{split}\end{align}
Then by \cite[Th.\,4.1.7]{MM07},  there exist
second order differential operators $\mO_{r}$
such that we have an asymptotic expansion in $t$ when $t\to 0$,
\begin{align}\label{bk2.22}
    \cL_{t} = \cL_{0} + \sum_{r=1}^{m} t^r \mO_{r} + \cO(t^{m+1}).
\end{align}
\noindent
From \cite[Th.\,\,4.1.21,\,4.1.25]{MM07}
(cf. also Theorem \ref{bkt2.21}), we obtain
\begin{align}\label{bk2.30}
    \cL_0=& \sum_j b_jb^+_j=\cL, \qquad \mO_1=0.
\end{align}

\noindent
\textbf{\emph{Resolvent analysis.\/}}
We define by recurrence
$f_r(\lambda)\in \End(L^2(\R^{2n}, E_{x_0}))$ by
\begin{align}\label{bk2.23}
f_{0}(\lambda) = (\lambda -\cL_0)^{-1},
\quad f_r(\lambda)=(\lambda -\cL_0)^{-1}
 \sum_{j=1}^{r} \mO_j   f_{r-j}(\lambda).
\end{align}
Let $\delta$ be the counterclockwise
oriented circle in $\C$ of center $0$ and radius $\pi/2$.
Then by \cite[(1.110)]{MM08a} (cf. also \cite[(4.1.91)]{MM07})
\begin{equation}\label{bk2.77}
\cF_{r,\,x_{0}}= \frac{1}{2\pi \sqrt{-1}} \int_{\delta}  f_r (\lambda)d \lambda.
\end{equation}

Since the spectrum of $\cL$ is well understood we can calculate  
the coefficients $\cF_{r,\,x_{0}}$. Recall than $\cP^\bot= \Id - \cP$.
From Theorem \ref{bkt2.17}, \eqref{bk2.30} and \eqref{bk2.77}, we get
\begin{align}\label{bk2.31}
    \begin{split}
	\cF_{0,\,x_{0}}=& \cP,  \quad \cF_{1,\,x_{0}}= 0,\\
\cF_{2,\,x_{0}}=&- \cL^{-1} \cP^\bot\mO_2 \cP
- \cP \mO_2\cL^{-1}\cP^\bot,\\
\cF_{3,\,x_{0}}=&- \cL^{-1} \cP^\bot\mO_3\cP
- \cP \mO_3\cL^{-1}\cP^\bot,
\end{split}\end{align}
and
\begin{align}\nonumber
    \cF_{4,\,x_{0}}&= \frac{1}{2 \pi \sqrt{-1}} \int_{\delta}
  \Big[(\lambda -\cL)^{-1} \cP^\bot (\mO_2 f_2 + \mO_4 f_0)(\lambda)
+ \frac{1}{\lambda}  \cP  (\mO_2 f_2+ \mO_4 f_0)(\lambda)\Big]d\lambda \\
\label{bk2.32}
    \begin{split}&=\cL^{-1}\cP^\bot \mO_{2}\cL^{-1}\cP^\bot \mO_{2} \cP
- \cL^{-1}\cP^\bot \mO_{4} \cP\\
&\quad+ \cP\mO_{2}\cL^{-1}\cP^\bot \mO_{2}\cL^{-1}\cP^\bot
- \cP \mO_{4}\cL^{-1}\cP^\bot \\
&\quad+ \cL^{-1}\cP^\bot \mO_{2} \cP \mO_{2} \cL^{-1}\cP^\bot
-  \cP \mO_{2} \cL^{-2} \cP^\bot  \mO_{2} \cP\\
&\quad- \cP\mO_{2} \cP \mO_{2} \cL^{-2}\cP^\bot
- \cP^\bot \cL^{-2}\mO_{2}  \cP \mO_{2} \cP.
\end{split}\end{align}
In particular, the first two identities of \eqref{bk2.31} imply
\begin{align}\label{bk2.33}
J_{0,\,x_{0}}= 1, \quad J_{1,\,x_{0}}=0.
\end{align}

\begin{rem} \label{toet2.7}
$\cL_t$ is a formally self-adjoint elliptic operator on
$\cC^\infty(\R^{2n}, E_{x_0})$ with respect to the norm
$\|\cdot\|_{L^2}$
induced by $h^{E_{x_0}}$, $dZ$. Thus $\cL_0$
and $\mO_r$ are also formally self-adjoint with respect to
$\|\cdot\|_{L^2}$.
Therefore the third and fourth terms in \eqref{bk2.32} are the adjoints 
of the first and second terms, respectively.
In Lemma \ref{bkt3.1}, we will show that $\cP\mO_2\cP=0$, hence the
last two terms in \eqref{bk2.32} vanish.
Set
\begin{align}\label{bk2.34}
\cF_{41}= \cL^{-1}\cP^\bot \mO_{2}\cL^{-1}\cP^\bot \mO_{2} \cP
- \cL^{-1}\cP^\bot \mO_{4} \cP.
\end{align}
\end{rem}

%\newpage
%%%%%%%%%%%%%%%%%%%%%%%%%%%%%%%%%%%%%%%%%
%%%%%%%%%%%%%%%%%%%%%%%%%%%%%%%%%%%%%%%%%
\section{Taylor expansion of the rescaled operator $\cL_t$}\label{toes3}
%%%%%%%%%%%%%%%%%%%%%%%%%%%%%%%%%%%%%

In this section we compute
the operators $\cL_{0}$ and $\mO_{i}$ (for $1\leqslant i \leqslant 4$) 
from \eqref{bk2.22} (see Theorem \ref{bkt2.21}),  which will
be used in Sections \ref{bks3}, \ref{toes4} for the evaluation of the
coefficients of the expansion of the kernels of the Berezin-Toeplitz
operators.

 We denote by $\left\langle\cdot,\cdot\right\rangle$ the
$\C$-bilinear form on $TX\otimes_\R \C$ induced by $g^{TX}$.
 Let $R^{TX}$ be the curvature of the Levi-Civita connection $\nabla^{TX}$.
Let $\ric$ and $\br$ be the Ricci and scalar curvature of $\nabla^{TX}$.
Then we have the following well know facts:
for $U,V, W,Y$ vector fields on $X$, %(cf. \cite[\S 1.2]{BeGeVe})
\begin{align}\label{alm01.0}
\begin{split}
&R^{TX}(U,V)W+ R^{TX}(V,W)U + R^{TX}(W,U)V=0,\\
&\langle R^{TX}(U,V)W,Y\rangle= \langle R^{TX}(W,Y)U,V\rangle
 = - \langle R^{TX}(V,U)W,Y\rangle.
\end{split}
\end{align}

Now we work on $T_{x_0}X\simeq \R^{2n}$ as in \eqref{n11}.
Recall that we have trivialized $L, E$.
Let $\nabla_U$ denote the ordinary differentiation
operator on $T_{x_0}X$ in the direction $U$.

We adopt the convention that all tensors will be evaluated at the base
point $x_0\in X$ and most of time, we will omit the subscript $x_0$.

For $W\in T_{x_0}X$, $Z\in\R^{2n}$,
  let $\wi{W}(Z)$ be the parallel transport of $W$
with respect to $\nabla^{TX}$ along the curve $[0,1]\ni u\to uZ$.
Because the complex structure $J$ is parallel with respect to
$\nabla^{TX}$, we know that
\begin{align}\label{lm01.1}
J_Z \wi{W}(Z) = \wi{J_{x_0}W}(Z).
\end{align}

Recall that $\{e_i\}$ is a fixed orthonormal basis of
$(T_{x_0}X, g^{TX})$. Then for $U,V\in T_{x_0}X$,
\begin{align}\label{alm01.1}
\ric_{x_0}(U,V)=- \langle  R^{TX}(U,e_j)V,e_j\rangle_{x_0},\quad
\br_{x_0}= - \langle  R^{TX}(e_i,e_j)e_i,e_j\rangle_{x_0}.
\end{align}
We define 
\begin{align*}
&R^{TX}_{\, ;\, \bullet}\in (T^*X\otimes \Lambda^2(T^*X) \otimes
\End(TX))_{x_0}\,,\\
&R^{TX}_{\, ; \,(\bullet,\bullet)}\in
((T^*X)^{\otimes 2}\otimes  \Lambda^2(T^*X) \otimes \End(TX))_{x_0}\,,\\
&\ric_{\, ; \,\bullet}\in (T^*X\otimes (T^*X)^{\otimes 2})_{x_0}\,,\\
&R^E_{\, ;\, \bullet}\in (T^*X\otimes\Lambda^2(T^*X)\otimes \End(E))_{x_0}\,,\\
&R^{E}_{\, ; \,(\bullet,\bullet)}\in ((T^*X)^{\otimes 2}\otimes
 \Lambda^2(T^*X) \otimes \End(E))_{x_0}\,,
\end{align*}
by
\begin{align} \label{lm01.2}
 \begin{split}
&\left \langle R^{TX}_{\, ;\, e_{k}}(e_{m}, e_j) e_{q}, e_i\right\rangle
= \Big(\nabla_{e_k}\left \langle  R^{TX} (\wi{e}_{m}, \wi{e}_j) \wi{e}_{q},
\wi{e}_i\right\rangle\Big)_{x_0},\\
&\left \langle R^{TX}_{\, ; (e_k, \,e_\ell)}
(e_{m}, e_j) e_{q}, e_i\right\rangle
= \Big(\nabla_{e_\ell}\nabla_{e_k}\left \langle
R^{TX} (\wi{e}_{m}, \wi{e}_j) \wi{e}_{q}, \wi{e}_i\right\rangle\Big)_{x_0},\\
&\ric_{\, ; \,e_{k}}(e_i,e_j)= (\nabla_{e_k} \ric(\wi{e}_i,\wi{e}_j))_{x_0},\\
&R^E_{\, ; \,e_{k}}(e_i,e_j)= (\nabla_{e_k} R^E(\wi{e}_i,\wi{e}_j))_{x_0},\\ 
&R^E_{\, ;\, (e_{k},\, e_{\ell})}(e_i,e_j)
= (\nabla_{e_\ell}\nabla_{e_k} R^E(\wi{e}_i,\wi{e}_j))_{x_0}.
\end{split}
\end{align}
We will also use the complex coordinates $z=(z_1,\ldots,z_n)$. Note that
\begin{equation}\label{lm01.3}
e_{2j-1}=\frac{\partial}{\partial Z_{2j-1}}
=\frac{\partial}{\partial z_j}+ \frac{\partial}{\partial \ov{z}_j},
\quad   e_{2j}=\frac{\partial}{\partial Z_{2j}}
=\sqrt{-1}\Big(\frac{\partial}{\partial z_j}
- \frac{\partial}{\partial \ov{z}_j}\Big),
 \quad \left|\frac{\partial}{\partial z_j}\right|^2=\frac{1}{2}.
\end{equation}
Set
%\begin{equation} \label{lm01.4} \begin{split}
%&R_{k\ov{m}\ell\ov{q}} =\left \langle  R^{TX} (\tfrac{\partial}{\partial
%  z_k},
%\tfrac{\partial}{\partial \ov{z}_m}) \tfrac{\partial}{\partial z_\ell},
%\tfrac{\partial}{\partial  \ov{z}_q}\right\rangle_{x_0},
%\quad \ric_{k\ov{\ell}}= \ric_{x_0}(\tfrac{\partial}{\partial z_k},
%\tfrac{\partial}{\partial \ov{z}_\ell}), \quad
%R^E_{k\ov{\ell}}= R^E_{x_0}(\tfrac{\partial}{\partial z_k},
%\tfrac{\partial}{\partial \ov{z}_\ell}),\\
%&R_{k\ov{m}\ell\ov{q};\,s} =\Big \langle  R^{TX}_{\, ;
%\tfrac{\partial}{\partial z_{s}}}(\tfrac{\partial}{\partial z_k},
%\tfrac{\partial}{\partial\ov{z}_m}) \tfrac{\partial}{\partial z_\ell},
%\tfrac{\partial}{\partial\ov{z}_q}\Big\rangle,\quad
%R^E_{k\ov{q}; \,s}=
%R^E_{\, ; \tfrac{\partial}{\partial z_{s}}}(\tfrac{\partial}{\partial z_k},
%\tfrac{\partial}{\partial \ov{z}_q}),
%\end{split}\end{equation}
\begin{equation} \label{lm01.4} \begin{split}
&R_{k\ov{m}\ell\ov{q}} =\left \langle  R^{TX} \Big(\frac{\partial}{\partial z_k},
\frac{\partial}{\partial \ov{z}_m}\Big) \frac{\partial}{\partial z_\ell},
\frac{\partial}{\partial  \ov{z}_q}\right\rangle_{x_0},\quad
R^E_{k\ov{\ell}}= R^E_{x_0}\Big(\frac{\partial}{\partial z_k},
\frac{\partial}{\partial \ov{z}_\ell}\Big),\\
&\ric_{k\ov{\ell}}= \ric_{x_0}\Big(\frac{\partial}{\partial z_k},
\frac{\partial}{\partial \ov{z}_\ell}\Big)\,,\\ 
&R_{k\ov{m}\ell\ov{q};\,s} =\left \langle  R^{TX}_{\, ;
\frac{\partial}{\partial z_{s}}}\Big(\frac{\partial}{\partial z_k},
\frac{\partial}{\partial\ov{z}_m}\Big) \frac{\partial}{\partial z_\ell},
\frac{\partial}{\partial\ov{z}_q}\right\rangle\,,\quad
R^E_{k\ov{q}; \,s}=
R^E_{\, ; \frac{\partial}{\partial z_{s}}}\Big(\frac{\partial}{\partial z_k},
\frac{\partial}{\partial \ov{z}_q}\Big),
\end{split}\end{equation}
and in the same way, we define $R_{k\ov{m}\ell\ov{q};\, \ov{s}}$\,, 
$R_{k\ov{m}\ell\ov{q};\, t\ov{s}}$\,,
$\ric_{k\ov{q};\, \ov{s}}$\,, $R^E_{k\ov{q};\, \ov{s}}$, 
$R^E_{k\ov{q};\, t\ov{s}}$\,.

Since $R^{TX}$ is a $(1,1)$-form
and $\nabla^E$ is the Chern connection on $(E,h^E)$, we deduce from
\eqref{alm01.0}--\eqref{alm01.1} the following.
%Then by \eqref{alm01.0} and \eqref{lm01.1}, and using
%the fact that $R^{TX}$ is a $(1,1)$-form
%and $\nabla^E$ is the Chern connection on $(E,h^E)$, we get
\begin{lemma}\label{lmt1.6}
\noindent
\\[2pt]
{\rm(1)} $R_{k\ov{m}\ell\ov{q}} =R_{\ell\ov{m}k\ov{q}}=
R_{k\ov{q}\ell\ov{m}}= R_{\ell\ov{q}k\ov{m}}$\,, 
$\br= 8\,  R_{m\ov{q}q\ov{m}}$\,, $(R^E_{k\ov{q}})^* 
=R^E_{q\ov{k}}$\,, \\[2pt]
{\rm(2)} $R_{k\ov{m}\ell\ov{q}\,;\, \ov{s}} =R_{\ell\ov{m}k\ov{q}\,;\, \ov{s}}=
R_{k\ov{q}\ell\ov{m}\,; \,\ov{s}}= R_{\ell\ov{q}k\ov{m}\,; \,\ov{s}}$\,,
\\[2pt]
{\rm(3)} $\ric$ is a symmetric $(1,1)$-tensor and 
$\ric_{m\ov{q}}=2 \,  R_{m\ov{k}k\ov{q}}$\,,
$\ric_{m\ov{q}\,; \,\ov{s}}=2 \,  R_{m\ov{k}k\ov{q}\,;\, \ov{s}}$\,,
\\[2pt]
{\rm(4)} $R^{TX}_{\, ;\, e_k}$, $R^{TX}_{\, ;\, (e_k,\,e_\ell)}$
 are
$(1,1)$-forms  with values in $\End(T_{x_0}X)$ which commute with
$J_{x_0}$\,,
\\[2pt]
{\rm(5)} 
$R^{E}_{\, ; \,e_k}$\,, $R^{E}_{\, ; \,(e_k,\,e_\ell)}\in \End(E_{x_0})$.
%are $(1,1)$-forms  with values in $\End(E_{x_0})$.
\end{lemma}
%\begin{align} \label{lm01.6}\begin{split}
%&R_{k\ov{m}\ell\ov{q}} =R_{\ell\ov{m}k\ov{q}}=
%R_{k\ov{q}\ell\ov{m}}= R_{\ell\ov{q}k\ov{m}},\quad
%\br= 8\,  R_{m\ov{q}q\ov{m}} \, ,\quad (R^E_{k\ov{q}})^* =R^E_{q\ov{k}}, \\
%&R_{k\ov{m}\ell\ov{q}; \ov{s}} =R_{\ell\ov{m}k\ov{q}; \ov{s}}=
%R_{k\ov{q}\ell\ov{m}; \ov{s}}= R_{\ell\ov{q}k\ov{m}; \ov{s}},\\
%&\ric \text{  is a symmetric } (1,1)-\text{tensor and }
%\, \,  \ric_{m\ov{q}}=2 \,  R_{m\ov{k}k\ov{q}}\, ,
%\quad  \ric_{m\ov{q}; \ov{s}}=2 \,  R_{m\ov{k}k\ov{q}; \ov{s}}\, , \\
%&R^{TX}_{\, ; e_k}, R^{TX}_{\, ; (e_k,e_\ell)}
%\text{  (resp. $R^{E}_{\, ; e_k}, R^{E}_{\, ; (e_k,e_\ell)}$) are
%$(1,1)$-form  with values} \\
%&\text{ in $\End(T_{x_0}X)$ which commute with
%$J_{x_0}$
%(resp. in $\End(E_{x_0})$)}.
%\end{split}\end{align}
Let $\text{div}(\ric)$ be the divergence of $\ric$.
By \cite[\S 2.3.4, Prop. 6]{Petersen06},
\begin{align}\label{alm01.4}
d\br = 2\, \text{div} (\ric) = 2 \,
(\nabla^{T^{*}X}_{\wi{e}_{m}}\ric)(\wi{e}_{m}, \cdot)\,.
\end{align}
Lemma \ref{lmt1.6} and \eqref{alm01.4} entail
\begin{align}\label{alm01.5}\begin{split}
&R_{\ell\ov{\ell}m\ov{m}\,;\, \ov{k}}
=  R_{\ell\ov{\ell}m\ov{k}\,;\, \ov{m}}\,,\quad
R_{\ell\ov{\ell}m\ov{m}\,;\, k}= R_{\ell\ov{\ell}k\ov{m}\,;\, m},\\
-(&\Delta\br)_{\,x_{0}} =2 e_qe_m(\ric(\wi{e}_q,\wi{e}_m))_{\,x_{0}}
= 32  R_{k\ov{m} q\ov{q}\,;\,m\ov{k}}
= 32  R_{m\ov{m} q\ov{q}\,;\, k\ov{k}}\;\;.
\end{split}\end{align}
Set
\begin{equation}\label{bk2.62}
\mR:= \sum_i Z_i e_i =Z, \qquad
\nabla_{0,\scriptstyle\bullet}:= \nabla_{\scriptstyle\bullet} 
+ \frac{1}{2} R^L_{x_0}(\mR, \,{\scriptstyle\bullet}\,).
\end{equation}
Thus $\mR$ is the radial vector field on
$\R^{2n}$. We also introduce the vector fields
$z= \sum_{i} z_{i} \tfrac{\partial}{\partial z_i}$ and
$\ov{z}= \sum_{i} \ov{z}_{i} \tfrac{\partial}{\partial \ov{z}_i}$.
By \cite[Prop. 1.2.2]{MM07}, \eqref{lm01.1}, we have
\begin{align} \label{lm01.5}
\mR= \sum_i Z_i \wi{e}_i,\quad
z= \sum_{i} z_{i} \wi{\tfrac{\partial}{\partial z_i}},\quad
\ov{z}= \sum_{i} \ov{z}_{i} \wi{\tfrac{\partial}{\partial \ov{z}_i}}.
\end{align}
By \eqref{toe2.1} and \eqref{toe1.1}, we get
\begin{equation}\label{bk2.64}
 b_i=-2\nabla_{0,\tfrac{\partial}{\partial z_i}},\quad
b^{+}_i=2\nabla_{0,\tfrac{\partial}{\partial \overline{z}_i}},\quad
 R^L_{x_0} = -2\pi\sqrt{-1}\left \langle J\,\scriptstyle\bullet\,,\,
 \scriptstyle\bullet\, \right \rangle_{x_0}.
\end{equation}
Let $A_{1Z}, A_{2Z}\in (T^*X)^{\otimes 2}$ be polynomials in $Z$ with
values symmetric tensors, defined by
\begin{align} \label{lm01.7} \begin{split}
&A_{1Z}(e_i,e_j)=\left\langle R^{TX}_{\, ;(Z,Z)} (\mR,e_i) \mR,
e_j\right\rangle_{x_0},\\
&A_{2Z}(e_i,e_j)=\left \langle R^{TX}_{x_0} (\mR,e_i) \mR,
R^{TX}_{x_0} (\mR,e_j) \mR\right \rangle_{x_0}.
\end{split}\end{align}
\comment{
Then by (\ref{alm01.0}), $R^{TX}$ is a $(1,1)$-form
and $\mR= z + \ov{z}$ as vector fields, we get
\begin{align} \label{lm01.8} \begin{split}
A_{1Z}(e_j,e_j)=& 4A_{1Z}(\tfrac{\partial}{\partial z_j},
\tfrac{\partial}{\partial \ov{z}_j})= 4 \left\langle R^{TX}_{\, ;(Z,Z)}
(z,\tfrac{\partial}{\partial \ov{z}_j})\ov{z} ,
\tfrac{\partial}{\partial z_j}\right\rangle_{x_0},\\
A_{2Z}(e_j,e_j)=& 4A_{2Z}(\tfrac{\partial}{\partial z_j},
\tfrac{\partial}{\partial \ov{z}_j})=4 \left \langle R^{TX}_{x_0}
(\ov{z}, \tfrac{\partial}{\partial z_j}) \mR,
R^{TX}_{x_0} (z,\tfrac{\partial}{\partial \ov{z}_j}) \mR\right \rangle_{x_0}.
\end{split}\end{align}
}
Recall that the operator $\cL$ was defined in \eqref{toe1.1}. Set
\begin{align} \label{lm01.9a}
 \begin{split}
    \mO^{\,\prime}_2=&\frac{1}{3} \left \langle R^{TX}_{x_0} (\mR,e_i)
    \mR, e_j\right \rangle
     \nabla_{0,e_i} \nabla_{0,e_j}
    - 2 R^E_{k\ov{k}} \\
&\hspace{10mm}    +  \left(\left \langle
   \frac{\pi}{3} R^{TX}_{x_0} (z,\ov{z}) \mR,  e_j\right \rangle
  +  \frac{2}{3}  \ric_{x_0} (\mR,e_j)  -R^E_{x_0} (\mR, e_j)\right)
    \nabla_{0,e_j},
\end{split}\end{align}
and
\begin{align} \label{lm01.9}
 \begin{split}
\mO_{41}= &  \frac{1}{20}\Big( A_{1Z}
-  \frac{4}{3}A_{2Z}\Big) (e_{i}, e_{j})
\nabla_{0,e_i}\nabla_{0,e_j},\\
\mO_{42}=&\Big[\cL, -\Big(\frac{1}{80}A_{1Z}
 -  \frac{1}{360}A_{2Z}\Big)(e_{j}, e_{j})
- \frac{1}{288} \ric(\mR,\mR)^{2} \Big]
+  \frac{\cL}{144} \ric(\mR,\mR)^{2},\\
\mO_{43}=&-\frac{1}{144} \ric(\mR,\mR) \cL  \ric(\mR,\mR),\\
\mO_{44}=&\Big \{ \frac{\pi}{30} A_{1Z}(\ov{z}, e_{i})
%\left \langle R^{TX}_{\, ;(Z,Z)} (z,\ov{z})\mR,
%e_i\right \rangle_{x_0}
- \frac{\pi}{10}  A_{2Z}(\ov{z}, e_{i})
%\left \langle R^{TX}_{x_0} (z,\ov{z})\mR,
%R^{TX}_{x_0} (\mR, e_i)\mR\right \rangle_{x_0}\\
+  \tfrac{\partial}{\partial Z_j} \Big(\frac{1}{20}A_{1Z}
+  \frac{2}{45}A_{2Z} \Big)(e_{i}, e_{j}) \\
&
- \tfrac{\partial}{\partial Z_i}  \Big(\frac{1}{40}A_{1Z}
 +  \frac{1}{45}A_{2Z} \Big)(e_{j}, e_{j})\Big\}
\nabla_{0,e_i},
\end{split}\end{align}
and
\begin{align} \label{lm01.9b}
 \begin{split}
\mO_{45}=&\Big \{\frac{2}{9}\left \langle R^{TX}_{x_0} (\mR,e_k) \mR,
R^{TX}_{x_0} (\mR,e_k)e_\ell\right\rangle_{x_0}
- \frac{1}{9} \left \langle R^{TX}_{x_0} (\mR,e_\ell) \mR,
e_k\right\rangle_{x_0} \ric(\mR,e_k)\\
&+\frac{1}{4}
\left \langle R^{TX}_{x_0} (\mR,e_\ell) \mR, e_m\right \rangle_{x_0}
R^E_{x_0} (\mR,e_m)
-  \frac{1}{4} R^E_{\, ;(Z,Z)} (\mR,  e_\ell)\Big\}
\nabla_{0,e_\ell},\\
\mO_{46}=&-\frac{\pi^2}{36} A_{2Z}(\ov{z}, \ov{z})
%\left\langle  R^{TX}_{x_0}(z,\ov{z})\mR,
% R^{TX}_{x_0}(z,\ov{z})\mR\right\rangle
+ \frac{\pi}{30} \left \langle
   R^{TX}_{\, ;(Z,e_\ell)} (z,\ov{z})\mR, e_\ell\right \rangle_{x_0}\\
&\hspace{-3mm}
%- \frac{11\pi}{180} \left \langle R^{TX}_{x_0} (z,\ov{z})e_l,
%R^{TX}_{x_0} (\mR, e_l)\mR\right \rangle_{x_0}
- \frac{\pi}{20}\left \langle R^{TX}_{x_0} (z,\ov{z})\mR,
e_m \right \rangle_{x_0}  \ric(\mR,e_m)
 + \frac{4}{9} \ric_{k\ov{m}}\ric_{m\ov{\ell}}z_{k}\ov{z}_{\ell}
- \frac{4}{9} R_{k\ov{\ell}m\ov{q}}\ric_{\ell\ov{m}}z_{k}\ov{z}_{q}\\
%+\frac{1}{9} \ric(\mR,e_m)\ric(\mR,e_m)
%+\frac{1}{18}\left \langle R^{TX}_{x_0} (\mR,e_\ell) \mR,
%e_m\right\rangle_{x_0} \ric(e_\ell,e_m) \\
&\hspace{-3mm}+ \frac{1}{6}\left\langle 
\pi R^{TX}_{x_0}(z,\ov{z})\mR , e_m\right\rangle
 R^E_{x_0} (\mR,e_m)
+ \frac{1}{8} \ric(\mR,e_m)R^E_{x_0} (\mR,e_m)\\
&\hspace{-3mm}+ \frac{1}{2} (R^E_{k\ov{m}}R^E_{m\ov{\ell}}
+ R^E_{m\ov{\ell}}R^E_{k\ov{m}}) z_{k}\ov{z}_{\ell}
%-  \frac{1}{4}R^E_{x_0} (\mR,e_m)^2\\
 - \frac{1}{4}R^E_{\, ;(Z,e_\ell)} (\mR,  e_\ell)
- R^E_{\, ;(Z,Z)}
(\tfrac{\partial}{\partial z_\ell},\tfrac{\partial}{\partial \ov{z}_\ell}).
% - \frac{4}{3} R_{k\ov{\ell}l\ov{m}} z_k \ov{z}_m R^E_{i\ov{l}}.
\end{split}\end{align}
The following result extends \cite[Th.\,4.1.25]{MM07} where
$\cL_0, \mO_1,\mO_2$ were computed.

\begin{thm} \label{bkt2.21} The following identities hold for the operators 
 $\mO_r$ introduced in \eqref{bk2.22}\,\rm{:}
\begin{subequations} 
 \begin{align}
\cL_0=& \sum_j b_jb^+_j=\cL
=-\sum_{i} \nabla_{0,e_{i}}\nabla_{0,e_{i}}- 2\pi n, \qquad \mO_1=0,
\label{lm01.101}\\
\mO_2=& \mO^{\,\prime}_2 - \frac{1}{3} \ric_{x_0}(\mR,e_j)\nabla_{0,e_j}
- \frac{\br_{x_0}}{6},\label{lm01.102}
\end{align}
\end{subequations}
and
\begin{subequations}
\begin{align}\label{lm01.11}
\mO_3=& \frac{1}{6}\left \langle R^{TX}_{\, ; Z}  (\mR,e_i) \mR, e_j\right
\rangle_{x_0} \nabla_{0,e_i} \nabla_{0,e_j}\\
&+  \Big[\frac{2\pi}{15}\left \langle
   R^{TX}_{\, ; Z}  (z,\ov{z})\mR, e_i\right\rangle_{x_0}
+\frac{1}{6} \ric_{\, ; Z}  (\mR,e_i)\nonumber\\
%+ \frac{1}{12}\nabla_{e_i} \ric (\mR,\mR)\\
%\left\langle  R^{TX}_{\, ; Z} (\mR, e_j) e_j, e_i\right\rangle_{x_0}\\
&\hspace{10mm}+ \frac{1}{6}
\left\langle  R^{TX}_{\, ; e_j} (\mR, e_j) \mR, e_i\right\rangle_{x_0}
%\left\langle \nabla_{e_i} R^{TX}(\mR, e_j)e_j,\mR \right\rangle_{x_0}
- \frac{2}{3} R^E_{\, ; Z}  (\mR, e_i)\Big]\nabla_{0,e_i} \nonumber\\
&+  \frac{\pi}{15}\left \langle
  R^{TX}_{\, ; e_{j}} (z,\ov{z})\mR, e_j\right \rangle_{x_0}
-\frac{1}{6} \ric_{\, ; e_{i}} (\mR, e_i)\nonumber\\
& -\frac{1}{12} \ric_{\, ; Z} (e_i, e_i)
%- \frac{1}{24} \frac{\partial^2}{\partial Z_i^2} \ric_{\, ; Z}  (\mR,\mR)
% \left\langle  R^{TX}_{\, ; Z}  (\mR, e_j) \mR, e_j\right\rangle_{x_0}
-\frac{1}{3}  R^E_{\, ; e_{i}} (\mR, e_i)- \frac{\sqrt{-1}}{2}
R^E_{\, ; Z} (e_i, Je_i),\nonumber\\
\mO_4=&\mO_{41}+ \mO_{42}+\mO_{43}+\mO_{44}
+\mO_{45}+\mO_{46}.\label{lm01.12}
\end{align}
\end{subequations}
\end{thm}

\begin{proof}
%Let $\{e_i\}_i$ be an oriented orthonormal basis of $T_{x_0}X$.
%We also denote by $\{e ^i\}_i$ the dual basis of $\{e_i\}$.
 Recall that $\wi{e}_i (Z)$ is the parallel
transport of ${e}_i$ with respect to $\nabla^{TX}$
 along the curve $[0,1]\ni u \to uZ$.
%Note that $e_j=\frac{\partial}{\partial Z_j}$.
Let $\wi{\theta} (Z) = (\theta _j^i (Z))_{i,j=1}^{2n}$
be the $2 n \times 2n$-matrix such that
\be\label{lm01.29}
e_i = \sum_j \theta ^j_i(Z) \wi{e}_j (Z), \quad
\wi{e}_j (Z)= (\wi{\theta} (Z)^ {-1})_j^k e_k.
\ee
Taking into account the Taylor expansion
of $\theta ^i_j$ at $0$ we have (cf.\,\cite[(1.2.27)]{MM07}) 
\begin{align}\label{lm01.35}
\sum_{|\alpha| \geqslant 1} ( |\alpha|^2 + |\alpha|)
(\partial ^\alpha\theta ^i_j)(0) \frac{Z^\alpha}{\alpha !}=
\left \langle R^{TX} ( \mR,e_j) \mR, \wi{e}_i\right \rangle_Z .
\end{align}
From this equation, we obtain first that
\be\label{lm01.27}
e_j(Z)  = \wi{e}_j(Z) + \frac{1}{6}
\left \langle R^{TX}_{x_0} (\mR,e_j) \mR, e_k\right \rangle_{x_0}
\wi{e}_k(Z) + \cO(|Z|^3).
\ee
From \eqref{lm01.2}, \eqref{lm01.5}, \eqref{lm01.7},
\eqref{lm01.35} and \eqref{lm01.27}, we get further
\begin{equation}\label{lm01.36}
\begin{split}
\theta ^i_j  = \delta_{ij}
& + \frac{1}{6}
\left \langle R^{TX}_{x_0} (\mR,e_i) \mR, e_j\right \rangle_{x_0}
+  \frac{1}{12}
\left \langle R^{TX}_{\, ; Z}  (\mR,e_i) \mR, e_j\right
\rangle_{x_0}\\
& +  \frac{1}{20} \Big(\frac{1}{2}A_{1Z}(e_i,e_j)
%\left \langle R^{TX}_{\, ;(Z,Z)} (\mR,e_i) \mR,
%e_j\right\rangle_{x_0}
+  \frac{1}{6}A_{2Z}(e_i,e_j)
%\left \langle R^{TX}_{x_0} (\mR,e_i) \mR,
%R^{TX}_{x_0} (\mR,e_j) \mR\right \rangle_{x_0}
\Big) + \cO(|Z|^5)\,.
\end{split}
\end{equation}
Set $g_{ij}(Z)= g^{TX}(e_i,e_j)(Z) =\langle e_i,e_j\rangle_Z$
and let $(g^{ij}(Z))$ be the inverse of the matrix $(g_{ij}(Z))$.
Then by \eqref{lm01.7}, \eqref{lm01.29} and \eqref{lm01.36}, we have
\begin{equation}\label{lm01.37}
\begin{split}
g_{ij}(Z)= \theta^k_i(Z) \theta^k_j(Z)
= \delta_{ij} &+  \frac{1}{3}
\left \langle R^{TX}_{x_0} (\mR,e_i) \mR, e_j\right \rangle_{x_0}
 + \frac{1}{6}\left \langle R^{TX}_{\, ; Z}  (\mR,e_i) \mR, e_j\right
\rangle_{x_0}\\
& +  \frac{1}{20} A_{1Z}(e_i,e_j)
+  \frac{2}{45}A_{2Z}(e_i,e_j)
+ \cO (|Z|^5).
\end{split}
\end{equation}
In view of the expansion $(1+a)^{-1}= 1-a +a^2+\ldots$\,, we obtain
\begin{equation}\label{lm01.38}
\begin{split}
g^{ij}(Z)
= \delta_{ij} & - \frac{1}{3}
\left \langle R^{TX}_{x_0} (\mR,e_i) \mR, e_j\right \rangle_{x_0}
 - \frac{1}{6}\left \langle R^{TX}_{\, ; Z}  (\mR,e_i) \mR, e_j\right
\rangle_{x_0}\\
&-  \frac{1}{20} A_{1Z}(e_i,e_j)
 +  \frac{1}{15}A_{2Z}(e_i,e_j)
+ \cO (|Z|^5).
\end{split}
\end{equation}
If $\Gamma _{ij}^\ell$ are the Christoffel symbols
% connection form
of $\nabla ^{TX}$ with respect to the frame $\{e_i\}$, then
$(\nabla ^{TX}_{e_i}e_j)(Z)$ $= \Gamma _{ij}^\ell (Z) e_\ell$.
By %(\ref{lm01.1}) and (\ref{alm01.27}),
the explicit formula for $\nabla^{TX}$, we get 
(cf.\,\cite[(4.1.102)]{MM07}) with $\partial_j:=\frac{\partial}{\partial Z_j}$
\begin{equation}\label{bk2.83}
\begin{split}
\Gamma _{ij}^\ell (Z)& =  \frac{1}{2}  g^{\ell k} (\partial_i g_{jk}
+ \partial_j g_{ik}-\partial_k g_{ij})(Z)\\
&= \frac{1}{3}\Big  [ \left \langle R^{TX}_{x_0} 
(\mR, e_j) e_i, e_\ell\right \rangle _{x_0}
+ \left \langle R^{TX}_{x_0} (\mR, e_i) e_j, e_\ell\right \rangle_{x_0}\Big ]
 + \cO(|Z|^2)\,.
 \end{split}
\end{equation}
For $j$ fixed, $\Gamma _{jj}^\ell (Z)=\frac{1}{2}  g^{\ell k} ( 2\partial_j g_{jk}
-\partial_k g_{jj})(Z)$, thus by \eqref{lm01.37} and \eqref{lm01.38},
\begin{equation}\label{abk2.83}
\begin{split}
\Gamma&_{jj}^\ell (Z) = \frac{2}{3}
\left \langle R^{TX}_{x_0} (\mR, e_j) e_j, e_\ell\right \rangle _{x_0}\\
&+  \frac{1}{12} \Big[
4 \left \langle  R^{TX}_{\, ; Z}  (\mR, e_j) e_j, e_\ell\right
\rangle_{x_0} 
+ 2\left \langle  R^{TX}_{\, ; e_j} (\mR, e_j) \mR, e_\ell\right
\rangle_{x_0}
+ \left \langle R^{TX}_{\, ; e_\ell} (\mR, e_j) e_j,\mR \right
\rangle_{x_0} \Big ] \\
&- \frac{2}{9} \left \langle R^{TX}_{x_0} (\mR,e_\ell) \mR, R^{TX}_{x_0}
  (\mR,e_j) e_j\right\rangle_{x_0}
+  \frac{\partial}{\partial Z_j} \Big(
\frac{1}{20}A_{1Z}
 +  \frac{2}{45}A_{2Z}\Big)(e_j,e_\ell) \\
&-  \frac{1}{2} \frac{\partial}{\partial Z_l}\Big (
\frac{1}{20}A_{1Z}
+  \frac{2}{45}A_{2Z}\Big)(e_j,e_j)
 + \cO(|Z|^4).
 \end{split}
\end{equation}
Note that 
\[
\det(\delta_{ij}+ a_{ij})= 1 + \sum_{i} a_{ii} +
\sum_{i<j}(a_{ii}a_{jj}- a_{ij}a_{ji}) + \ldots
\] 
and
\[(1+a)^{1/4}= 1+ \frac{1}{4} a - \frac{3}{32} a^{2} + \ldots
\]
By \eqref{alm01.1} and \eqref{lm01.37}, we get
\begin{equation}\label{bk2.82}
\begin{split}
\kappa(&Z)^{1/2}= |\det (g_{ij}(Z))|^{1/4}  = 1 + \frac{1}{12}
\left \langle R^{TX}_{x_0} (\mR,e_j) \mR, e_j\right \rangle_{x_0}\\ 
&-  \frac{1}{24}   \ric_{\, ; Z}(\mR,\mR)
 + \frac{1}{80} A_{1Z}(e_j,e_j) +  \frac{1}{90} A_{2Z}(e_j,e_j)\\
&+ \frac{1}{36} \sum_{i<j} \Big( \left \langle R^{TX}_{x_0} (\mR,e_i)
  \mR, e_i\right \rangle_{x_0} \left \langle R^{TX}_{x_0} (\mR,e_j)
  \mR, e_j\right \rangle_{x_0}- \left \langle R^{TX}_{x_0} (\mR,e_i) \mR,
  e_j\right \rangle_{x_0} ^2 \Big)\\
&-  \frac{1}{96} \Big(\sum_j\left \langle R^{TX}_{x_0} (\mR,e_j) \mR,
  e_j\right \rangle_{x_0}\Big) ^2+ \cO (|Z|^5)\\
=&\, 1-\frac{1}{12} \ric(\mR,\mR) -\frac{1}{24}   \ric_{\, ; Z}(\mR,\mR)
+ \Big(\frac{1}{80} A_{1Z}
 -  \frac{1}{360} A_{2Z}\Big)(e_j,e_j)\\
 &+  \frac{1}{288} \ric(\mR,\mR)^2 + \cO (|Z|^5).
 \end{split}
\end{equation}
Thus
\begin{equation}\label{abk2.82}
\begin{split}
\kappa(Z)^{-1/2}
= 1&+\frac{1}{12} \ric(\mR,\mR) +\frac{1}{24}   \ric_{\, ; Z}(\mR,\mR)
- \Big(\frac{1}{80} A_{1Z}
 -  \frac{1}{360} A_{2Z}\Big)(e_j,e_j) \\
&+ \Big(\frac{1}{144} -\frac{1}{288}\Big) \ric(\mR,\mR)^2 + \cO (|Z|^5)\,.
\end{split}
\end{equation}

%\newpage

Observe that $J$ is parallel with respect to $\nabla ^{TX}$,
thus $\left \langle J\wi{e}_i,\wi{e}_j \right \rangle_{Z}=\left
\langle Je_i,e_j \right \rangle_{x_0}$.
% By \cite[Prop. 1.2.2]{MM07}, $\mR=Z_i \wi{e}_i$, thus
From (\ref{bk2.62}), (\ref{lm01.5}),
 (\ref{lm01.29}) and (\ref{lm01.36}), we get
\begin{equation}\label{bk2.84}
\begin{split}
\frac{\sqrt{-1}}{2 \pi} R^L_Z(\mR,e_\ell)=&\theta ^j_\ell(Z)
\left \langle J\wi{e}_i,\wi{e}_j \right \rangle_{Z} Z_i
=\theta ^j_\ell(Z)
\left \langle J\mR,e_j \right \rangle_{x_0}  \\
= &\left \langle J\mR, e_\ell \right \rangle_{x_0}
+\frac{1}{6}\left \langle R^{TX}_{x_0}(\mR, J\mR)\mR, 
e_\ell\right \rangle_{x_0}\\
 +&\frac{1}{12}\left \langle  R^{TX}_{\, ; Z}  
 (\mR, J\mR)\mR, e_\ell\right \rangle_{x_0}
 + \frac{1}{40}\left \langle R^{TX}_{\, ;(Z,Z)} (\mR,
   J\mR)\mR, e_\ell\right \rangle_{x_0}\\
+& \frac{1}{120}\left \langle R^{TX}_{x_0} (\mR, J\mR)\mR,
R^{TX}_{x_0} (\mR, e_\ell)\mR\right \rangle_{x_0}
 + \cO(|Z|^6).
 \end{split}
\end{equation}

Let $\Gamma ^\bullet= \Gamma ^E, \Gamma ^L$ and
$R^\bullet= R^E, R^L$, respectively.
By \cite[Lemma 1.2.4]{MM07}, the Taylor coefficients of 
$\Gamma ^\bullet (e_\ell) (Z)$ at $x_0$ up to order $r$ are only 
determined by those of $R^\bullet$ up to order $r-1$, and
\begin{align}\label{0c39}
\sum_{|\alpha|=r}  (\partial^\alpha
 \Gamma ^\bullet ) _{x_0} (e_\ell) \frac{Z^\alpha}{\alpha !}
=\frac{1}{r+1} \sum_{|\alpha|=r-1}
(\partial^\alpha R^\bullet ) _{x_0}(\mR, e_\ell)
  \frac{Z^\alpha}{\alpha !}.
\end{align}
Thus by (\ref{bk2.84}), (\ref{0c39}) and since $R^{TX}$ is a $(1,1)$-form,
we obtain 
\[
R^{TX} (\mR,J\mR)=  -2\sqrt{-1}R^{TX}(z,\ov{z})\,,\quad
 \left \langle R^{TX}_{\, ;(Z,Z)} (\mR, J\mR)\mR,
e_i\right \rangle_{x_0}= -2\sqrt{-1} A_{1Z}(\ov{z}, e_i)
\] and
\begin{equation}\label{0c40}
\begin{split}
t^{-1} \Gamma ^L &(e_i)(tZ)
= -\pi \sqrt{-1} \left \langle J\mR, e_i \right \rangle_{x_0}
 -t^2\frac{\pi}{6}\left \langle R^{TX}_{x_0} (z,\ov{z})\mR,
 e_i\right \rangle_{x_0}\\
&-t^3\frac{\pi}{15}\left \langle  R^{TX}_{\, ; Z} (z,\ov{z})\mR,
e_i\right \rangle_{x_0}
 -t^4 \frac{\pi}{60}A_{1Z}(\ov{z}, e_i)
%\left\langle R^{TX}_{\, ;(Z,Z)}(z,\ov{z})\mR, e_i\right\rangle_{x_0}\\
- t^4 \frac{\pi}{180}A_{2Z}(\ov{z}, e_i)
%\left \langle R^{TX}_{x_0} (z,\ov{z})\mR,
%R^{TX}_{x_0} (\mR, e_i)\mR\right \rangle_{x_0}
 + \cO(t^5).
 \end{split}
\end{equation}
By (\ref{bk2.21}), (\ref{lm01.2}),  (\ref{lm01.27}), (\ref{0c39}) 
and (\ref{0c40}), for $t=\frac{1}{\sqrt{p}}$, we get
\begin{align}\label{bk2.85}\begin{split}
& \Gamma ^E (e_i)(Z) = \frac{1}{2} R^E_{x_0} (\mR,e_i)
+ \frac{1}{3} R^E_{\, ; Z} (\mR,e_i)  \\
&\hspace*{10mm}+\frac{1}{8}
\Big( R^E_{\, ; (Z,Z)} (\mR,e_i)
+ \frac{1}{3}\left \langle  R^{TX} (\mR,e_i)\mR,
e_k\right \rangle_{x_0} R^E_{x_0} (\mR,e_k) \Big) +  \cO(|Z|^4),\\
&\nabla'_{t, e_i} =\nabla_{e_i}
+ \frac{1}{t} \Gamma ^L (e_i)(tZ) + t \Gamma ^E (e_i)(tZ) \\
%- \frac{t}{2\kappa} \nabla_{e_i}\kappa\\
&\hspace*{10mm}
= \nabla_{0,e_i} -\frac{t^2}{6} \left\langle \pi R^{TX}_{x_0}(z,\ov{z})\mR
, e_i\right\rangle_{x_0}
+ \frac{t^2}{2} R^E_{x_0} (\mR,e_i) + \cO(t^3).
\end{split}\end{align}
\comment{
we get
\begin{align}\label{bk2.85}
&\nabla_{t, e_i} = \kappa(tZ)^{1/2} \Big [\nabla_{e_i}
+ \Big(\frac{1}{t} \Gamma ^L (e_i) + t \Gamma ^E (e_i)
%- \frac{t}{2\kappa} \nabla_{e_i}\kappa
\Big)(tZ) \Big]\kappa(tZ)^{-1/2}\\
&\hspace*{3mm}
= \nabla_{0,e_i} -\frac{t^2}{6} \left\langle \pi R^{TX}_{x_0}(z,\ov{z})\mR
- R^{TX}_{x_0} (\mR, e_k)e_k , e_i\right\rangle
+ \frac{t^2}{2} R^E_{x_0} (\mR,e_i) + \cO(t^3). \nonumber
\end{align}
}
%Let $\{\wi{w}_i(Z)\}_i$ be the parallel transport of $\{w_i\}_i$
%along the curve $[0,1]\ni u\to uZ$.
%Recall that $\cL_t$ is the restriction of $\cL^t_2$
%on $\cC^{\infty}(\R^{2n}, E_{x_0})$.
By \eqref{toe2.5} and \eqref{bk2.21},
% when we restrict on $\cC^{\infty}(\R^{2n}, E_{x_0})$,
we get,
\begin{equation}\label{bk2.87}
\begin{split}
\cL_t= -\kappa(tZ)^{1/2} g^{ij} (tZ) \Big [ \nabla'_{t, e_i}\nabla'_{t, e_j}
- t \Gamma _{ij}^l(t\,{\scriptscriptstyle\bullet}) \nabla'_{t, e_l} \Big ](Z)
\kappa(tZ)^{-1/2} \\
-   t^2\tfrac{\sqrt{-1}}{2} R^E (\wi{e}_i, J \wi{e}_i) (tZ) -2 \pi n.
\end{split}
\end{equation}
We will derive now \eqref{lm01.101} and \eqref {lm01.102} (they 
were already obtained in \cite[Th.\,4.1.25]{MM07}).
By using the Taylor expansion of the expressions from \eqref{bk2.87} 
(see \eqref{lm01.38}, \eqref{bk2.83}, \eqref{bk2.82}, \eqref{abk2.82}, 
\eqref{bk2.85}) we obtain immediately the formulas for 
$\cL_0$ and $\mO_1$ given in \eqref{lm01.101}. 

In order to compute $\mO_2$, observe first that by \eqref{alm01.0} 
and the fact that $R^{TX}$ is a (1,1)-form with values in
$\End(TX)$, we get
\begin{multline}\label{abk2.87}
\nabla_{e_j}\left \langle  R^{TX}_{x_0} (z,\ov{z}) \mR,
e_j\right \rangle
= 2\left(\tfrac{\partial}{\partial\ov{z}_j}
\left \langle  R^{TX}_{x_0} (z,\ov{z})\ov{z},
\tfrac{\partial}{\partial z_j}\right \rangle
+ \tfrac{\partial}{\partial z_j}\left \langle  R^{TX}_{x_0} (z,\ov{z})z,
\tfrac{\partial}{\partial\ov{z}_j}\right \rangle\right) =0.
\end{multline}
Thus from \eqref{lm01.9a}, \eqref{lm01.38}, \eqref{bk2.83}, \eqref{bk2.82},
\eqref{bk2.85}--\eqref{abk2.87}, we have
\begin{align}\label{0c41}
\mO_{2}= \mO^{\,\prime}_{2} + \Big[\cL_{0}, \frac{1}{12} \ric(\mR,\mR)\Big].
\end{align}
By the formula of $\cL_0$ (see \eqref{lm01.101}) and since
\[
[\cL_{0},  \ric(\mR,\mR)]= -4 \ric(\mR,e_{j})\nabla_{0,e_j}
-2 \ric(e_{j},e_{j})\,,
\]
we get from \eqref{0c41} the formula for $\mO_{2}$ given in \eqref{lm01.102}.

{}From \eqref{bk2.87}, we have also
\begin{equation}\label{0c43}
\begin{split}
\mO_3\;\;= &\;\;\;\frac{1}{6}\left \langle R^{TX}_{\, ; Z}  
(\mR,e_i) \mR, e_j\right
\rangle_{x_0} \nabla_{0,e_i} \nabla_{0,e_j}\\
&-  \Big[-\frac{2\pi}{15}\left \langle
   R^{TX}_{\, ; Z} (z,\ov{z})\mR, e_i\right\rangle_{x_0}
+ \frac{2}{3} R^E_{\, ; Z} (\mR, e_i)\Big ]\nabla_{0,e_i} \\
&-\frac{\partial}{\partial Z_i}\Big[ -\frac{\pi}{15}\left \langle
   R^{TX}_{\, ; Z} (z,\ov{z})\mR, e_i\right \rangle_{x_0}
+ \frac{1}{3} R^E_{\, ; Z}  (\mR, e_i)\Big]\\
&+ \frac{1}{12} \Big[ 4   \ric_{\, ; Z}(\mR,e_l)
+ 2 \left \langle  R^{TX}_{\, ; e_j} (\mR, e_j) \mR, e_l\right
\rangle_{x_0}
+  \ric_{\, ; e_l} (\mR,\mR)\Big ]\nabla_{0,e_l} \\
&+ \Big[\cL_0,  \frac{1}{24}  \ric_{\, ; Z}(\mR,\mR)\Big]
- \frac{\sqrt{-1}}{2} R^E_{\, ; Z} (e_i, Je_i).
\end{split}
\end{equation}
In \eqref{0c43}, the first (resp. second and third, resp. fourth,
resp. fifth) term is the contribution of the coefficient
of $t^3$ in $g^{ij}(tZ)$ (resp. $\nabla'_{t,e_i}$,
resp. $t\Gamma^l_{ii}(t\cdot)$, resp. $\kappa^{1/2}(tZ)$).
By the same argument in (\ref{abk2.87})
and  the formula of $\cL_0$ given in \eqref{lm01.101}, we get
\begin{align}\label{0c44}\begin{split}
 \frac{\partial}{\partial Z_i}
\left \langle   R^{TX}_{\, ; Z} (z,\ov{z}) \mR, e_i\right \rangle_{x_0}
=& \left \langle  R^{TX}_{\, ; e_{j}} (z,\ov{z}) \mR,
    e_j\right \rangle_{x_0},\\
[\cL_0,   \ric_{\, ; Z}(\mR,\mR)] =& -2 ( \ric_{\, ; e_i}(\mR,\mR)
+2 \ric_{\, ; Z}(\mR,e_i))\nabla_{0,e_i} \\
&- 4  \ric_{\, ; e_i}(\mR,e_i)-2  \ric_{\, ; Z}(e_i,e_i).
\end{split}\end{align}
From \eqref{0c43} and \eqref{0c44} we get the formula for $\mO_{3}$
asserted in \eqref{lm01.11}. Moreover,
\begin{equation}\label{0c45}
\begin{split}
\mO_4\;\;= &\;\;\;\mO_{42}
+\Big[\mO^{\,\prime}_2, \frac{1}{12} \ric(\mR,\mR)\Big] +\mO_{43}
+ \mO_{41}\\
 &+\frac{1}{3}
\left \langle R^{TX}_{x_0} (\mR,e_i) \mR, e_j\right \rangle_{x_0}
\Big\{ \Big(-\frac{1}{3} \left\langle \pi R^{TX}_{x_0}(z,\ov{z})\mR ,
e_j\right\rangle_{x_0}
+  R^E_{x_0} (\mR,e_j)\Big)\nabla_{0,e_i} \\
&\hspace{20mm}-\frac{1}{6} \frac{\partial}{\partial Z_i}
\left\langle \pi R^{TX}_{x_0}(z,\ov{z})\mR , e_j\right\rangle_{x_0}
- \frac{2}{3} \left \langle R^{TX}_{x_0} (\mR, e_i) e_j,
e_l\right \rangle _{x_0} \nabla_{0,e_l}\Big\} \\
&- \Big(-\frac{1}{6} \left\langle \pi R^{TX}_{x_0}(z,\ov{z})\mR, 
e_i\right\rangle_{x_0}
+ \frac{1}{2} R^E_{x_0} (\mR,e_i)
- \frac{2}{3} \left \langle R^{TX}_{x_0} (\mR, e_j) e_j,
e_i\right \rangle _{x_0}\Big)\\
&\hspace{45mm}
\times  \Big(-\frac{1}{6} \left\langle \pi R^{TX}_{x_0}(z,\ov{z})\mR
, e_i\right\rangle_{x_0}
+ \frac{1}{2} R^E_{x_0} (\mR,e_i)\Big)\\
&-\Big \{ -\mO_{44} + \Big[- \frac{\pi}{9} A_{2Z}(\ov{z},e_i)
%-\frac{\pi}{30} A_{1Z}(\ov{z},e_i)
%- \frac{\pi}{90} A_{2Z}(\ov{z},e_i)
+ \frac{1}{12}\left \langle  R^{TX} (\mR,e_i)\mR,
e_m\right \rangle_{x_0} R^E_{x_0} (\mR,e_m)\\
&\hspace{20mm}+ \frac{2}{9} \left \langle R^{TX}_{x_0} (\mR,e_i) \mR,
e_k\right\rangle_{x_0} \ric(\mR,e_k)+  \frac{1}{4} R^E_{\, ;(Z,Z)} (\mR,  e_i)
%-  \tfrac{\partial}{\partial Z_j} \Big(\frac{1}{20}A_{1Z}
%+  \frac{2}{45}A_{2Z} \Big)(e_{j}, e_{i})
%+  \frac{1}{2} \tfrac{\partial}{\partial Z_i} \Big(\frac{1}{20}A_{1Z}
% +  \frac{2}{45}A_{2Z} \Big)(e_{j}, e_{j})\Big\}
\Big]\nabla_{0,e_i} \Big\}\\
&-  \frac{\partial}{\partial Z_i}
\Big(-\frac{\pi}{60}A_{1Z}(\ov{z},e_i)
%\left\langle R^{TX}_{\, ;(Z,Z)} (z,\ov{z})\mR, e_i\right\rangle_{x_0}
- \frac{\pi}{180}A_{2Z}(\ov{z},e_i)
%\left\langle R^{TX}_{x_0} (z,\ov{z})\mR,
%R^{TX}_{x_0} (\mR, e_i)\mR\right\rangle_{x_0}\Big)\\
+  \frac{1}{8} R^E_{\, ;(Z,Z)} (\mR,  e_i)\\
&\hspace{55mm}+ \frac{1}{24}\left \langle  R^{TX}_{x_0} (\mR,e_i)\mR,
e_k\right \rangle_{x_0} R^E_{x_0} (\mR,e_k)\Big)\\
&-\frac{\sqrt{-1}}{4} R^E_{\, ;(Z,Z)} (e_i, J e_i)\,.
% + \frac{\sqrt{-1}}{6} \left \langle R^{TX}_{x_0} (\mR,e_i) \mR,
%e_l\right \rangle_{x_0} R^E_{x_0} (e_i, J e_l),
\end{split}
\end{equation}%{multline}
Here
\begin{itemize}
\item $\mO_{42}$ is the contribution of
the coefficients of $t^4$ in $\kappa^{1/2}(tZ)$ and $\kappa^{-1/2}(tZ)$,
\item the second term is the contribution of
the coefficients of $t^2$ in
$\kappa ^{1/2}(tZ)$, $\kappa ^{-1/2}(tZ)$ and in
$-g^{ij}(tZ) (\nabla'_{t,e_i}\nabla'_{t,e_j}
- t \Gamma^\ell_{ij}(tZ)\nabla'_{t,e_\ell})$,
\item $\mO_{43}$ is the contribution of
the coefficients of $t^2$ in
$\kappa ^{1/2}(tZ)$ and $\kappa ^{-1/2}(tZ)$,
\item $\mO_{41}$ is the contribution of
the coefficients of $t^4$ in $g^{ij}(tZ)$,
\item the fifth term is the contribution of
the coefficients of $t^2$ in
$-g^{ij}(tZ)$ and in
$ (\nabla'_{t,e_i}\nabla'_{t,e_j}
- t \Gamma^\ell_{ij}(tZ)\nabla'_{t,e_\ell})$,
\item  the sixth, seventh and eight terms are the contributions of
the coefficients of $t^4$ in $-(\nabla'_{t,e_i}\nabla'_{t,e_i}
- t \Gamma^\ell_{ii}(tZ)\nabla'_{t,e_\ell})$:
 the sixth term is the contribution of the coefficients of $t^2$ in
$\nabla'_{t,e_i}$, $- t \Gamma^\ell_{ii}(tZ)$ and $t^2$ in $\nabla'_{t,e_i}$;
 the seventh and eighth terms are the contributions 
of the coefficients of $t^4$ in $\nabla'_{t,e_i}$ and $-t \Gamma^\ell_{ii}(tZ)$.
%\item the eighth is the contribution of the coefficients of 
%$t^4$ in $-t \Gamma^l_{ii}(tZ)$

\end{itemize}
%Here $\mO_{42}$ (resp. the second term, resp. $\mO_{43}$)
%is the contribution of
%the coefficients of $t^4$ in $\kappa ^{1/2}(tZ)$ (resp. $t^2$ in
%$\kappa ^{1/2}(tZ)$, $\kappa ^{-1/2}(tZ)$ and  $t^2$ in
%$-g^{ij}(tZ) (\nabla'_{t,e_i}\nabla'_{t,e_j}
%- t \Gamma^l_{ij}(tZ)\nabla'_{t,e_l})$, resp. $t^2$ in
%$\kappa ^{1/2}(tZ)$ and $t^2$ in $\kappa ^{-1/2}(tZ)$ );
%$\mO_{41}$ (resp. the fifth term) is the contribution of
%the coefficients of $t^4$ in $g^{ij}(tZ)$ (resp. $t^2$ in
%$-g^{ij}(tZ)$ and  $t^2$ in
%$ (\nabla'_{t,e_i}\nabla'_{t,e_j}
%- t \Gamma^l_{ij}(tZ)\nabla'_{t,e_l})$),
%the rest terms until the last third term are  the contribution of
%the coefficients of $t^4$ in $-(\nabla'_{t,e_i}\nabla'_{t,e_i}
%- t \Gamma^l_{ii}(tZ)\nabla'_{t,e_l})$: the sixth (resp. seventh and
%eighth) term is the contributions of the coefficients of $t^2$ in
%$\nabla'_{t,e_i}$, $- t \Gamma^l_{ii}(tZ)$ and $t^2$ in $\nabla'_{t,e_i}$
%(resp. $t^4$ in $\nabla'_{t,e_i}$, $-t \Gamma^l_{ii}(tZ)$).

\noindent
Now by \eqref{lm01.9a},
\begin{equation}\label{0c46}
\begin{split}
\frac{1}{12}[\mO^{\,\prime}_2,  \ric(\mR,\mR)] =
\frac{1}{36}\left \langle R^{TX}_{x_0} (\mR,e_i) \mR, e_j\right
\rangle_{x_0}
( 4 \ric(\mR,e_j) \nabla_{0,e_i} + 2 \ric(e_i,e_j))\\
+ \frac{1}{6} \left(\left \langle
   \frac{\pi}{3} R^{TX}_{x_0} (z,\ov{z}) \mR,  e_j\right \rangle
  +  \frac{2}{3}  \ric (\mR,e_j)  -R^E_{x_0} (\mR, e_j)\right)\ric(\mR,e_j).
  \end{split}
\end{equation}
and the same argument used to obtain \eqref{abk2.87} shows that
\begin{align} \label{0c47}
 \begin{split}
&\left \langle R^{TX}_{x_0} (\mR,e_i) \mR, e_j\right \rangle_{x_0}
\frac{\partial}{\partial Z_i}
\left\langle R^{TX}_{x_0}(z,\ov{z})\mR , e_j\right\rangle
= \left \langle R^{TX}_{x_0} (\mR,e_i) \mR,
R^{TX}_{x_0}(z,\ov{z}) e_i\right\rangle,\\
& \frac{\partial}{\partial Z_i}(A_{1Z}(\ov{z},e_i))
= \frac{\partial}{\partial Z_i}\left \langle
   R^{TX}_{\, ;(Z,Z)} (z,\ov{z})\mR, e_i\right \rangle_{x_0}
= 2 \left \langle   R^{TX}_{\, ;(Z,e_i)} (z,\ov{z})\mR, e_i\right \rangle_{x_0},\\
&\frac{\partial}{\partial Z_i}(A_{2Z}(\ov{z},e_i))
=  \frac{\partial}{\partial Z_i}\left \langle R^{TX}_{x_0} (z,\ov{z})\mR,
R^{TX}_{x_0} (\mR, e_i)\mR\right \rangle_{x_0}\\
&\hspace{10mm}= \left \langle R^{TX}_{x_0} (z,\ov{z})e_i,
R^{TX}_{x_0} (\mR, e_i)\mR\right \rangle_{x_0}
+ \left \langle R^{TX}_{x_0} (z,\ov{z})\mR, e_j \right \rangle_{x_0}
\ric(\mR,e_j) .
\end{split}\end{align}
Finally, by \eqref{alm01.0} and since $R^{TX}$ is a $(1,1)$-form, we obtain
\begin{align}\label{0c50}\begin{split}
&    \left \langle R^{TX}_{x_0} (z,\ov{z})e_\ell,
    R^{TX}_{x_0} (\mR, e_\ell)\mR\right \rangle
=  \left \langle R^{TX}_{x_0} (z,\ov{z})e_\ell, e_m \right \rangle
  \left \langle   R^{TX}_{x_0} (\mR, e_\ell)\mR, e_m\right \rangle =0,\\
&\left \langle R^{TX}_{x_0} (\mR,e_\ell) \mR, e_m\right
\rangle_{x_0} \ric(e_\ell,e_m)
= -8 R_{k\ov{\ell}m\ov{q}} \ric_{\ell\ov{m}} \, z_{k} \ov{z}_{q}.
\end{split}\end{align}
Thus
\begin{equation}\label{0c49}
\begin{split}
\mO_4=&\sum_{\alpha=1}^{5}\mO_{4\alpha}
 -\frac{\pi^2}{36} \left\langle R^{TX}_{x_0}(z,\ov{z})\mR,
 R^{TX}_{x_0}(z,\ov{z})\mR\right\rangle -  \frac{1}{4}R^E_{x_0} (\mR,e_j)^2 \\
&+ \frac{1}{6}\left\langle \pi R^{TX}_{x_0}(z,\ov{z})\mR , e_j\right\rangle
 R^E_{x_0} (\mR,e_j)
+ \frac{\pi}{30} \left \langle
   R^{TX}_{\, ;(Z,e_i)} (z,\ov{z})\mR, e_i\right \rangle_{x_0}\\
%+  (-\frac{2\pi}{18}+\frac{\pi}{18} + \frac{\pi}{180})
&-\frac{\pi}{20} \left \langle R^{TX}_{x_0} (z,\ov{z})\mR,
e_j \right \rangle_{x_0} \ric(\mR,e_j)
 + \frac{1}{9} \ric(\mR,e_j)\ric(\mR,e_j)\\
&+  \frac{1}{18}\left \langle R^{TX}_{x_0} (\mR,e_i) \mR, e_j\right
\rangle_{x_0} \ric(e_i,e_j)
+    %(\frac{1}{3}- \frac{1}{6}- \frac{1}{24})
\frac{1}{8} \ric(\mR,e_j)R^E_{x_0} (\mR,e_j)\\
& - \frac{1}{4} R^E_{\, ;(Z,e_i)} (\mR,  e_i)
- R^E_{\, ;(Z,Z)}
(\tfrac{\partial}{\partial z_i},\tfrac{\partial}{\partial \ov{z}_i}).
 %- \frac{4}{3} R_{k\ov{i}l\ov{j}} z_k \ov{z}_j R^E_{i\ov{l}}.
 \end{split}
\end{equation}
\noindent
{}Putting together Lemma \ref{lmt1.6}, \eqref{lm01.7}, \eqref{lm01.9},  \eqref{0c50},
\eqref{0c49} and the fact that $R^{TX}$ is a $(1,1)$-form,
we infer \eqref{lm01.12}.
The proof of Theorem \ref{bkt2.21} is completed.
\end{proof}

%\newpage
%%%%%%%%%%%%%%%%%%%%%%%%%%%%%%%%%%%%%%%%%%
%%%%%%%%%%%%%%%%%%%%%%%%%%%%%%%%%%%%%%%%%%
\section{Evaluation of $\cF_4$ from \eqref{bk2.32}}\label{bks3}
%%%%%%%%%%%%%%%%%%%%%%%%%%%%%%%%%%%%%%%%%%

We calculate in this section an explicit formula for the operator $\cF_4$, 
defined by \eqref{bk2.24} and appearing in the Bergman kernel expansion 
\eqref{toe2.9}. This
is necessary in Section \ref{toes4} in order to evaluate the
expansion of the kernel of the Berezin-Toeplitz operators.
We use formula \eqref{bk2.32} to achieve our aim.
Recall that explicit formulas for the operators $\cL$, $\mO_2$, $\mO_4$ 
appearing in \eqref{bk2.32} were given in Theorem~\ref{bkt2.21}. 

%We will analyze the contribution
% of $\mO_2$, $\mO_4$ in \eqref{bk2.32}. This
%%The analysis of  Sections \ref{bks3.2}, \ref{bks3.3}
%is necessary in Section \ref{toes4} to evaluate the
%expansion of the kernel of the Berezin-Toeplitz operators.

This section is organized as follows.
%In Section \ref{bks3.1}, we state a formula for $\cF_4(0,0)$
%in Theorem \ref{bkt3.0}, and we verify it is compatible to
%Riemann-Roch-Hirzebruch Theorem. The proof of Theorem \ref{bkt3.0}
%is given in Sections \ref{bks3.2}, \ref{bks3.3}.
In Section \ref{bks3.2}, we determine the terms in \eqref{bk2.32} which 
involve $\mO_{2}$. In 
Section \ref{bks3.3}, we calculate the terms in \eqref{bk2.32} which involve
$\mO_{4}$.
In Section \ref{bks3.4}, we
obtain the formula for $\cF_4(0,0)$ (cf.\;Theorem \ref{bkt3.0}).

We adopt the convention that all tensors will be evaluated at the base
point $x_0\in X$ and most of time, we will omit the subscript $x_0$.

%%%%%%%%%%%%%%%%%%%%%%%%%%%%%%%%%%%%%%%
\subsection{Contribution of
$\mO_{2}$ to $\cF_4(0,0)$}\label{bks3.2}
%%%%%%%%%%%%%%%%%%%%%%%%%%%%%%%%%%%%%%%%
\begin{lemma} \label{bkt3.1}  The following identities hold:
\begin{subequations}
\begin{align}
&\cP \mO_2 \cP=0,\label{bk3.0a}\\
&(\cL^{-1} \mO_2 \cP\mO_2\cL^{-1})(0,0)
=\frac{1}{4\pi^2} \Big(\sum_{km}R_{m\ov{m}k\ov{k}}
+ \sum_k R^E_{k\ov{k}}\Big)^2,\label{bk3.0b}\\
& \label{bk3.0c}(\cP \mO_2\cL^{-2} \mO_2 \cP)(0,0)\\
&\hspace{10mm} =
\frac{1}{36\pi^2}  R_{m\ov{k}q\ov{\ell}}  R_{k\ov{m}\ell\ov{q}}
+ \frac{1}{4\pi^2} \Big(\frac{4}{3}R_{q\ov{m}m\ov{\ell}} 
+  R^E_{q\ov{\ell}}\Big)
\Big(\frac{4}{3} R_{\ell\ov{k}k\ov{q}} + R^E_{\ell\ov{q}}\Big).\nonumber
\end{align}
\end{subequations}
\end{lemma}
\begin{proof}
Note that by \eqref{toe1.1} and (\ref{toe1.3}),
\begin{equation} \label{bk3.1}
(b^+_i\cP)(Z,Z^{\prime}) =0\,,\quad \,
(b_i\cP)(Z,Z^{\prime})=2\pi (\ov{z}_i-\ov{z}^{\,\prime}_i)\cP(Z,Z^{\prime}).
\end{equation}
For $\phi\in T^{*}X$, by (\ref{lm01.3}),
(\ref{bk2.64}) and (\ref{bk3.1}), we have
\begin{align}\label{bk3.2}\begin{split}
&\phi(e_{i}) e_{i}= 2 \phi(\tfrac{\partial}{\partial z_j})
\tfrac{\partial}{\partial\ov{z}_j}
+ 2 \phi(\tfrac{\partial}{\partial\ov{z}_j})
\tfrac{\partial}{\partial z_j},\quad
\phi(e_{i}) \nabla_{0, e_{i}}=  \phi(\tfrac{\partial}{\partial z_j})
b^{+}_{j}
-  \phi(\tfrac{\partial}{\partial\ov{z}_j}) b_{j},\\
&\phi(e_{i}) \nabla_{0, e_{i}}\cP(Z,0)
= -2\pi \phi(\ov{z}) \cP(Z,0) .
\end{split}\end{align}
\noindent
By Lemma \ref{lmt1.6}, \eqref{alm01.0},  \eqref{lm01.4}, \eqref{lm01.102},
  \eqref{bk3.2} and the fact that $R^{TX}$ is a $(1,1)$-form, we get
\begin{equation} \label{bk3.3}
\begin{split}
    \mO_2 & = \frac{1}{3} R_{k\ov{m}\ell\ov{q}} z_{k} z_{\ell} b_{m} b_{q}
    + \frac{1}{3} R_{k\ov{q}\ell\ov{m}} z_{k} \ov{z}_{m} (b_{q} b^+_{\ell}+
    b^+_{\ell}b_{q} )
  +  \frac{1}{3} R_{k\ov{m}\ell\ov{q}}  \ov{z}_{m} \ov{z}_{q} b^+_{k} b^+_{\ell}
  -  \frac{4}{3}  R_{m\ov{m}q\ov{q}}\\
&\quad+ \Big( \frac{2}{3}R_{k\ov{m}\ell\ov{k}}  \ov{z}_{m}
- \frac{\pi}{3} R_{k\ov{m}\ell\ov{q}} z_k \ov{z}_{m} \ov{z}_{q}\Big) b^+_l
-  \Big( \frac{2}{3}R_{\ell\ov{k}k\ov{q}}  z_\ell
+ \frac{\pi}{3} R_{k\ov{m}\ell\ov{q}} z_k \ov{z}_{m} z_{\ell}\Big) b_q\\
&\quad- 2 R^E_{q\ov{q}}
- R^E(\ov{z},\tfrac{\partial}{\partial z_\ell}) b^+_\ell
+  R^E(z,\tfrac{\partial}{\partial\ov{z}_q}) b_q\,.
\end{split}
\end{equation}
By Lemma \ref{lmt1.6}, \eqref{bk2.66} and \eqref{bk3.3}, we get
$ [R_{k\ov{m}\ell\ov{q}} z_{k} z_{\ell}, b_{m} b_{q}]
= 8 b_q R_{k\ov{k}\ell\ov{q}}  z_{\ell} + 8 R_{m\ov{m}q\ov{q}}$, and
\begin{equation} \label{bk3.4}
\begin{split}
    \mO_2 & = \frac{1}{3} b_{m} b_{q} R_{k\ov{m}\ell\ov{q}} z_{k} z_{\ell}
    +  b_{q} \Big(- \frac{\pi}{3} R_{k\ov{m}\ell\ov{q}} z_k  z_{\ell}\ov{z}_{m}
+ 2R_{\ell\ov{k}k\ov{q}}  z_\ell+R^E_{\ell\ov{q}}z_\ell\Big)\\
&\quad+ \Big(  \frac{2 b_q}{3}  R_{k\ov{m}\ell\ov{q}}  z_{k} \ov{z}_{m}
- \frac{\pi}{3} R_{k\ov{m}\ell\ov{q}} z_k \ov{z}_{m}\ov{z}_{q}
+ 2R_{k\ov{k}\ell\ov{m}} \ov{z}_{m} 
+R^E_{\ell\ov{m}}\ov{z}_m \Big) b^+_{\ell}\\
&\quad+  \frac{1}{3}  R_{k\ov{m}\ell\ov{q}}  \ov{z}_{m} \ov{z}_{q} 
b^+_{k} b^+_{\ell}.
\end{split}
\end{equation}
Thus Lemma \ref{lmt1.6}, \eqref{bk2.66}, \eqref{bk3.1} and \eqref{bk3.4} yield
\begin{equation} \label{bk3.5}
\begin{split}
\mO_2\cP  &= \Big\{\frac{1}{3} b_{m} b_{q} R_{k\ov{m}\ell\ov{q}} z_{k} z_{\ell}
 +  b_{q} \Big(- \frac{\pi}{3} R_{k\ov{m}\ell\ov{q}} z_k z_{\ell}
(\frac{b_m}{2\pi} +\ov{z}_{m}')
+ 2R_{\ell\ov{k}k\ov{q}}  z_\ell +R^E_{\ell\ov{q}}z_\ell\Big)\Big\}\cP\\
&= \Big\{\frac{1}{6} b_{m} b_{q} R_{k\ov{m}\ell\ov{q}} z_{k} z_{\ell}
    +  \frac{4}{3}b_{q} R_{\ell\ov{k}k\ov{q}} z_\ell
- \frac{\pi}{3}b_{q}  R_{k\ov{m}\ell\ov{q}} z_kz_\ell \ov{z}_{m}'
+ b_q  R^E_{\ell\ov{q}}z_\ell\Big\}\cP .
\end{split}
\end{equation}
Now, \eqref{bk3.0a} follows from Theorem \ref{bkt2.17}, \eqref{toe1.4} and \eqref{bk3.5}. These imply also
\begin{equation} \label{bk3.6}
\begin{split}
&\cL^{-1} \mO_2 \cP =
\Big\{\frac{b_{m} b_{q}}{48\pi}  R_{k\ov{m}\ell\ov{q}}\, z_{k}\, z_{l}
    +  \frac{b_{q}}{3\pi} R_{\ell\ov{k}k\ov{q}}\, z_\ell
- \frac{b_{q}}{12} R_{k\ov{m}\ell\ov{q}}\, z_k\,z_\ell\, \ov{z}_{m}^{\,\prime}
+ \frac{b_q}{4\pi}  R^E_{\ell\ov{q}}\,z_\ell\Big\}\cP .
\end{split}
\end{equation}
Due to \eqref{bk2.66} and \eqref{bk3.6} we have
\begin{equation} \label{bk3.7}
\begin{split}
&(\cL^{-1} \mO_2 \cP)(Z,0) =
\Big\{\frac{b_{m} b_{q}}{48\pi}  R_{k\ov{m}\ell\ov{q}} z_{k} z_{\ell}
    + \frac{b_q}{4\pi}\Big( \frac{4}{3} R_{\ell\ov{k}k\ov{q}}
+   R^E_{\ell\ov{q}}\Big)z_\ell\Big\}\cP(Z,0) ,\\
&(\cL^{-1} \mO_2 \cP)(0,Z) =
- \frac{1}{2\pi} (R_{m\ov{m}q\ov{q}} + R^E_{q\ov{q}})\cP(0,Z).
\end{split}
\end{equation}
Since  $\mO_2$, $\cL$ are symmetric (as explained in Remark
\ref{toet2.7}) and $(R^E_{\ell\ov{q}})^{*} = R^E_{q\ov{\ell}}$\,, we get
by \eqref{toe1.1}, \eqref{lm01.4} and \eqref{bk3.7},
\begin{equation} \label{bk3.8}
\begin{split}
&(\cP \mO_2 \cL^{-1})(0,Z) =  ((\cL^{-1} \mO_2 \cP)(Z,0))^* \\
&\hspace{32mm}
= \Big\{\cP \Big[ R_{u\ov{v}q\ov{s}} \ov{z}_{v} \ov{z}_{s}
\frac{b^+_{u} b^+_{q}}{48\pi}
  + \Big(\frac{4}{3} R_{q\ov{v}v\ov{s}} + R^E_{q\ov{s}}\Big)
\ov{z}_{s}\frac{b^+_{q}}{4\pi}\Big]\Big\}(0,Z) ,\\
&(\cP \mO_2 \cL^{-1})(Z,0) = ( (\cL^{-1} \mO_2 \cP)(0,Z))^*
=- \frac{1}{2\pi} (R_{k\ov{k}q\ov{q}} + R^E_{q\ov{q}})\cP(Z,0).
\end{split}
\end{equation}
Note that $\cP(0,0)=1$ by \eqref{toe1.3}.
From \eqref{toe1.4}, \eqref{bk2.66}, \eqref{bk3.1}, \eqref{bk3.7} and \eqref{bk3.8},
we get  \eqref{bk3.0b}, and
\begin{equation} \label{bk3.9}
\begin{split}
(\cP \mO_2\cL^{-2}  \mO_2 \cP)(0,0)
=  \cP \Big\{&\frac{32\pi^2}{(48\pi)^2}
R_{m\ov{s}q\ov{t}} \ov{z}_{s} \ov{z}_{t}
R_{k\ov{m}\ell\ov{q}} z_{k} z_{\ell} \\
&+  \frac{1}{4\pi} \Big(\frac{4}{3}R_{q\ov{s}s\ov{t}}
+  R^E_{q\ov{t}}\Big)  \ov{z}_{t}
\Big(\frac{4}{3} R_{\ell\ov{k}k\ov{q}} +
R^E_{\ell\ov{q}}\Big)z_\ell\Big\}\cP(0,0) .
\end{split}
\end{equation}
Let $\phi\in\C[b,z]$ be a polynomial in
$b,z$. By \eqref{bk2.66},  \eqref{bk3.1}, we have
\begin{align}\label{bk3.10}\begin{split}
&(\ovz_k \phi(b,z) \cP)(Z,0) =
\phi(b,z)\frac{b_k}{2\pi} \cP(Z,0)
= \Big(\frac{b_k}{2\pi}\phi
+ \frac{1}{\pi}\frac{\partial \phi}{\partial z_k}\Big)\cP(Z,0),\\
&(\ovz_k\ovz_l \phi(b,z) \cP)(Z,0)=(\phi(b,z)\ovz_k\ovz_l \cP)(Z,0) \\
&\hspace{30mm} =  \Big(\frac{b_k b_l}{4\pi^2}\phi
+\frac{b_k}{2\pi^2}\frac{\partial \phi}{\partial z_l}
+\frac{b_l}{2\pi^2}\frac{\partial \phi}{\partial z_k}
+ \frac{1}{\pi^2}\frac{\partial^2 \phi}{\partial z_k\partial z_l}
\Big)\cP(Z,0).
\end{split}\end{align}

\noindent
Let $F(Z)$ be a homogeneous degree $2$ polynomial in $Z$.
By \eqref{bk3.1},
\begin{align}\label{bk3.15}
\left(F(Z)\cP\right)(Z,0) = \Big(\frac{1}{2}
\frac{\partial ^2 F}{\partial z_i\partial z_j} z_iz_j +
\frac{\partial ^2 F}{\partial z_i\partial \ov{z}_j} z_i
\frac{b_j}{2\pi} + \frac{1}{2}\frac{\partial ^2 F}{\partial
\ov{z}_i\partial \ov{z}_j}
 \frac{b_ib_j}{4\pi^2} \Big) \cP(Z,0).
\end{align}
Thus from \eqref{toe1.4}, \eqref{bk3.10} and \eqref{bk3.15},
we have
%\begin{subequations}
\begin{align}\label{abk3.8}\begin{split}
&(\cP F \cP)(Z,0)=
\Big(\sum_{|\alpha|=2}\frac{\partial^2 F}{\partial {z}^\alpha}
\frac{z^\alpha}{\alpha !}
+\frac{1}{\pi} \frac{\partial ^2 F}{\partial z_{i} \partial \ov{z}_{i}}
 \Big) \cP(Z,0),\\
&(\cP F \ov{z}_i\ov{z}_j \cP)(0,0)= \frac{1}{\pi^2}
 \frac{\partial ^2 F}{\partial z_{i} \partial z_{j}}.
\end{split}\end{align}
%\end{subequations}
By  \eqref{bk3.9} and \eqref{abk3.8}, we get \eqref{bk3.0c}.
The proof of Lemma \ref{bkt3.1} is completed.
\end{proof}

\begin{lemma}\label{bkt3.2} The following identity holds
\begin{equation}\label{bk3.12}
\begin{split}
\pi^2 (\cL^{-1}\cP^{\bot}\mO_2\cL^{-1} \mO_2 \cP)(0,0)=
- \frac{25}{2^3\cdot 3^3} R_{m\ov{k}q\ov{\ell}}R_{k\ov{m}\ell\ov{q}}
-  \frac{47}{54} R_{k\ov{k}q\ov{\ell}} R_{m\ov{m}\ell\ov{q}}  \\
+  \frac{1}{8} R_{k\ov{k}\ell\ov{\ell}} R_{m\ov{m}q\ov{q}}
+ \frac{1}{4}R^E_{\ell\ov{\ell}}  R_{m\ov{m}q\ov{q}}
- \frac{7}{6} R^E_{q\ov{\ell}}  R_{m\ov{m}\ell\ov{q}}
+ \frac{1}{8} (R^E_{\ell\ov{\ell}} R^E_{q\ov{q}}
-3  R^E_{q\ov{\ell}} R^E_{\ell\ov{q}}).
\end{split}
\end{equation}
\end{lemma}
\begin{proof}
Set
\begin{align}\label{bk3.13}\begin{split}
\bI_1 =& \Big\{\frac{1}{3} b_{k} b_{\ell} R_{s\ov{k}t\ov{\ell}} z_{s} z_{t}
    +  b_{\ell} \Big[- \frac{\pi}{3} R_{s\ov{k}t\ov{\ell}}
z_{s} z_{t} \ov{z}_{k}
+ \big(2R_{t\ov{s}s\ov{\ell}} + R^E_{t\ov{\ell}}) z_{t}\big) \Big]\Big\}
\\
&\times \Big[\frac{b_{m} b_{q}}{48\pi}  R_{i\ov{m}j\ov{q}} z_{i} z_{j}
    +  \frac{b_{q}}{4\pi} \big( \frac{4}{3}R_{j\ov{i}i\ov{q}} +
    R^E_{j\ov{q}}\big) z_j\Big],\\
\bI_2=& \Big[  \frac{2 b_{\ell}}{3}  R_{s\ov{k}q\ov{\ell}}  z_{s} \ov{z}_{k}
- \frac{\pi}{3} R_{s\ov{k}q\ov{\ell}} z_{s}\ov{z}_{k}\ov{z}_{\ell}
+ \big(2R_{s\ov{s}q\ov{k}} + R^E_{q\ov{k}}\big) \ov{z}_{k} \Big]\\
&\times \Big(\frac{b_{m}}{6}  R_{i\ov{m}j\ov{q}} z_{i} z_{j}
    +  \frac{4}{3} R_{j\ov{i}i\ov{q}} z_j + R^E_{j\ov{q}}z_j \Big).
\end{split}\end{align}
Then by Lemma \ref{lmt1.6}, \eqref{bk2.66}, \eqref{bk3.1}, \eqref{bk3.4}
and  \eqref{bk3.7}, we get as in \eqref{bk3.9} that
\begin{align} \label{bk3.14}
(\mO_2\cL^{-1} \mO_2 \cP)(Z,0)= \Big\{\bI_1 +
\bI_2
+  \frac{2\pi}{9} R_{m\ov{s}q\ov{t}}\,\ov{z}_{s}\ov{z}_{t}
 R_{k\ov{m}\ell\ov{q}} z_{k} z_{\ell}\Big\} \cP (Z,0).
\end{align}
Let $h(z)=\sum_i h_i z_i $, $h'(z)=\sum_i h'_i z_i$
with $h_i, h'_i\in\C$, and
let $F(Z)$ be a homogeneous degree $2$ polynomial in $Z$.
\comment{
By \eqref{bk3.1},
\begin{align}\label{bk3.15}
\left(F(Z)\cP\right)(Z,0) = \Big(\frac{1}{2}
\frac{\partial ^2 F}{\partial z_i\partial z_j} z_iz_j +
\frac{\partial ^2 F}{\partial z_i\partial \ov{z}_j} z_i
\frac{b_j}{2\pi} + \frac{1}{2}\frac{\partial ^2 F}{\partial
\ov{z}_i\partial \ov{z}_j}
 \frac{b_ib_j}{4\pi^2} \Big) \cP(Z,0).
\end{align}
}By Theorem \ref{bkt2.17}, \eqref{toe1.4}, \eqref{bk2.66},
\eqref{bk3.1}, \eqref{bk3.10} and \eqref{bk3.15}, we have
\begin{subequations}
\begin{align}
&(\cP^{\bot} F\cP)(0,0)
=-\frac{1}{\pi} \frac{\partial ^2 F}{\partial z_i\partial \ov{z}_i}\;,
%\quad (b_{i}b_{j}F\cP)(0,0)
%=4  \frac{\partial ^2 F}{\partial z_i\partial z_j},
\label{bk3.16a}\\
%&\left(\cL^{-1} b_ib_j \cP\right)(0,0) =
%\left(\cL^{-1} \cP^{\bot} h\cP\right)(0,0) = 0,\label{bk3.16b}\\
&\left(\cL^{-1} \cP^{\bot} h b_i \cP\right)(0,0)=
\left(\cL^{-1} b_ih \cP\right)(0,0)
=-\frac{1}{2\pi} h_i\;,\label{bk3.16c}\\
%\frac{\partial h}{\partial z_i},\\
&\left(\cL^{-1} \cP^{\bot} F \cP\right)(0,0)= -\frac{1}{4\pi ^2}
\frac{\partial ^2 F}{\partial z_i\partial \ov{z}_i}\;,\label{bk3.16d}\\
&\left(\cL^{-1} b_j F b_i \cP\right)(0,0)=
-\left(\cL^{-1}  b_i b_j F  \cP\right)(0,0)= -\frac{1}{2\pi}
\frac{\partial ^2 F}{\partial z_i\partial z_j}\;,\label{bk3.16e}\\
&\left(\cL^{-1} \cP^{\bot}F b_i b_j \cP\right)(0,0)
= \cL^{-1} \cP^{\bot}\Big( b_jF  b_i
+ 2 \frac{\partial F}{\partial z_j}  b_i\Big) \cP(0,0)
= -\frac{3}{2\pi} \frac{\partial ^2 F}{\partial z_i\partial z_j}\;,
\label{bk3.16f}\\
&\left(\cL^{-1} \cP^{\bot}F \ov{z}_{i} \ov{z}_{j}\cP\right)(0,0)
=\cL^{-1} \cP^{\bot}F\,\frac{b_{i} b_{j}}{4\pi^2}\cP(0,0)
= -\frac{3}{8\pi^3} \frac{\partial ^2 F}{\partial z_i\partial z_j}\;,
\label{bk3.16g}\\
%&\left(\cL^{-1} \cP^{\bot}  b_i h  b_j h'
%\cP\right)(0,0)= -\frac{1}{2\pi} (h_j h'_i- h_i h'_j),\\
%\Big(\frac{\partial h}{\partial z_j} \frac{\partial h'}{\partial z_i}
%- \frac{\partial h}{\partial z_i}\frac{\partial h'}{\partial z_j}
%\Big),\\
 &\left(\cL^{-1} b_k F \ov{z}_{l} \cP\right)(0,0)
=\cL^{-1}  b_k F  \frac{b_{l}}{2\pi}\cP(0,0)
= -\frac{1}{4\pi^{2}} \frac{\partial ^2 F}{\partial z_k\partial  z_l}\;.
\label{bk3.16h}
\end{align}
\end{subequations}
In \eqref{bk3.16d} and \eqref{bk3.16e}, we have used
$F b_i= b_i F + 2 \frac{\partial F}{\partial z_i}$.
Observe that \eqref{toe1.3},  \eqref{bk2.66} imply that for every 
homogeneous degree $k$ polynomial $G$ in $Z$, and every
$\alpha\in\N^n$, we have
\begin{align}\label{abk3.16}
(b^\alpha G\cP)(0,0) = \begin{cases} 0\,, &\quad \text{ if  $|\alpha|\neq k$},\\
(-2)^k  \displaystyle\frac{\partial^\alpha G}{\partial
  z^\alpha}\,,&\quad \text{ if  $|\alpha|=k$}\;.
\end{cases}
\end{align}
By Theorem \ref{bkt2.17},  \eqref{toe1.4} and \eqref{bk2.66}, we also have
\begin{subequations}
\begin{align}
 &\left(\cL^{-1}  b_{i} h b_{j} h' \cP\right)(0,0)
= \left(\frac{b_j b_{i}}{8\pi} h h' + \frac{b_{i}}{2\pi} h_{j}
h'\right) \cP(0,0) = \frac{1}{2\pi} ( h_{i}h'_{j}- h_{j}h'_{i}),\label{bk3.17a} \\
%&\hspace{30mm}= \frac{1}{2\pi} ( h_{j'}h'_{j}- h_{j}h'_{j'}),\nonumber\\
&\left(\cL^{-1} \cP^{\bot}  h b_i  h' b_j
\cP\right)(0,0)
= -\frac{1}{2\pi}h_j h'_i  -\frac{3}{2\pi}h_i h'_j.\label{bk3.17b}
%\frac{\partial h}{\partial z_j} \frac{\partial h'}{\partial z_i}
%-\frac{3}{2\pi} \frac{\partial h}{\partial z_i}
%\frac{\partial h'}{\partial z_j}  ,\\
 %   &\left(\cL^{-1} \cP^{\bot} \ov{z}_{k}  \ov{z}_{l} F \cP\right)(0,0)
%    =\Big(\frac{b_k b_l}{32\pi^3}\phi
%+\frac{b_k}{8\pi^3}\tfrac{\partial F}{\partial z_l}
%+\frac{b_l}{8\pi^3}\tfrac{\partial F}{\partial z_k}\Big)\cP(0,0)
%=-\frac{3}{8\pi^{3}} \frac{\partial ^2 F}{\partial z_k\partial z_l}.
\end{align}
\end{subequations}
where we use in the last equation that $ h b_i  h' b_j= b_i h h' b_j + 2
h_i h' b_j$ and further \eqref{bk3.16c}, \eqref{bk3.16e}.

Let $\phi(z)= \phi_{ij} z_iz_j,\psi=\psi_{ij}
z_iz_j$ be degree $2$ polynomials in $z$ with symmetric matrices 
$(\phi_{ij})$, $(\psi_{ij})$. Then
by Theorem \ref{bkt2.17}, (\ref{bk2.66}), (\ref{bk3.10}), (\ref{bk3.16e}),
 (\ref{bk3.16f}) and (\ref{bk3.17a}), we obtain
\begin{subequations}
\begin{align}
    \left(\cL^{-1} b_{\ell }b_{k} \phi b_{j}h \cP\right) (0,0)
    &=  \left(\frac{b_{\ell }b_{k} b_{j}}{12\pi} \phi  h +
 \frac{b_{\ell }b_{k}}{4\pi} \frac{\partial \phi}{\partial z_{j}}
 h\right)  \cP (0,0)\label{bk3.18a}\\
& = \frac{-2}{3\pi}\frac{\partial^{3} (\phi h)}{\partial z_{\ell }\partial
z_{k}\partial z_{j}} + \frac{2}{\pi}(\phi_{\ell j}h_{k}
+\phi_{kj}h_{\ell }),\nonumber\\
     \left(\cL^{-1} b_{k} \phi b_{\ell } b_{j}h \cP\right) (0,0)
& =  \cL^{-1} \left(b_{\ell }b_{k}  \phi + 2 b_{k}
\frac{\partial \phi}{\partial z_{\ell }} \right)
b_{j}h \cP (0,0)\label{bk3.18b}\\	
   & = \frac{-2}{3\pi}\frac{\partial^{3} (\phi h)}{\partial z_{\ell }\partial
z_{k}\partial z_{j}} + \frac{2}{\pi}(\phi_{kj}h_{\ell }
+ \phi_{\ell k}h_{j}),\nonumber\\
- \pi \left(\cL^{-1} b_{k} \phi \ov{z}_{\ell } b_{j}h \cP\right) (0,0)
&= - \frac{1}{2} \cL^{-1}b_{k} \phi b_{j} \left(b_{\ell } h
+ 2 h_{\ell }  \right)\cP (0,0)\label{bk3.18c}\\
&= \frac{1}{3\pi}\frac{\partial^{3} (\phi h)}{\partial z_{\ell }\partial
z_{k}\partial z_{j}} - \frac{1}{\pi}\phi_{\ell k}h_{j},\nonumber
\end{align}
\end{subequations}
and
\begin{subequations}
\begin{align}
    \left( \cL^{-1}b_{k} h  b_{i}b_{j}\psi \cP\right) (0,0)
   & = \left(\frac{b_{k}b_{i} b_{j}}{12\pi} h \psi+
 \frac{b_{k}}{4\pi}(b_{i}h_{j}+ b_{j} h_{i})
\psi\right) \cP (0,0)\label{bk3.19a}\\
  &  = \frac{-2}{3\pi}\frac{\partial^{3} (h\psi)}
{\partial z_{k}\partial z_{i}\partial z_{j}}
 + \frac{2}{\pi}(\psi_{ki}h_{j} + \psi_{kj}h_{i}),\nonumber\\
    (\cL^{-1} b_{k} h \ov{z}_{\ell } b_{j} \psi\cP) (0,0)
    &=  \cL^{-1}b_{k} h b_{j} \left( \frac{b_{\ell }}{2\pi}\psi
+ \frac{1}{\pi}\frac{\partial \psi}{\partial z_{\ell }} \right)\cP (0,0)
\label{bk3.19b}\\
   & = \frac{-1}{3\pi^{2}}
\frac{\partial^{3} (h\psi)}{\partial z_{\ell }\partial z_{k}\partial z_{j}}
+ \frac{1}{\pi^{2}}(\psi_{kj}h_{\ell }+ \psi_{\ell j}h_{k}),\nonumber\\
- \pi  \left(\cL^{-1} h \ov{z}_{k}\ov{z}_{\ell } b_{j}\psi \cP\right) (0,0)
&= \frac{1}{6\pi^{2}}\frac{\partial^{3} (h \psi)}{\partial z_{\ell }\partial
z_{k}\partial z_{j}} + \frac{3}{2\pi^{2}}\psi _{\ell k}h_{j},\label{bk3.19c}
\end{align}
\end{subequations}
In fact, by \eqref{bk3.1}, \eqref{bk3.10} implies
\[
\pi (h \ov{z}_{k}\ov{z}_{\ell } b_{j}\psi \cP) (Z,0)
= \frac{1}{2}\Big(b_{k} h b_{j} \psi \ov{z}_{\ell }
+ \Big(2 h_{k}b_{j} \psi + 2 h b_j \frac{\partial \psi}{\partial
  z_{k}}\Big)  \frac{b_{\ell }}{2 \pi} \Big)\cP(Z,0)\,,
 \] 
so from \eqref{bk3.16e}, \eqref{bk3.17b} and  \eqref{bk3.19b}
we get \eqref{bk3.19c}.

Set $R_{\ov{m}\,\ov{q}}(z)= R_{k\ov{m}l\ov{q}}\, z_kz_l$.
By \eqref{bk2.66} and Lemma \ref{lmt1.6},
\[
R_{\ov{s}\ov{k}}(z) b_m b_qR_{\ov{m}\,\ov{q}}(z)
= \big(b_m b_q R_{\ov{s}\ov{k}}(z) + 8 b_m R_{q\ov{s}\ell\ov{k}}\, z_{\ell}
+ 8 R_{m\ov{s}q\ov{k}} \big) R_{\ov{m}\,\ov{q}}(z)\,,
\] 
thus Theorem \ref{bkt2.17}, \eqref{bk3.10},
\eqref{bk3.18b} and \eqref{bk3.19a} show that
\begin{equation}\label{bk3.21}
\begin{split}
 &     \left( \cL^{-1}b_{s}b_{k} R_{\ov{s}\ov{k}}
     b_{m}b_{q}R_{\ov{m}\,\ov{q}} \cP\right) (0,0)\\
 & \hspace{5mm}= b_{s}b_{k} \Big(\frac{b_m b_q}{16\pi} R_{\ov{s}\ov{k}}
+ \frac{2 b_m}{3\pi} R_{q\ov{s}\ell \ov{k}} z_{\ell }
+ \frac{1}{\pi} R_{m\ov{s}q\ov{k}}  \Big) R_{\ov{m}\,\ov{q}}\cP(0,0)\\
& \hspace{5mm}
= \frac{1}{\pi}\frac{\partial^{4} (R_{\ov{s}\ov{k}}R_{\ov{m}\,\ov{q}})}
{\partial z_{s}\partial z_{k}\partial z_{m}\partial z_{q}}
-  \frac{16}{3\pi}\frac{\partial^{3} (R_{q\ov{s}\ell \ov{k}}z_{\ell }
R_{\ov{m}\,\ov{q}})}
{\partial z_{s}\partial z_{k}\partial z_{m}}
+  \frac{8}{\pi}R_{m\ov{s}\,q\ov{k}} R_{s\ov{m}k\ov{q}}\, ,\\
& - \left( \pi \cL^{-1}b_{k} \ov{z}_{s} R_{\ov{s}\ov{k}}
     b_{m}b_{q}R_{\ov{m}\,\ov{q}} \cP\right) (0,0)\\
&   \hspace{5mm}
= - \cL^{-1}b_{k}  \Big( \frac{1}{2}b_{s} R_{\ov{s}\ov{k}}
     b_{m}b_{q}R_{\ov{m}\,\ov{q}}
+ \frac{\partial R_{\ov{s}\ov{k}}}{\partial z_{s}} b_m
  b_qR_{\ov{m}\,\ov{q}}
+   R_{\ov{s}\ov{k}}b_{m}b_{q}\frac{\partial
    R_{\ov{m}\,\ov{q}}}{\partial z_{s}} \Big) \cP(0,0)\\
& \hspace{5mm}
= \frac{1}{6\pi}\frac{\partial^{4} (R_{\ov{s}\ov{k}}R_{\ov{m}\,\ov{q}})}
{\partial z_{s}\partial z_{k}\partial z_{m}\partial z_{q}}
+  \frac{8}{3\pi}\frac{\partial^{3} (R_{q\ov{s}\ell \ov{k}}z_{\ell }
R_{\ov{m}\,\ov{q}})}
{\partial z_{s}\partial z_{k}\partial z_{m}}
-  \frac{4}{\pi}R_{m\ov{s}q\ov{k}} R_{s\ov{m}k\ov{q}}
- \frac{16}{\pi}R_{s\ov{s}q\ov{k}} R_{m\ov{m}k\ov{q}}\, .
\end{split}
\end{equation}
Due to \eqref{bk3.13}, \eqref{bk3.16a}--\eqref{bk3.21}, we get
 \begin{multline}\label{bk3.22}
 \pi^2   (\cL^{-1}\bI_1\cP)(0,0)=
\frac{1}{144} \Big[\frac{7}{6}
 \frac{\partial^{4} (R_{\ov{s}\ov{t}}R_{\ov{m}\ov{q}})}
{\partial z_{s}\partial z_{t}\partial z_{m}\partial z_{q}}
- \frac{8}{3}\frac{\partial^{3} (R_{q\ov{s}\ell \ov{t}}z_{v}
R_{\ov{m}\ov{q}})}
{\partial z_{s}\partial z_{t}\partial z_{m}}
+ 4 R_{m\ov{s}q\ov{t}} R_{s\ov{m}t\ov{q}}\\
 - 16 R_{s\ov{s}q\ov{t}} R_{m\ov{m}t\ov{q}}
 - 2 \frac{\partial^{3} ( ( 2R_{v\ov{u}u\ov{t}} + R^E_{v\ov{t}})z_{v}
R_{\ov{m}\ov{q}})}{\partial z_{t}\partial z_{m}\partial z_{q}}
+ 12 (2R_{q\ov{u}u\ov{t}} + R^E_{q\ov{t}}) R_{t\ov{m}m\ov{q}}\Big]\\
 + \frac{1}{12} \Big[ - \frac{1}{3}
\frac{\partial^{3} (  R_{\ov{s}\ov{t}}(\frac{4}{3}R_{\ell\ov{k}k\ov{q}}
+ R^E_{\ell\ov{q}})z_{\ell} )}
{\partial z_{s}\partial z_{t}\partial z_{q}}
+ 4 R_{s\ov{s}q\ov{t}}\big(\frac{4}{3}R_{t\ov{k}k\ov{q}}
+ R^E_{t\ov{q}}\big) - R_{s\ov{s}t\ov{t}}\big(\frac{4}{3}R_{q\ov{k}k\ov{q}}
+ R^E_{q\ov{q}}\big)\Big]\\
 +  \frac{1}{8} \Big[\big(2R_{t\ov{u}u\ov{t}}
+ R^E_{t\ov{t}}\big)\big(\frac{4}{3}R_{q\ov{k}k\ov{q}}
+ R^E_{q\ov{q}}\big)-\big(2R_{q\ov{u}u\ov{t}} + R^E_{q\ov{t}}\big)
\big(\frac{4}{3}R_{t\ov{k}k\ov{q}} + R^E_{t\ov{q}}\big)\Big].
\end{multline}
But from Lemma \ref{lmt1.6}, we have
\begin{equation}\label{bk3.23}
\begin{split}
& \frac{\partial^{4} (R_{\ov{s}\ov{t}}R_{\ov{m}\ov{q}})}
{\partial z_{s}\partial z_{t}\partial z_{m}\partial z_{q}}
= 4 R_{m\ov{s}q\ov{t}} R_{s\ov{m}t\ov{q}}
+ 4R_{s\ov{s}t\ov{t}} R_{m\ov{m}q\ov{q}}
+ 16 R_{s\ov{s}m\ov{t}} R_{q\ov{m}t\ov{q}},\\
&\frac{\partial^{3} (R_{s\ov{s}v\ov{t}} z_{v} R_{\ov{m}\ov{q}})}
{\partial z_{t}\partial z_{m}\partial z_{q}}
=  2R_{s\ov{s}t\ov{t}} R_{m\ov{m}q\ov{q}}
+ 4 R_{s\ov{s}q\ov{t}} R_{m\ov{m}t\ov{q}},\\
&\frac{\partial^{3} (R_{m\ov{s}v\ov{t}}z_{v} R_{\ov{m}\ov{q}})}
{\partial z_{s}\partial z_{t}\partial z_{q}}
= 2 R_{m\ov{s}q\ov{t}} R_{s\ov{m}t\ov{q}}
+ 4 R_{m\ov{s}s\ov{t}} R_{t\ov{m}q\ov{q}}.
\end{split}
\end{equation}
Plugging \eqref{bk3.23} in \eqref{bk3.22} we see that
the coefficient of $R_{s\ov{s}q\ov{t}} R_{m\ov{m}t\ov{q}}$ in
the term $\frac{1}{144}[\cdots]$ of \eqref{bk3.22} is
$\frac{1}{144}(  \frac{56}{3}- \frac{32}{3}-16 -16+24)=0$ and
\begin{equation}\label{bk3.24}
\begin{split}
 \pi^2 &(\cL^{-1}\bI_1\cP)(0,0)= \frac{1}{144}
%( \frac{14}{3}- \frac{16}{3}+4)
\frac{10}{3} R_{m\ov{s}q\ov{t}} R_{s\ov{m}t\ov{q}}
%+  \Big[\frac{1}{144}(  \frac{56}{3}- \frac{32}{3}-16 -16+24) +
%\frac{1}{12}( -\frac{16}{9} +\frac{16}{3}) -\frac{1}{3}\Big ]
-  \frac{1}{27} R_{s\ov{s}q\ov{t}} R_{m\ov{m}t\ov{q}}\\
&+   \Big[\frac{1}{144}\cdot\frac{-10}{3}+ \frac{1}{12}\cdot\frac{-20}{9}
%\frac{1}{144}(\frac{14}{3}-8)+ \frac{1}{12}(-\frac{8}{9}-\frac{4}{3})
+\frac{1}{3}\Big ]  R_{s\ov{s}t\ov{t}} R_{m\ov{m}q\ov{q}}
+ \Big[\frac{-4}{144}+ \frac{1}{12}\cdot\frac{-5}{3}%(-\frac{2}{3}-1)
+ \frac{5}{12} \Big]
R^E_{t\ov{t}}  R_{m\ov{m}q\ov{q}} \\
&+ \Big [\frac{4}{144} +  \frac{1}{12} \cdot\frac{8}{3}
%( - \frac{4}{3} +4)
-\frac{10}{3\cdot 8}  \Big ]  R^E_{q\ov{t}}  R_{m\ov{m}t\ov{q}}
+ \frac{1}{8} (R^E_{t\ov{t}} R^E_{q\ov{q}}
- R^E_{q\ov{t}} R^E_{t\ov{q}}) \\
= &\,  \frac{5}{2^3\cdot 3^3} R_{m\ov{s}q\ov{t}} R_{s\ov{m}t\ov{q}}
-  \frac{1}{27} R_{s\ov{s}q\ov{t}} R_{m\ov{m}t\ov{q}}
+  \frac{1}{8} R_{s\ov{s}t\ov{t}} R_{m\ov{m}q\ov{q}}   \\
&+ \frac{1}{4} R^E_{t\ov{t}}  R_{m\ov{m}q\ov{q}}
- \frac{1}{6} R^E_{q\ov{t}}  R_{m\ov{m}t\ov{q}}
+ \frac{1}{8} (R^E_{t\ov{t}} R^E_{q\ov{q}}
- R^E_{q\ov{t}} R^E_{t\ov{q}}).
\end{split}
\end{equation}
From  \eqref{bk3.13}, \eqref{bk3.16d}, \eqref{bk3.16g},
\eqref{bk3.16h}, \eqref{bk3.19b}, \eqref{bk3.19c} and \eqref{bk3.23}, we get
 \begin{equation}\label{bk3.25}
 \begin{split}
 \pi^2 (\cL^{-1} \bI_2 & \cP)(0,0)=
\frac{1}{18} \Big[ -\frac{1}{2}
\frac{\partial^{3} (R_{u\ov{s}q\ov{t}}z_{u}
R_{\ov{m}\ov{q}})}
{\partial z_{s}\partial z_{t}\partial z_{m}}
+ 4 R_{s\ov{s}q\ov{t}} R_{t\ov{m}m\ov{q}}\\
&+ \frac{3}{2} R_{m\ov{s}q\ov{t}} R_{s\ov{m}t\ov{q}}
+ 3 (2 R_{u\ov{u}q\ov{s}}
+ R^E_{q\ov{s}}) (-\frac{2}{4}) R_{m\ov{m}s\ov{q}} \Big] \\
&+ \Big[ %(-\frac{1}{3} +\frac{1}{4})
-\frac{1}{12}R_{t\ov{s}q\ov{t}}
- \frac{1}{4}( 2 R_{u\ov{u}q\ov{s}}
+ R^E_{q\ov{s}}) \Big] \big(\frac{4}{3}R_{s\ov{k}k\ov{q}}
+ R^E_{s\ov{q}}\big) \\
=&\frac{1}{36} R_{m\ov{s}q\ov{t}}
R_{s\ov{m}t\ov{q}}
-  \frac{5}{6} R_{s\ov{s}q\ov{t}} R_{m\ov{m}t\ov{q}}
- R^E_{s\ov{q}}  R_{t\ov{s}q\ov{t}}
- \frac{1}{4}  R^E_{q\ov{t}} R^E_{t\ov{q}}.
\end{split}
\end{equation}
By \eqref{bk3.16g}, we get
 \begin{align}\label{bk3.26}
  \Big(\cL^{-1}\frac{2\pi}{9} R_{m\ov{s}q\ov{t}}\ov{z}_{s}\ov{z}_{t}
 R_{k\ov{m}\ell\ov{q}} z_{k} z_{\ell} \cP\Big) (0,0)
= - \frac{1}{6\pi^2} R_{m\ov{s}q\ov{t}}R_{s\ov{m}t\ov{q}}.
\end{align}
{}Relations \eqref{bk3.14}, \eqref{bk3.24}, \eqref{bk3.25} and \eqref{bk3.26}
imply the desired formula \eqref{bk3.12}.
\end{proof}

%%%%%%%%%%%%%%%%%%%%%%%%%%%%%%%%%%%%%%%%%%
\subsection{Contribution of
$\mO_{4}$ to  $\cF_4(0,0)$}\label{bks3.3}%\label{kernel calculus}
%%%%%%%%%%%%%%%%%%%%%%%%%%%%%%%%%%%%%%%%%%%%

We will use the following remark repeatedly in our computation.

\begin{rem} \label{bkt3.5}
Let $\Phi$ be a polynomial in $b^+,z, b, \ov{z}$. Due to
\eqref{bk2.66} and \eqref{bk3.1}, the value of the kernels of
$\cP \Phi \cP$, $\cP^\bot \Phi \cP$,
$\cL^{-1}\cP^\bot \Phi \cP$ at $(0,0)$ consists of the terms of $\Phi$ whose
total degree in $b$ and $\ov{z}$ is the same as the total degree in 
$b^+$ and $z$.
\end{rem}

\begin{lemma}\label{bkt3.4} We have the following identity\,\rm{:}
\begin{align}\label{lm3.31}
\begin{split}
-\pi^2 (\cL^{-1} \mO_4 \cP)(0,0)= &\,-\frac{\Delta \br}{96}
%\frac{2}{45} R_{i\ov{i} j\ov{j}; j'\ov{j'}}
%+\frac{13}{45} R_{j'\ov{i} j\ov{j};  i\ov{j'}}
+ \frac{23}{108}
R_{m\ov{s} q\ov{t}}  R_{s\ov{m}t\ov{q}}
+ \frac{41}{54}
R_{s\ov{s} q\ov{t}}  R_{m\ov{m}t\ov{q}}\\
&+ R_{m\ov{m}q\ov{k}}  R^E_{k\ov{q}}
+ \frac{1}{8}(-R^E_{m\ov{m}; q\ov{q}}+3 R^E_{q\ov{m}; m\ov{q}})
+ \frac{1}{4}  R^E_{k\ov{q}} R^E_{q\ov{k}}\,.
\end{split}
\end{align}
\end{lemma}

\begin{proof} By (\ref{bk2.66}), \eqref{lm01.9}, (\ref{bk3.1}) and
    (\ref{bk3.2}), as in (\ref{bk3.3}), we have
\begin{align}\label{lm3.33}
\begin{split}
-&(\mO_{41} \cP)(Z,0)= - \Big\{  \frac{1}{20}\Big(A_{1Z}
-  \frac{4}{3}A_{2Z}\Big) (\tfrac{\partial}{\partial \ov{z}_i},
\tfrac{\partial}{\partial \ov{z}_j}) b_i b_j \\
&\hspace{45mm}- \frac{1}{20}\Big(A_{1Z}
-  \frac{4}{3}A_{2Z}\Big) (\tfrac{\partial}{\partial z_i},
\tfrac{\partial}{\partial \ov{z}_j})b^+_i b_j\Big\}  \cP(Z,0)\\
&=   \Big[- \frac{\pi^2}{5}\Big(A_{1Z}
-  \frac{4}{3}A_{2Z}\Big)  (\ov{z},\ov{z})
+  \frac{\pi}{5}\Big(A_{1Z}
-  \frac{4}{3}A_{2Z}\Big)   (\tfrac{\partial}{\partial z_j},
\tfrac{\partial}{\partial \ov{z}_j})\Big] \cP(Z,0).
\end{split}
\end{align}
and
\begin{align} \label{lm3.34}
-\cL^{-1} \mO_{42} \cP= &   \cP^\bot\Big[4 \Big(\frac{1}{80}A_{1Z}
 -  \frac{1}{360}A_{2Z}\Big)(\tfrac{\partial}{\partial z_j},
\tfrac{\partial}{\partial \ov{z}_j})
- \frac{1}{72} \ric(z, \ov{z})^{2} \Big]\cP.
\end{align}
From \eqref{lm01.7}, \eqref{lm01.9}, \eqref{bk3.1} and
    \eqref{bk3.2}, and since $R^{TX}$ is $(1,1)$-form, we have
    \begin{align} \label{lm3.35}
    \begin{split}
   - (\mO_{44} \cP)(Z,&0)
    = 2\pi \Big \{ \frac{\pi}{30}A_{1Z}(\ov{z}, \ov{z})
    +   \Big(\frac{\partial}{\partial  Z_j} \Big(\frac{1}{20}A_{1Z}
    +  \frac{2}{45}A_{2Z} \Big)\Big)( \ov{z}, e_j) \\
    &- \frac{\pi}{10}A_{2Z}(\ov{z}, \ov{z})-   \ov{z}_i 
    \frac{\partial}{\partial \ov{z}_i}  \Big(\frac{1}{10}A_{1Z}
     +  \frac{4}{45}A_{2Z} \Big)(\tfrac{\partial}{\partial z_j},
    \tfrac{\partial}{\partial \ov{z}_j})\Big\}\cP(Z,0),\\
  - (\mO_{45} \cP)(Z,&0) = 2\pi \Big \{
\frac{2}{9}\left \langle R^{TX} (\mR,e_k) z,  R^{TX}_{x_0}
      (\mR,e_k) \ov{z}\right\rangle\\
& - \frac{1}{9} \left \langle R^{TX} (z, \ov{z}) \mR,
    e_k\right\rangle \ric(\mR,e_k)\\
    &+\frac{1}{4}
    \left \langle R^{TX} (z, \ov{z}) \mR, e_j\right \rangle
    R^E (\mR,e_j)
    -  \frac{1}{4}R^E_{\, ;(Z,Z)} (z,
     \ov{z})\Big\}\cP(Z,0).
    \end{split}\end{align}
Let $\psi_{ijk}$ be degree $3$ polynomials in $z$  
 which are symmetric in $i,j$. By \eqref{bk2.66}, \eqref{bk3.1}, we get
\begin{align} \label{lm3.37}
\begin{split}
(&\psi_{ijk} \ov{z}_i \ov{z}_j  \ov{z}_k\cP)(Z,0)
= \frac{1}{8\pi^3} ( \psi_{ijk}  b_i b_j b_k \cP)(Z,0)
= \frac{1}{8\pi^3} \Big\{ b_i b_j b_k \psi_{ijk}  \\
& + 2 b_i b_j  \frac{\partial\psi_{ijk}}{\partial z_k}
+ 4 b_i b_k  \frac{\partial\psi_{ijk}}{\partial z_j}
+ 8 b_i  \frac{\partial^2\psi_{ijk}}{\partial z_j\partial z_k}
+ 4 b_k  \frac{\partial^2\psi_{ijk}}{\partial z_i\partial z_j}
+ 8 \frac{\partial^3\psi_{ijk}}
{\partial z_i\partial z_j\partial z_k}\Big\}\cP(Z,0).
\end{split}
\end{align}
Thus by Theorem \ref{bkt2.17}, (\ref{toe1.1}) and (\ref{lm3.37}), we obtain
\begin{align} \label{lm3.38}
\begin{split}
\pi^2 (\cL^{-1} & \cP^\bot \psi_{ijk} \ov{z}_i \ov{z}_j  \ov{z}_k
\cP)(0,0)
= \frac{1}{8\pi^2} \Big\{  \frac{b_i b_j b_k}{12} \psi_{ijk}
+ \frac{1}{4} b_i b_j  \frac{\partial\psi_{ijk}}{\partial z_k} \\
&+ \frac{1}{2} b_i b_k  \frac{\partial\psi_{ijk}}{\partial z_j}
+ 2 b_i  \frac{\partial^2\psi_{ijk}}{\partial z_j\partial z_k}
+  b_k  \frac{\partial^2\psi_{ijk}}{\partial z_i\partial z_j}\Big\}\cP(0,0)
= - \frac{11}{24\pi^2}\frac{\partial^3\psi_{ijk}}
{\partial z_i\partial z_j\partial z_k} \;.
\end{split}
\end{align}
Since $|\tfrac{\partial}{\partial z_j} |^2 = \frac{1}{2}$\,, 
 Lemma \ref{lmt1.6}, \eqref{lm01.4},  \eqref{lm01.7} and 
 the fact that $R^{TX}$ is a $(1,1)$-form entail
\begin{align} \label{lm3.39}
\begin{split}
&A_{1Z}(\ov{z},\ov{z})=\left\langle R^{TX}_{\, ; (Z,Z)}
(z,\ov{z}) z,\ov{z}\right\rangle,\\
&A_{2Z}(\ov{z},\ov{z})=2 \left \langle R^{TX} (z,\ov{z}) z,
R^{TX} (z,\ov{z}) \ov{z}\right \rangle
= -4 R_{u\ov{s}v\ov{t}} R_{k\ov{m} t\ov{q}}
  z_{u}  \ov{z}_{s}  z_{v}  z_{k}\ov{z}_{m}\ov{z}_{q},\\
&A_{1Z}(\tfrac{\partial}{\partial z_q},\tfrac{\partial}{\partial\ov{z}_q})
=\left\langle R^{TX}_{\, ;(Z,Z)} (\ov{z},
\tfrac{\partial}{\partial z_q}) z,
\tfrac{\partial}{\partial\ov{z}_q}\right\rangle,\\
&A_{2Z}(\tfrac{\partial}{\partial z_q},\tfrac{\partial}{\partial\ov{z}_q})
= \left \langle R^{TX} (\ov{z},
\tfrac{\partial}{\partial z_q}) \mR,
R^{TX} ( z,\tfrac{\partial}{\partial\ov{z}_q})\mR\right \rangle\\
&\hspace{23mm} = 2(R_{q\ov{s}k\ov{t}} R_{\ell\ov{q} t\ov{u}}
+ R_{q\ov{s}t \ov{u}} R_{\ell\ov{q} k\ov{t}} )
  \ov{z}_{s}\ov{z}_{u}   z_{k} z_{\ell}.
\end{split}\end{align}
By Remark \ref{bkt3.5}, we can replace $A_{1Z}(\ov{z},\ov{z})$
by $2R_{k\ov{m} \ell\ov{q}; s\ov{t}}  z_{k}  \ov{z}_{m}  z_{\ell}
\ov{z}_{q} z_{s}\ov{z}_{t}$ in our computation.
We deduce from \eqref{lm3.38} and \eqref{lm3.39} that
\begin{align} \label{lm3.40}
\begin{split}
\pi^2 (\cL^{-1}\cP^\bot A_{1Z}(\ov{z},\ov{z}) \cP)(0,0)
=& - \frac{11}{6\pi^2} (R_{m\ov{m} q\ov{q}; t\ov{t}}
+ 2  R_{t\ov{m} q\ov{q};  m\ov{t}}),\\
\pi^2 (\cL^{-1}\cP^\bot A_{2Z}(\ov{z},\ov{z}) \cP)(0,0)
=&  \frac{11}{3\pi^2} ( R_{m\ov{s} q\ov{t}}  R_{s\ov{m} t\ov{q}}
+ 2  R_{s\ov{s} q\ov{t}}  R_{m\ov{m} t\ov{q}}).
\end{split}\end{align}
Let $F_{ij}$
be homogeneous degree $2$ polynomials in $Z$. Then by
\eqref{toe1.1}, \eqref{bk3.10} and  \eqref{bk3.16g}, we obtain
\begin{subequations}
    \begin{align}
	&(\cP^\bot F_{ij}\ov{z}_i \ov{z}_j \cP)(0,0)
= -  \frac{1}{\pi^2}\frac{\partial^2 F_{ij}}{\partial z_i \partial z_j}\,,
  \label{lm3.36a}\\
&(\cL^{-1}\cP^\bot \ov{z}_k \frac{\partial (F_{ij}\ov{z}_i \ov{z}_j)}
{\partial \ov{z}_k}\cP)(0,0)
= -  \frac{3}{4\pi^3}\frac{\partial^2 F_{ij}}{\partial z_i \partial z_j}\,.
  \label{lm3.36d}
\end{align}\end{subequations}
(Note that by Remark \ref{bkt3.5} the contributions of
$\ov{z}_k \frac{\partial (F_{ij}\ov{z}_i \ov{z}_j)}{\partial \ov{z}_k}$
and  $2 F_{ij}\ov{z}_i \ov{z}_j $
to \eqref{lm3.36d} are the same, so \eqref{lm3.36d}
follows from \eqref{bk3.16g}.)

By Remark \ref{bkt3.5} and \eqref{lm3.39}, only the term
$-2 R_{q\ov{m} k\ov{q}; \ell\ov{t}}\ov{z}_{m}z_{k}\ov{z}_{t}z_{\ell}$
from $A_{1Z}(\tfrac{\partial}{\partial z_q},
\tfrac{\partial}{\partial\ov{z}_q})$
has a nontrivial contribution in our computation at $(0,0)$,
and from \eqref{bk3.16g}, \eqref{lm3.36a} and \eqref{lm3.36d}, we get
\begin{align} \label{lm3.45}
\begin{split}
&\left(\cP^{\bot}A_{1Z}(\tfrac{\partial}{\partial z_q},
   \tfrac{\partial}{\partial \ov{z}_q})\cP\right) (0,0)
   =\frac{2}{\pi ^2} (R_{q\ov{m} m\ov{q}; t\ov{t}}
+  R_{q\ov{m} t\ov{q}; m\ov{t}}),\\
&\frac{1}{2} \left(\cL^{-1}\cP^{\bot} \ov{z}_m
\frac{\partial}{\partial\ov{z}_m} A_{1Z}(\tfrac{\partial}{\partial z_q},
\tfrac{\partial}{\partial \ov{z}_q})\cP\right) (0,0) =
\left(\cL^{-1} \cP^{\bot}A_{1Z}(\tfrac{\partial}{\partial z_q},
\tfrac{\partial}{\partial \ov{z}_q})\cP\right) (0,0)\\
&\hspace{20mm}=  \frac{3}{4\pi ^3} (R_{q\ov{m} m\ov{q}; t\ov{t}}
+  R_{q\ov{m} t\ov{q}; m\ov{t}}),\\
&\left(\cP^{\bot}A_{2Z}(\tfrac{\partial}{\partial z_q},
\tfrac{\partial}{\partial \ov{z}_q})\cP\right) (0,0)
= - \frac{2}{\pi ^2} (R_{q\ov{s}s\ov{t}} R_{u\ov{q} t\ov{u}}
+ 3 R_{q\ov{s}u\ov{t}} R_{s\ov{q} t\ov{u}} ),\\
& \frac{1}{2} \left(\cL^{-1} \cP^{\bot}\ov{z}_m
\frac{\partial}{\partial\ov{z}_m }
A_{2Z}(\tfrac{\partial}{\partial z_q},
\tfrac{\partial}{\partial \ov{z}_q})\cP\right) (0,0)
=\left(\cL^{-1} \cP^{\bot}A_{2Z}(\tfrac{\partial}{\partial z_q},
\tfrac{\partial}{\partial \ov{z}_q})\cP\right) (0,0)\\
&\hspace{20mm}
= - \frac{3}{4\pi ^3} (R_{q\ov{s}s\ov{t}} R_{u\ov{q} t\ov{u}}
+ 3 R_{q\ov{s}u\ov{t}} R_{s\ov{q} t\ov{u}} ).
\end{split}\end{align}
\comment{
From (\ref{lm3.36d}) and (\ref{lm3.39}), we get
\begin{align} \label{lm3.46}
\begin{split}
    \left(\cL^{-1}\cP^{\bot} \ov{z}_i \tfrac{\partial}{\partial\ov{z}_j}
    A_{1Z}(\tfrac{\partial}{\partial z_i},
    \tfrac{\partial}{\partial \ov{z}_j})\cP\right) &(0,0) =
    \frac{3}{4\pi ^2} (R_{j\ov{j} i\ov{i}; l\ov{l}} + R_{j\ov{j}
    l\ov{i}; i\ov{l}} + R_{j\ov{l} i\ov{i}; l\ov{j}}+  R_{j\ov{l}
    l\ov{i}; \ov{j}} )\\
&=  \frac{3}{4\pi ^2} (R_{j\ov{j} i\ov{i}; l\ov{l}} + 3 R_{j\ov{j}
    l\ov{i}; i\ov{l}}),\\
\left(\cL^{-1} \cP^{\bot}\ov{z}_i \tfrac{\partial}{\partial\ov{z}_j}
A_{2Z}(\tfrac{\partial}{\partial z_i},
\tfrac{\partial}{\partial \ov{z}_j})\cP\right) &(0,0)
= - \frac{3}{4\pi ^2} (  R_{j\ov{l}i\ov{j'}} R_{l\ov{i} j'\ov{j}}
+ 3 R_{j\ov{j}i\ov{j'}} R_{l\ov{i} j'\ov{l}}  + 4
R_{j\ov{l}j'\ov{j}} R_{i\ov{i} l\ov{j'}} )\\
&= - \frac{3}{4\pi ^2} (  R_{j\ov{l}i\ov{j'}} R_{l\ov{i} j'\ov{j}}
+ 7 R_{j\ov{j}i\ov{j'}} R_{l\ov{i} j'\ov{l}}  ).
\end{split}\end{align}
}
Remark \ref{bkt3.5}, \eqref{bk3.2}, \eqref{bk3.16d}
and \eqref{bk3.16g} yield
\begin{align} \label{lm3.47}
\begin{split}
&\left(\cL^{-1}\cP^{\bot}  \left \langle  R^{TX}_{\, ;(Z,e_m)} (z,\ov{z})\mR,
e_m\right \rangle_{x_0}\cP\right) (0,0) \\
 &\hspace{10mm} = 2 \left(\cL^{-1}\cP^{\bot}  (R_{k\ov{q} \ell\ov{m}; m\ov{s}}
- R_{k\ov{q} m\ov{s}\,;\, \ell\ov{m}})z_{k} \ov{z}_{q} 	
z_{\ell}\ov{z}_{s}\cP\right) (0,0)=0,\\
&    \left(\cL^{-1}\cP^{\bot}   \left \langle R^{TX}_{x_0} (z, \ov{z}) \mR,
e_k\right\rangle_{x_0} \ric(\mR,e_k)\cP\right) (0,0) \\
&\hspace{10mm}= 2 \left(\cL^{-1}\cP^{\bot}  
(-R_{\ell\ov{s}k\ov{q}} \ric_{m\ov{k}}
+ R_{\ell\ov{s}m\ov{k}} \ric_{k\ov{q}}) z_{\ell} \ov{z}_{s} 	
\ov{z}_{q}z_{m}\cP\right) (0,0)=0,\\
&\left(\cL^{-1}\cP^{\bot} \ric(\mR,e_q)R^E_{x_0} (\mR,e_q) \cP\right)
(0,0) \\
&\hspace{10mm}=2 \left(\cL^{-1}\cP^{\bot} ( - \ric_{k\ov{q}} R^E_{q\ov{s}}
+ \ric_{q\ov{s}} R^E_{k\ov{q}}) z_{k} \ov{z}_{s} \cP\right) (0,0) = 0,\\
&\left(\cL^{-1}\cP^{\bot} R^E_{\, ;(Z,e_s)} (\mR,  e_s)
\cP\right) (0,0) \\
&\hspace{10mm}
= 2\left(\cL^{-1}\cP^{\bot} (R^E_{k\ov{s}; s\ov{q}}
- R^E_{s\ov{q}; k\ov{s}})z_{k} \ov{z}_{q}\cP\right) (0,0) =0.
\end{split}\end{align}
By \eqref{0c47}, \eqref{0c50}, (\ref{bk3.2}) and \eqref{lm3.47},
we know that for $\alpha=1,2$,
\begin{align} \label{lm3.48}\begin{split}
&\Big(\cL^{-1}\cP^{\bot} \frac{\partial}{\partial  Z_q}
\Big(A_{\alpha Z}( \ov{z}, e_q) \Big)\cP\Big) (0,0)=0,\\
&\Big(\cL^{-1}\cP^{\bot} \Big(\frac{\partial}{\partial  Z_q}
A_{\alpha Z}\Big)( \ov{z}, e_q) \cP\Big) (0,0)
= -2  \Big(\cL^{-1}\cP^{\bot}
A_{\alpha Z}(\tfrac{\partial}{\partial  \ov{z}_q},
 \tfrac{\partial}{\partial z_q} ) \cP\Big) (0,0).
\end{split}\end{align}
By \eqref{lm01.3} and  \eqref{bk3.16g}, we have
\begin{align} \label{lm3.49}
\begin{split}
&\left(\cL^{-1}\cP^{\bot} \left \langle R^{TX}_{x_0} (\mR,e_k) z,  R^{TX}_{x_0}
  (\mR,e_k) \ov{z}\right\rangle_{x_0}\cP\right) (0,0) \\
&\hspace{5mm}= 4 \left(\cL^{-1}\cP^{\bot}
 (R_{k\ov{s} \ell \ov{u}} R_{m\ov{k}u\ov{q}}
+ R_{m\ov{k} \ell \ov{u}}R_{k\ov{s} u\ov{q}}) \ov{z}_{s} 	z_{\ell }
z_{m} \ov{z}_{q}\cP\right) (0,0) \\
&\hspace{5mm}= -  \frac{3}{2\pi ^3} (R_{k\ov{s}s\ov{u}}  R_{q\ov{k}u\ov{q}}
+ 3 R_{s\ov{k}q\ov{u}}  R_{k\ov{s}u\ov{q}}),\\
&    \left(\cL^{-1}\cP^{\bot}   \left\langle \pi
R^{TX}_{x_0}(z,\ov{z})\mR , e_k\right\rangle
 R^E_{x_0} (\mR,e_k) \cP\right) (0,0) \\
&\hspace{5mm}= 2\pi \left(\cL^{-1}\cP^{\bot}
(-R_{\ell \ov{s}k\ov{q}} R^{E}_{m\ov{k}}
- R_{\ell \ov{s}m\ov{k}} R^{E}_{k\ov{q}}) z_{\ell } \ov{z}_{s} 	
\ov{z}_{q}z_{m}\cP\right) (0,0)
= \frac{3}{\pi ^2}R_{s\ov{s}q\ov{k}} R^{E}_{k\ov{q}},\\
&\left(\cL^{-1}\cP^{\bot}R^E_{\, ;(Z,Z)} (z,
     \ov{z})\cP\right) (0,0)
%= 2\left(\cL^{-1}\cP^{\bot}R^E_{k\ov{i}; l\ov{j}}
%z_k  \ov{z}_i z_l  \ov{z}_j\cP\right) (0,0)
=- \frac{3}{4\pi ^3} (R^E_{s\ov{s}\,;\, q\ov{q}}+ R^E_{q\ov{s}\,;\, s\ov{q}}) .
\end{split}\end{align}
Note that by Lemma \ref{lmt1.6}, $\ric(\mR,\mR)= 2
\ric_{k\ov{q}}z_{k}\ov{z}_{q}$, and by (\ref{toet1.1}), (\ref{bk2.66})
and (\ref{bk3.1}),
\begin{align}\label{alm3.49}
\cL \ric(\mR,\mR) \cP = 2 b_{m} b^+_{m} \ric_{k\ov{q}}z_{k}\ov{z}_{q} \cP
= 4 b_{m}\ric_{k\ov{m}}\,z_{k}\cP.
\end{align}
Thus by \eqref{lm01.9}, \eqref{bk3.17b} and \eqref{alm3.49}, we obtain
\begin{align} \label{alm3.50}
\begin{split}
    -(\cL^{-1} \cP^\bot \mO_{43} \cP)(0,0)
    &= \frac{1}{18} (\cL^{-1} \cP^\bot \ric_{l\ov{q}}z_{l} b_{m}
    \ric_{k\ov{m}}z_{k}\frac{b_{q}}{2\pi} \cP)(0,0)\\
 &   = - \frac{1}{72\pi ^2} (\ric_{m\ov{m}}\ric_{q\ov{q}}
+ 3 \ric_{m\ov{q}}\ric_{q\ov{m}}).
\end{split}\end{align}
From \eqref{lm3.33}--\eqref{lm3.48} and \eqref{alm3.50}, we get
\begin{align} \label{lm3.50}
\begin{split}
-\pi^2 (\cL^{-1}&(\mO_{41}+ \mO_{42}+ \mO_{43}+ \mO_{44}) \cP)(0,0)
=  -\frac{2}{15} \cdot \frac{-11}{6}(R_{\ell \ov{\ell } q\ov{q}\,;\, u\ov{u}}
+ 2  R_{u\ov{\ell } q\ov{q}\,;\,  \ell \ov{u}})\\
&+ \frac{1}{15} \cdot \frac{11}{3}
( R_{\ell \ov{s} q\ov{u}}  R_{s\ov{\ell} u\ov{q}}
+ 2  R_{s\ov{s} q\ov{u}}  R_{\ell\ov{\ell} u\ov{q}})  \\
&+ \Big( - \frac{2}{5}\cdot \frac{3}{4}
+  \frac{1}{10}\Big)   (R_{q\ov{m} m\ov{q}; u\ov{u}}
+  R_{q\ov{m} u\ov{q}; m\ov{u}})\\
& - \Big( \frac{-12}{15}\cdot \frac{3}{4}
%(\frac{-4}{15}- \frac{8}{45}- \frac{16}{45})\times \frac{3}{4}
-  \frac{2}{90}\Big)  (R_{q\ov{s}s\ov{u}} R_{v\ov{q} u\ov{v}}
+ 3 R_{q\ov{s}v\ov{u}} R_{s\ov{q} u\ov{v}} )\\
&+ \frac{1}{72} ( \ric_{\ell\ov{\ell}}\ric_{q\ov{q}}
+\ric_{\ell\ov{q}}\ric_{q\ov{\ell}}) 
- \frac{1}{72} ( \ric_{\ell\ov{\ell}}\ric_{q\ov{q}}
+ 3 \ric_{\ell\ov{q}}\ric_{q\ov{\ell}})\\
=& \frac{2}{45} R_{\ell\ov{\ell} q\ov{q}\,;\, u\ov{u}}
+\frac{13}{45} R_{u\ov{\ell} q\ov{q}\,;\,  \ell\ov{u}}
+ \frac{19}{9} R_{\ell\ov{s} q\ov{u}}  R_{s\ov{\ell}u\ov{q}}
+  R_{s\ov{s} q\ov{u}}  R_{\ell\ov{\ell}u\ov{q}}\;.
\end{split}\end{align}
Moreover, \eqref{lm3.35}, \eqref{lm3.47} and \eqref{lm3.49} yield
\begin{align} \label{lm3.51}
\begin{split}
-\pi^2 (\cL^{-1}\mO_{45}\cP)(0,0)
= &-\frac{2}{3} (R_{k\ov{\ell}\ell\ov{u}}  R_{q\ov{k}u\ov{q}}
+ 3 R_{\ell\ov{k}q\ov{u}}  R_{k\ov{\ell}u\ov{q}})\\
&+ \frac{3}{2}  R_{\ell\ov{\ell}q\ov{k}} R^{E}_{k\ov{q}}
+ \frac{3}{8}(R^E_{\ell\ov{\ell}; q\ov{q}}+ R^E_{q\ov{\ell}\,;\, \ell\ov{q}}) .
\end{split}\end{align}
Further, \eqref{lm01.9b}, \eqref{bk3.16d}, \eqref{lm3.40},
\eqref{lm3.47} and \eqref{lm3.49} imply
\begin{align} \label{lm3.52}
\begin{split}
&-\pi^2 (\cL^{-1}\mO_{46}\cP)(0,0)=\frac{11}{108}
( R_{\ell\ov{s} q\ov{u}}  R_{s\ov{\ell} u\ov{q}}
+ 2  R_{s\ov{s} q\ov{u}}  R_{\ell\ov{\ell} u\ov{q}})
+ \frac{1}{9} \ric_{k\ov{q}} \ric_{q\ov{k}}\\
&\hspace{45mm}-  \frac{1}{9} R_{k\ov{\ell}q\ov{k}} \ric_{\ell\ov{q}}
-  \frac{1}{2} R_{\ell\ov{\ell}q\ov{k}}  R^E_{k\ov{q}}
+\frac{1}{4}  R^E_{k\ov{q}} R^E_{q\ov{k}}
- \frac{1}{2}  R^E_{\ell\ov{\ell}; q\ov{q}}\\
%- \frac{1}{3} R_{k\ov{\ell}l\ov{k}}R^E_{\ell\ov{l}}\\
&\hspace{3mm}= \frac{11}{108} R_{\ell\ov{s} q\ov{u}}  R_{s\ov{\ell} u\ov{q}}
+ \frac{23}{54}R_{s\ov{s} q\ov{u}}  R_{\ell\ov{\ell}
  u\ov{q}}
- \frac{1}{2} R_{\ell\ov{\ell}q\ov{k}}  R^E_{k\ov{q}}
+ \frac{1}{4}  R^E_{k\ov{q}} R^E_{q\ov{k}}
 -\frac{1}{2}  R^E_{\ell\ov{\ell}\,;\, q\ov{q}}.
\end{split}\end{align}
\comment{
By \cite[\S 2.3.4, Prop. 6]{Petersen06},
\begin{align}\label{lm3.56}
d\br = 2div (\ric),
\end{align}
thus by Lemma \ref{lmt1.6} and \eqref{lm3.56}, we have
\begin{align}\label{lm3.57}
-\Delta\br =2 e_je_i(\ric(e_j,e_i)) = 32  R_{j'\ov{i} j\ov{j};
  i\ov{j'}} ,
\quad -\Delta\br =32  R_{i\ov{i} j\ov{j}; j'\ov{j'}}.
\end{align}
}
By \eqref{alm01.5},  \eqref{lm01.12},  \eqref{lm3.50}--\,\eqref{lm3.52},
we get \eqref{lm3.31}.
 The proof of Lemma \ref{bkt3.4} is completed.
\end{proof}

%%%%%%%%%%%%%%%%%%%%%%%%%%%%%%%%%%%%%%%%%%
\subsection{Evaluation of $\cF_4(0,0)$}\label{bks3.4}
%%%%%%%%%%%%%%%%%%%%%%%%%%%%%%%%%%%%%%%%
\comment{
Set
\begin{align}\label{bk2.6}\begin{split}
&\bb_{2\C}=
\frac{4}{45} R_{i\ov{i} j\ov{j}; j'\ov{j'}}
+\frac{26}{45} R_{j'\ov{i} j\ov{j};  i\ov{j'}}
+ \frac{1}{6}
R_{i\ov{i'} j\ov{j'}}  R_{i'\ov{i}j'\ov{j}}
- \frac{2}{3}
R_{i'\ov{i'} j\ov{j'}}  R_{i\ov{i}j'\ov{j}}
+  \frac{1}{2} R_{i'\ov{i'}j'\ov{j'}} R_{i\ov{i}j\ov{j}}  ,\\
&\bb_{2E}= R^E_{j'\ov{j'}}  R_{i\ov{i}j\ov{j}}
- R^E_{j\ov{j'}}  R_{i\ov{i}j'\ov{j}}
+ \frac{1}{2} (R^E_{j'\ov{j'}} R^E_{j\ov{j}}
- R^E_{j\ov{j'}} R^E_{j'\ov{j}})
+\frac{1}{4}(-R^E_{i\ov{i}; j\ov{j}}+3 R^E_{j\ov{i}; i\ov{j}}),\\
&
\bb_{2} = \frac{1}{\pi^2}(\bb_{2\C} + \bb_{2E}).
\end{split}\end{align}
}
  \begin{thm}\label{bkt3.0} The following identity holds\rm{:}
\begin{align}\label{bk2.8}
    \cF_{4,\,x_{0}}(0,0)= \bb_{2} .
\end{align}
\end{thm}
\begin{proof}
By Lemmas \ref{bkt3.2}, \ref{bkt3.4} and (\ref{bk2.34}), we have
\begin{multline}\label{lm3.54}
\pi^2 \cF_{41}(0,0) =- \frac{\Delta \br}{96}
%\frac{2}{45} R_{m\ov{m} q\ov{q}; u\ov{u}}
%+\frac{13}{45} R_{u\ov{m} q\ov{q};  m\ov{u}}
+ \frac{7}{72}
R_{m\ov{s} q\ov{u}}  R_{s\ov{m}u\ov{q}}
- \frac{1}{9}
R_{s\ov{s} q\ov{u}}  R_{m\ov{m}u\ov{q}}
+  \frac{1}{8} R_{s\ov{s}u\ov{u}} R_{m\ov{m}q\ov{q}}  \\
+ \frac{1}{4}R^E_{u\ov{u}}  R_{m\ov{m}q\ov{q}}
- \frac{1}{6} R^E_{q\ov{u}}  R_{m\ov{m}u\ov{q}}
+ \frac{1}{8} (R^E_{u\ov{u}} R^E_{q\ov{q}}
-  R^E_{q\ov{u}} R^E_{u\ov{q}}
-R^E_{m\ov{m}; q\ov{q}}+3 R^E_{q\ov{m}\,;\, m\ov{q}})\,.
\end{multline}
Remark \ref{toet2.7}, Lemmas \ref{lmt1.6}, \ref{bkt3.1},
(\ref{bk2.6}),   \eqref{bk2.32},
\eqref{lm3.54} and formula $\cP(0,0)=1$ entail
 \begin{equation}\label{lm3.55}
 \begin{split}
 J_{4,x_0}(0,0) = &\,\cF_{41}(0,0)+ \cF_{41}(0,0)^*
+\frac{1}{4\pi^2} \Big[\sum_{mq}R_{m\ov{m}q\ov{q}}
+ \sum_q R^E_{q\ov{q}}\Big]^2\\
&- \frac{1}{36\pi^2}  R_{m\ov{k}q\ov{\ell}}  R_{k\ov{m}\ell\ov{q}}
-  \frac{1}{4\pi^2} \Big(\frac{4}{3}R_{q\ov{v}v\ov{\ell}} +  R^E_{q\ov{\ell}}\Big)
\Big(\frac{4}{3} R_{\ell\ov{k}k\ov{q}} + R^E_{\ell\ov{q}}\Big)
=\bb_{2}\,.
\end{split}
\end{equation}
The proof of Theorem \ref{bkt3.0} is completed.
\end{proof}

%\newpage
%%%%%%%%%%%%%%%%%%%%%%%%%%%%%%%%%%%%%%%%%%
%%%%%%%%%%%%%%%%%%%%%%%%%%%%%%%%%%%%%%%%%
\section{The first  coefficients of  the asymptotic
expansion}\label{toes4}
%%%%%%%%%%%%%%%%%%%%%%%%%%%%%%%%%

The lay-out of this section is as follows. In Section \ref{toes4.1}, we
explain the expansion of the kernel of Berezin-Toeplitz operators
and verify its compatibility with Riemann-Roch-Hirzebruch Theorem.
In Section \ref{toes4.2}, we establish Theorem \ref{toet4.1}. The
results from Sections \ref{bks3.2}, \ref{bks3.3} play an important role here.
In Section \ref{toes4.3}, we prove Theorems  \ref{toet4.6}, \ref{toet4.5}, 
i.e., the expansion of the composition of two Berezin-Toeplitz operators.

We use the notations and assumptions from Introduction and Section \ref{toes3}.

%%%%%%%%%%%%%%%%%%%%%%%%%%%%%%%%%%%%%%%%%%
\subsection{Expansion of the kernel of Berezin-Toeplitz operators}
\label{toes4.1}
%%%%%%%%%%%%%%%%%%%%%%%%%%%%%%%%%%%%%

%Let $\nabla^{TX}$ be the connection on $\Lambda(T^*X)$ 
%and $\Lambda(T^*X)\otimes \End(E)$ induced by $\nabla^{TX}$ 
%and $\nabla^E$.
For $U\in TX$, we have (cf. \eqref{lm01.1})
\begin{align} \label{abk4.1}
\nabla^{T^{*}X}_U\wi{dz_j} = 2 \langle \nabla^{TX}_U
\wi{\tfrac{\partial}{\partial \ov{z}_j}},
\wi{\tfrac{\partial}{\partial z_m}}\rangle \wi{dz_m},\quad
\nabla^{T^{*}X}_U\wi{d\ov{z}_j} = 2 \langle \nabla^{TX}_U
\wi{\tfrac{\partial}{\partial z_j}},
\wi{\tfrac{\partial}{\partial \ov{z}_m}}\rangle \wi{d\ov{z}_m}.
\end{align}
For 
$\sigma=\sum_{kq} \sigma_{k\ov{q}}\wi{dz_k}\wedge\wi{d\ov{z}_q}
\in \Omega^{1,1}(X,\End(E))$, 
by \cite[Lemma 1.4.4]{MM07}, \eqref{abk2.5} and \eqref{abk4.1}, we get
\begin{align} \label{abk4.2}
 \begin{split}
\nabla^{1,0*} \sigma =& 
- \Big(2 \nabla^E_{\wi{\tfrac{\partial}{\partial \ov{z}_m}}} \sigma_{m\ov{q}}
+ 4 \sigma_{k\ov{q}} 
\Big\langle \nabla^{TX}_{\wi{\tfrac{\partial}{\partial \ov{z}_m}}}
\wi{\tfrac{\partial}{\partial \ov{z}_k}},
\wi{\tfrac{\partial}{\partial z_m}}\Big\rangle 
+ 4 \sigma_{m\ov{l}}
\Big\langle \nabla^{TX}_{\wi{\tfrac{\partial}{\partial \ov{z}_m}}}
\wi{\tfrac{\partial}{\partial z_l}},
\wi{\tfrac{\partial}{\partial \ov{z}_q}}\Big\rangle \Big) \wi{d\ov{z}_q},\\
\ov{\partial}^{E*} \sigma = &
 \Big(2 \nabla^E_{\wi{\tfrac{\partial}{\partial z_m}}} \sigma_{k\ov{m}}
+ 4 \sigma_{k\ov{q}} 
\Big\langle \nabla^{TX}_{\wi{\tfrac{\partial}{\partial z_m}}}
\wi{\tfrac{\partial}{\partial z_q}},
\wi{\tfrac{\partial}{\partial \ov{z}_m}}\Big\rangle 
+ 4 \sigma_{l\ov{m}}
\Big\langle \nabla^{TX}_{\wi{\tfrac{\partial}{\partial \ov{z}_m}}}
\wi{\tfrac{\partial}{\partial \ov{z}_l}},
\wi{\tfrac{\partial}{\partial z_k}}\Big\rangle \Big) \wi{d z_k}.
\end{split}\end{align}
%By using \eqref{0c39} for $r=0,1$ associated with
% the vector bundle $T^{(1,0)}X$, from \eqref{abk4.2},
% we get at $x_0$ as the point $0$ in $\R^{2n}$ under our identification, 
We evaluate now \eqref{abk4.2} at the point $x_0$ 
(identified to $0\in\R^{2n}$). By using \eqref{0c39} applied for $r=0,1$ 
associated with the vector bundles $E$, $T^{(1,0)}X$, we get
\begin{align} \label{abk4.3}
\begin{split}
&(\nabla^{1,0*} \sigma)_{x_0}  = 
-2 \frac{\partial\sigma_{m\ov{q}}}{\partial \ov{z}_m}(0)\,d\ov{z}_q,
  \quad   (\ov{\partial}^{E*} \sigma)_{x_0}
= 2 \frac{\partial \sigma_{k\ov{m}}}{\partial z_m}(0)\, d z_k,\, \,   \\
&(\ov{\partial}^{E*}\nabla^{1,0*} \sigma)_{x_0}
=4  \frac{\partial  \sigma_{k\ov{m}} }{\partial z_m\partial \ov{z}_k}(0)
+ 2 \left[R^E_{m\ov{k}}\,, \sigma_{k\ov{m}}(0)\right].
\end{split}\end{align}
Note that by %\cite[Lemma 1.4.4]{MM07},
\eqref{bk2.6}, \eqref{lm01.4} and \eqref{abk4.3}, we have at $x_0$, 
\begin{align} \label{bk4.2a}
\begin{split}
&\om= \frac{\sqrt{-1}}{2} dz_q\wedge d\ov{z}_q,\quad 
 \tr[R^{T^{(1,0)}X}]= \ric_{k\ov{q}} dz_k \wedge d\ov{z}_q
= -\sqrt{-1} \ric_\om,\\
& R^E=R^E_{k\ov{q}} dz_k\wedge  d\ov{z}_q,\quad
\nabla^{1,0*} R^E= -2 R^E_{m\ov{q}; \ov{m}} d\ov{z}_q,\quad
\ov{\partial}^{E*} R^E= 2  R^E_{k\ov{q};q} dz_k.
\end{split}\end{align}
For $f\in\cC^\infty(X,\End(E))$, recall that  $T_{f,\,p}(x,x')$
is the smooth kernel of the Berezin-Toeplitz operator  $T_{f,\,p}$
defined according to \eqref{toe2.4}.
Then by \eqref{0c39}, at $x_0$, 
\begin{align} \label{bk4.3a}
\frac{\partial ^2 f_{\,x_0}}{\partial z_{q}
\partial \ov{z}_{\ell}} (0) =
(\nabla^E_{\wi{\tfrac{\partial}{\partial z_q}}}
\nabla^E_{\wi{\tfrac{\partial}{\partial \ov{z}_\ell}}} f)(x_0)
 - \frac{1}{2}  [R^E_{q\ov{\ell}}, f(x_0)],\quad 
\Delta^E f = -4 \frac{\partial ^2 f_{\,x_0}}{\partial z_{q}
\partial \ov{z}_{q}} (0).
\end{align}
In view of Lemma \ref{lmt1.6}, \eqref{bk2.5} and \eqref{bk4.3a},
we introduce the following coefficients:
\begin{align}\label{bk4.1}\begin{split}
\bb_{\C f}:= & R_{m\ov{m} q\ov{q}} \frac{\partial ^2 f_{\,x_0}}{\partial z_{k}
\partial \ov{z}_{k}} (0)
- R_{\ell\ov{k}k\ov{q}}\frac{\partial ^2 f_{\,x_0}}{\partial z_{q}
\partial \ov{z}_{\ell}} (0)\\
=&- \frac{\br}{32} \Delta^E f   - \frac{\sqrt{-1}}{8} 
\big\langle \ric_\om, \nabla^{1,0}\ov{\partial}^Ef
-\frac{1}{2}  [R^E, f] \big\rangle_{\om}\,,
\end{split}\end{align}
%where $\langle\cdot\,,\cdot \rangle$ acts $\C$-bilinearly
%(and pointwisely) on $TX$ and but not on $E$, 
and
\begin{align}\label{bk4.1a}\begin{split}
\bb_{E f1}:=& \frac{\partial f_{\,x_0}}{\partial z_{u}}(0)
\Big( \frac{1}{6} R^E_{k\ov{k} \,;\,  \ov{u}}
-\frac{5}{12} R^E_{q\ov{u}\,;\,\ov{q}}\Big)
+ \frac{1}{4} R^E_{m\ov{u}\,;\, \ov{m}}
\frac{\partial f_{\,x_0}}{\partial z_{u}}(0) \\
&+ \Big( \frac{1}{6} R^E_{k\ov{k}\,;\, u}
- \frac{5}{12} R^E_{u\ov{q}\,;\,q}\Big)
\frac{\partial f_{\,x_0}}{\partial \ov{z}_{u}}(0)
+ \frac{1}{4}\frac{\partial f_{\,x_0}}{\partial \ov{z}_{u}}(0)
R^E_{u\ov{m} \,;\, m}\\
=& \frac{1}{48} \Big\langle  \nabla^{1,0} f, 2\sqrt{-1}\,
 \ov{\partial}^E R^E_{\Lambda}
+ 5 \nabla^{1,0*} R^E\Big\rangle_{\om}
- \frac{1}{16} \langle \nabla^{1,0*} R^E, \nabla^{1,0} f \rangle_{\om}\\
&+ \frac{1}{16} \langle  \ov{\partial}^E f, \ov{\partial}^{E*}R^E \rangle_{\om}
+  \frac{1}{48}\Big\langle  2\sqrt{-1}\nabla^{1,0}R^E_{\Lambda}
 - 5  \ov{\partial}^{E*}R^E, \ov{\partial}^E f\Big\rangle_{\om},\\
\bb_{E f2}:=&\frac{1}{2} \frac{\partial ^2 f_{\,x_0}}{\partial z_{k}
\partial \ov{z}_{k}} (0) R^E_{q\ov{q}}
- \frac{1}{2} \frac{\partial ^2 f_{\,x_0}}{\partial z_{q}
\partial \ov{z}_{\ell}} (0)  R^E_{\ell\ov{q}}
+ \frac{1}{2} R^E_{q\ov{q}}
\frac{\partial ^2 f_{\,x_0}}{\partial z_{k}\partial \ov{z}_{k}} (0)
- \frac{1}{2}R^E_{\ell\ov{q}}
\frac{\partial ^2 f_{\,x_0}}{\partial z_{q}\partial \ov{z}_{\ell}} (0)\\
= &- \frac{\sqrt{-1}}{16} \Big[R^E_{\Lambda} \Delta^E f
+ (\Delta^E f)R^E_{\Lambda}   \Big]
+\frac{1}{8}  \Big\langle \nabla^{1,0} \ov{\partial}^E f
-\frac{1}{2}  [R^E, f], R^E\Big\rangle_{\om}\\
&+\frac{1}{8} \Big \langle  R^E, \nabla^{1,0} \ov{\partial}^E f
-\frac{1}{2}  [R^E, f]\Big\rangle_{\om}.
\end{split}\end{align}

The following result implies Theorem \ref{toet4.1}.
\begin{thm} \label{toet4.1a}
Let $f\in\cC^\infty(X,\End(E))$.
There exist smooth sections
$\bb_{r,f}(x)\in \End(E)_x$ such that \eqref{bk4.2} and \eqref{bk4.3}
hold and
\begin{equation}\label{bk4.4}
\begin{split}
\pi^2 \bb_{2,f} =  \bb_{2\C} f(x_0)
&+ \frac{1}{2} \Big(\bb_{2E}+ \frac{1}{16}(R^E_{\Lambda})^2\Big) f(x_0)
+ \frac{1}{2} f(x_0)  \Big(\bb_{2E}+ \frac{1}{16}(R^E_{\Lambda})^2\Big)\\
&- \frac{1}{16}R^E_{\Lambda}  f(x_0)R^E_{\Lambda}
+ \frac{1}{32}(\Delta^{E})^2 f + \bb_{\C f}
+\bb_{Ef1}+\bb_{Ef2}.
\end{split}
\end{equation}
\end{thm}

Before giving the proof, we verify that Theorem \ref{toet4.1} is 
compatible with the
Riemann-Roch-Hirzebruch Theorem. Note that by \eqref{toe2.1},
the first Chern class $c_1(L)$ of $L$ is represented by $\om$.
By the Kodaira vanishing Theorem and the
Riemann-Roch-Hirzebruch Theorem, we have for $p$ large enough:
\begin{equation}\label{bk2.10}
\begin{split}
\dim  H^0 (&X,  L^p\otimes E)= \int_X \td(T^{(1,0)}X)\ch( E) \, e^{p\, \om}\\
=&\,{\rm rk} (E) \int_X \frac{\om^n}{n!} p^n
+ \int_X \Big(c_1(E) + \frac{{\rm rk} (E)}{2} c_1(X)\Big)
\frac{\om^{n-1} p^{n-1}}{(n-1)!} \\
&+  \int_X \Big({\rm rk} (E)\{\td(T^{(1,0)}X)\}^{(4)}
+ \frac{1}{2} c_1(X)c_1(E)
+\{\ch( E)\}^{(4)}\Big)
\frac{\om^{n-2} p^{n-2}}{(n-2)!}\\
 &+ \cO(p^{n-3})\,.
 \end{split}
\end{equation}
As usual, $\ch(\cdot), c_1(\cdot), \td(\cdot)$ are the Chern character,
the first Chern class and the Todd class of the corresponding complex
vector bundles, $\{\cdot\}^{(4)}$ is the degree $4$-part of the
corresponding differential forms.
Note that \[\frac{x}{1-e^{-x}}= 1+ \frac{x}{2} + \frac{x^2}{12}
+\ldots\,,\] thus $\{\td(T^{(1,0)}X)\}^{(4)}= \frac{1}{12} (c_1(X)^2 + c_2(X))$.
%We like to find the coefficients of $p^{n-2}$ by using the curvatures,
Let $R^{T^{(1,0)}X}$ be the curvature of the Chern connection
on $T^{(1,0)}X$ which is the restriction of the Levi-Civita connection
in our case.
Then by \eqref{bk4.2a},
 we have the following identities at the cohomology level:
\begin{align} \label{bk2.11}
\begin{split}
&\{\ch( E)\}^{(4)} = -\frac{1}{8\pi^2} \tr[(R^E)^2],
\qquad c_1(X)= \frac{1}{2\pi} \ric_\om,\\
%\frac{\sqrt{-1}}{2\pi}\tr[R^{T^{(1,0)}X}],\\
%=\frac{1}{2\pi} \ric_{k\ov{q}} dz_k d\ov{z}_q,\\
&\{\td(T^{(1,0)}X)\}^{(4)}
 = \frac{1}{32\pi^2} (\ric_\om)^2
%-\frac{1}{32\pi^2} \Big( \tr[R^{T^{(1,0)}X}]\Big)^2
+ \frac{1}{96\pi^2} \tr\Big[(R^{T^{(1,0)}X})^2 \Big].
\end{split}\end{align}
%Now, at $x_0$ as the point $0$ in $\R^{2n}$ under our identification,
%by \eqref{toe2.1} and \eqref{bk2.64}, we have
%$\om_{x_0}= \frac{\sqrt{-1}}{2} dz_q\wedge d\ov{z}_q$, and 
%$\sqrt{-1}\tr[R^{T^{(1,0)}X}]=\ric_{k\ov{q}} dz_k \wedge d\ov{z}_q$,
By Lemma \ref{lmt1.6} and \eqref{bk4.2a}, we have 
\begin{align} \label{bk2.14}
\begin{split}
&  \frac{1}{32}\left\langle (\ric_\om)^2, \om^2/2\right\rangle
=  \frac{1}{2} (R_{u\ov{u}v\ov{v}}
R_{k\ov{k}q\ov{q}}- R_{u\ov{u} q\ov{v}}  R_{k\ov{k}v\ov{q}}),\\
&\frac{1}{96}\left\langle \tr\Big[(R^{T^{(1,0)}X})^2 \Big],
\om^2/2\right\rangle
= \frac{1}{6}(R_{k\ov{u} q\ov{v}}  R_{u\ov{k}v\ov{q}}
- R_{u\ov{u} q\ov{v}}  R_{k\ov{k}v\ov{q}}),\\
& \langle \{\ch( E)\}^{(4)}, \om^2/2\rangle
= \frac{1}{2 \pi^2} \tr \Big[R^E_{k\ov{k}} R^E_{q\ov{q}}
- R^E_{k\ov{q}} R^E_{q\ov{k}}\Big],\\
& \left\langle \frac{1}{2} c_1(X) c_1(E), \om^2/2\right\rangle = \frac{1}{2\pi^2}
\tr (R^E_{k\ov{k}}  \ric_{q\ov{q}} - \ric_{q\ov{k}}  R^E_{k\ov{q}}).
\end{split}\end{align}

\comment{
Let $\ov{\partial}^*$, $\partial^*$ be the adjoints of the operators
$\ov{\partial}$, $\partial$.
Let $\nabla^{TX}$ be the connection on
$\Lambda(T^{*}X)$ induced by the Levi-Civita connection  $\nabla^{TX}$.
Finally, for $\alpha$ a $(1,1)$-form, we have
\begin{align} \label{bk2.15}
 (\ov{\partial}^* \partial^* \alpha)_{x_0}
= 4 \Big(i_{\wi{\frac{\partial}{\partial \ov{z}_{q}}}}\,
\nabla^{TX}_{\wi{\frac{\partial}{\partial z_{q}}}}\,
i_{\wi{\frac{\partial}{\partial z_{i}}}}\,
\nabla^{TX}_{\wi{\frac{\partial}{\partial \ov{z}_{i}}}}\alpha\Big)_{x_0}
=4  \frac{\partial}{\partial z_{q}}
\frac{\partial}{\partial \ov{z}_{i}}
\Big(\alpha\Big(\wi{\tfrac{\partial}{\partial z_{i}}},
\wi{\tfrac{\partial}{\partial \ov{z}_{q}}}\Big)\Big)_{x_0},
\end{align}
by using \eqref{0c39} for $r=1$ and the bundle $T^{(1,0)}X$.
%Thus we get (\ref{bk2.16}).
}

We set now $f=1$ in Theorem \ref{toet4.1}, take the pointwise trace 
of the expansion \eqref{bk4.2} relative to $E$ and then integrate 
the result over $X$ with respect to the volume form $\om^{n}/n!$\,.
Taking into account \eqref{bk2.6} and \eqref{bk2.11}, \eqref{bk2.14}, 
we recover the expansion
up to $\cO(p^{n-3})$ given in \eqref{bk2.10} for the Hilbert polynomial.
%From  \eqref{bk2.10}-\eqref{bk2.14},
%by taking $f=1$ in Theorem \ref{toet4.1}
%and taking the trace of (\ref{bk4.2}) on $E$ and integrating on $X$,
%we refind  \eqref{bk2.10}. 
Thus the value of $\bb_{2}$ obtained in Theorem \ref{toet4.1}
is compatible with the Riemann-Roch-Hirzebruch Theorem.

%%%%%%%%%%%%%%%%%%%%%%%%%%%%%%%%%%%%%%%%%%
\subsection{Proof of Theorem \ref{toet4.1a}}\label{toes4.2}
%%%%%%%%%%%%%%%%%%%%%%%%%%%%%%%%%%%%%%%%%%

The first part of Theorem \ref{toet4.1a} follows from Lemma \ref{toet2.3}.
Moreover, by \eqref{toe2.13}, we have for any $r\in \N$,
\begin{align}\label{bk4.5}
 \bb_{r,f}(x_0) = Q_{2r,x_0}(f) (0,0).
\end{align}
Thus by \eqref{toe2.15}, the formula $\bb_{0,f}=f$,
and by \eqref{toe2.14} and \eqref{bk2.33}, we get
\begin{equation} \label{bk4.6}
  Q_{2,\,x_0}(f) =
 \cK\big[1, f(x_0)J_{2,x_0} \big]
+ \cK\big[J_{2,x_0}, f(x_0)\big]
+ \sum_{|\alpha|=2}
    \cK\Big[1\;,\;
  \frac{\partial ^\alpha f_{\,x_0}}{\partial Z^\alpha}(0)
  \frac{Z^\alpha}{\alpha !} \Big]\,.
  \end{equation}
Further, \eqref{toe1.6}, \eqref{bk2.24} and \eqref{bk2.31},
entail
\begin{align}\label{bk4.7}
\begin{split}
&\cK\big[1, f(x_0)J_{2,x_0} \big]\cP
= - f(x_0)\cP \mO_2 \cL^{-1}\cP^{\bot},\\
&\cK\big[J_{2,x_0}, f(x_0) \big]\cP
= - ( \cL^{-1}\cP^{\bot}\mO_2\cP ) f(x_0).
\end{split}\end{align}
From \eqref{toe1.6} and \eqref{abk3.8}, we deduce
\begin{align}\label{bk4.8}\begin{split}
&\sum_{|\alpha|=2}
    \cK\Big[1\;,\;
  \frac{\partial ^\alpha f_{\,x_0}}{\partial Z^\alpha}(0)
  \frac{Z^\alpha}{\alpha !} \Big]\cP(Z,0)
=\bigg(\sum_{|\alpha|=2}\frac{\partial^\alpha f_{x_0}}{\partial {z}^\alpha}(0)
\frac{z^\alpha}{\alpha !} 
+ \frac{1}{\pi} \frac{\partial ^2 f_{x_0}}
{\partial z_{i} \partial \ov{z}_{i}}(0)\!\bigg)\cP(Z,0).
\end{split}\end{align}
Lemma \ref{lmt1.6}, \eqref{bk3.7}, \eqref{bk3.8}, \eqref{bk4.3a}
and \eqref{bk4.5}--\eqref{bk4.8} yield the formula for $\bb_{1,f}$ 
from \eqref{bk4.2}.

It remains to compute $\bb_{2,f}$ for 
a self-adjoint section $f\in \cC^\infty(X,\End(E))$ in order to complete 
the proof of Theorem  \ref{toet4.1a}. 
%By (\ref{toe1.6}), set
Set
 \begin{align}\label{abk4.9}
 \cK_{2f}=    \sum_{|\alpha|=2}
     \cK\Big[1, \frac{\partial ^\alpha f_{\,x_0}}{\partial Z^\alpha}(0)
       \frac{Z^\alpha}{\alpha !}J_{2,x_0} \Big].
\end{align}
By  \eqref{toe2.14} and \eqref{bk2.33}, we get
\begin{equation} \label{bk4.9}
\begin{split}
  Q&_{4,\,x_0}(f) = \cK\Big[1, f(x_0)J_{4,x_0} \Big]
+\cK\Big[J_{2,x_0}, f(x_0)J_{2,x_0} \Big]
+ \cK\Big[J_{4,x_0}, f(x_0)\Big]\\
&+ \sum_{|\alpha|=2}  \cK\Big[J_{2,x_0}\;,\;
  \frac{\partial ^\alpha f_{\,x_0}}{\partial Z^\alpha}(0)
  \frac{Z^\alpha}{\alpha !} \Big]+ \cK_{2f}\\
&+  \cK\Big[1, \frac{\partial f_{\,x_0}}{\partial Z_i}(0) Z_i J_{3,x_0} \Big]
+ \cK\Big[J_{3,x_0}, \frac{\partial f_{\,x_0}}{\partial Z_i}(0) Z_i \Big]
+\sum_{|\alpha|=4}
\cK\Big[1,   \frac{\partial ^\alpha f_{\,x_0}}{\partial Z^\alpha}(0)
  \tfrac{Z^\alpha}{\alpha !}\Big] \,.
  \end{split}
  \end{equation}
Since $\cL_0$ and $\mO_r$ are formally self-adjoint, 
\eqref{bk2.23} and \eqref{bk2.77} show that $(\cF_{r,x_0})^* =\cF_{r,x_0}$. 
Hence, in the right hand side of \eqref{bk4.9},
the first, fourth and sixth terms
are adjoints of the third, fifth and seventh terms, respectively.
When we take $f=1$ in \eqref{bk4.9}, we get
\begin{align} \label{bk4.10}
  J_{4,x_0}=  \cK\big[1, J_{4,x_0} \big]
+\cK\big[J_{2,x_0}, J_{2,x_0} \big]
+ \cK\big[J_{4,x_0}, 1\big],
\end{align}
which is also a direct consequence of \eqref{bk2.31},
\eqref{bk2.32} and \eqref{bk3.0a}, as by  \eqref{toe1.6},
\begin{align}\label{bk4.11}
    \begin{split}
&\cK\big[1, J_{4,x_0} \big] \cP=
\cP\mO_{2}\cL^{-1}\cP^\bot \mO_{2}\cL^{-1}\cP^\bot
- \cP \mO_{4}\cL^{-1}\cP^\bot
-  \cP \mO_{2} \cL^{-2}   \mO_{2} \cP,\\
&\cK\big[J_{4,x_0}, 1\big]\cP=  \cL^{-1}\cP^\bot
\mO_{2}\cL^{-1}\cP^\bot \mO_{2} \cP- \cL^{-1}\cP^\bot \mO_{4} \cP-
\cP \mO_{2} \cL^{-2}   \mO_{2} \cP,\\
&\cK\big[J_{2,x_0}, J_{2,x_0} \big]\cP
= \cL^{-1}\cP^\bot \mO_{2} \cP \mO_{2} \cL^{-1}\cP^\bot
+  \cP \mO_{2} \cL^{-2}   \mO_{2} \cP.
\end{split}\end{align}
Set
\begin{align}\label{abk4.10}\begin{split}
&K_{41}:= -\frac{\Delta \br}{96}
+ \frac{5}{72} R_{m\ov{u} q\ov{v}}  R_{u\ov{m}v\ov{q}}
- \frac{5}{9} R_{u\ov{u} q\ov{v}}  R_{m\ov{m}v\ov{q}}
+  \frac{1}{8} R_{u\ov{u}v\ov{v}} R_{m\ov{m}q\ov{q}}  ,\\
&K_{42}:= \frac{1}{4}R^E_{v\ov{v}}  R_{m\ov{m}q\ov{q}}
- \frac{5}{6} R^E_{q\ov{v}}  R_{m\ov{m}v\ov{q}}
+ \frac{1}{8} \Big(R^E_{v\ov{v}} R^E_{q\ov{q}}
- 3 R^E_{q\ov{v}} R^E_{v\ov{q}}
-R^E_{m\ov{m}; q\ov{q}}+3 R^E_{q\ov{m}; m\ov{q}}\Big),\\
&K_{2f}:=  \frac{1}{4} (R_{m\ov{m}q\ov{q}}
+ R^E_{q\ov{q}}) f(x_{0})(R_{u\ov{u}v\ov{v}} + R^E_{v\ov{v}})
+\frac{1}{36}  R_{m\ov{k}q\ov{\ell}}  R_{k\ov{m}\ell\ov{q}} f(x_{0})  \\
&\hspace{15mm}
+ \frac{1}{4} \Big(\frac{4}{3}R_{q\ov{s}s\ov{\ell}} +  R^E_{q\ov{\ell}}\Big)
f(x_{0}) \Big(\frac{4}{3} R_{\ell\ov{k}k\ov{q}} + R^E_{\ell\ov{q}}\Big).
\end{split}\end{align}
By \eqref{bk2.34}, (\ref{bk3.0c}), (\ref{lm3.54}) and (\ref{bk4.11}), we have
\begin{align}\label{abk4.11}
    &\cK\big[J_{4,x_0},1 \big] (0,0)= \frac{1}{\pi^2} (K_{41}+ K_{42}).
\end{align}
By \eqref{toe1.6} and (\ref{bk2.31}), we see as in (\ref{bk4.11}) that
\begin{equation}\label{abk4.12}
\begin{split}
    \cK\big[J_{2,x_0}, f(x_{0})J_{2,x_0} \big] \cP
    = \cL^{-1} \cP^\bot\mO_2  f(x_{0})\cP \mO_2\cL^{-1}\cP^\bot\\
+ \cP \mO_2\cL^{-1} f(x_{0})\cL^{-1} \mO_2 \cP.
\end{split}
\end{equation}
Thus by (\ref{bk3.7}), (\ref{bk3.8}) and (\ref{abk4.12}),
as in Lemma \ref{bkt3.1}, we get
\begin{align}\label{abk4.13}
    \cK\big[J_{2,x_0}, f(x_{0})J_{2,x_0} \big] (0,0)
    =\frac{1}{\pi^2}  K_{2f}.
\end{align}
We next compute the fifth term in \eqref{bk4.9}.
 From  \eqref{bk2.31} and (\ref{abk4.9}), we get
\begin{align}\label{bk4.12}
\cK_{2f}=  \cP\sum_{|\alpha|=2} \frac{\partial^\alpha f_{\,x_0}}
{\partial Z^\alpha}(0)
  \frac{Z^\alpha}{\alpha !} \Big(-\cL^{-1} \mO_2 \cP
-\cP \mO_2\cL^{-1}\cP^\bot\Big)\,.
   \end{align}
For a degree $2$ polynomial $F(Z)$ we have
by Remark \ref{bkt3.5}, \eqref{toe1.4}, \eqref{bk2.66}, \eqref{bk3.1}
 and \eqref{bk3.7}
\begin{equation}\label{bk4.13}
\begin{split}
- (& \cP F \cL^{-1}  \mO_2 \cP) (0,0)\\
&=- \Big(\cP \frac{\partial ^2 F}{\partial z_{u}
\partial \ov{z}_{v}} z_{u}  \ov{z}_{v}
\Big\{\frac{b_{m} b_{q}}{48\pi}  R_{k\ov{m}l\ov{q}} z_{k} z_{l}
    + \frac{b_q}{4\pi}\Big( \frac{4}{3} R_{l\ov{k}k\ov{q}}
+   R^E_{l\ov{q}}\Big)z_l\Big\}
 \cP\Big) (0,0)\\
&= -  \Big(\cP \frac{\partial ^2 F}{\partial z_{q}
\partial \ov{z}_{v}}  \frac{ \ov{z}_{v}}{2\pi}
\Big( \frac{4}{3} R_{l\ov{k}k\ov{q}}+ R^E_{l\ov{q}}\Big)z_l
 \cP\Big) (0,0)
= - \frac{1}{2\pi^2} \frac{\partial ^2 F}{\partial z_{q}
\partial \ov{z}_{l}}\Big( \frac{4}{3} R_{l\ov{k}k\ov{q}}
+   R^E_{l\ov{q}}\Big),
\end{split}
  \end{equation}
where we have used \eqref{toe1.4}, \eqref{bk2.66}
and \eqref{abk3.8} in the last two equalities.

By  \eqref{bk3.8}, \eqref{abk3.8},  \eqref{bk4.12} and  \eqref{bk4.13},
we get
\begin{align}\label{bk4.14}
%\sum_{|\alpha|=2} \cK\Big[1,
%\tfrac{\partial ^\alpha f_{\,x_0}}{\partial Z^\alpha}(0)
%  \tfrac{Z^\alpha}{\alpha !}J_{2,x_0} \Big](0,0)\\
\cK_{2f}(0,0)= \frac{1}{2\pi^2} \frac{\partial ^2 f_{\,x_0}}{\partial z_{k}
\partial \ov{z}_{k}} (0) (R_{m\ov{m}q\ov{q}} + R^E_{q\ov{q}})
- \frac{1}{2\pi^2} \frac{\partial ^2 f_{\,x_0}}{\partial z_{q}
\partial \ov{z}_{l}} (0)
\Big( \frac{4}{3} R_{l\ov{k}k\ov{q}}+ R^E_{l\ov{q}}\Big).
   \end{align}
By Lemma \ref{lmt1.6}, (\ref{bk4.1}), (\ref{bk4.1a}) 
and (\ref{bk4.14}), we get
  \begin{align}\label{bk4.15}
      \cK_{2f}(0,0) + \cK_{2f}(0,0)^*
=\frac{1}{\pi^2} (\bb_{\C f}+ \bb_{Ef2})
 - \frac{1}{3\pi^2} R_{l\ov{k}k\ov{q}}
 \frac{\partial ^2 f_{\,x_0}}{\partial z_{q}\partial \ov{z}_{l}} (0) .
\end{align}
We compute now the last term in \eqref{bk4.9}.
By Remark \ref{bkt3.5}, \eqref{toe1.6} and \eqref{abk3.8},  we have
\begin{align}\label{bk4.16}
\sum_{|\alpha|=4}
\cK\Big[1, \frac{\partial ^\alpha f_{\,x_0}}{\partial Z^\alpha}(0)
  \frac{Z^\alpha}{\alpha !}\Big](0,0)
=  \frac{1}{2\pi^2}\frac{\partial ^4
f_{\,x_0}  }{\partial z_i \partial z_q\partial\ov{z}_i \partial\ov{z}_q}(0).
 \end{align}
We next turn to the computation of the sixth term in \eqref{bk4.9}. Set
\begin{equation}\label{bk4.17}
\begin{split}
K_{3f}: &= \Big(\frac{1}{6}R_{k\ov{k}m\ov{m}; \ov{u}}
- \frac{1}{3} R_{k\ov{k}m\ov{u}; \ov{m}}  \Big)
 \frac{\partial f_{\,x_0}}{\partial z_{u}}(0)
+ \Big(\frac{1}{6}R_{k\ov{k}m\ov{m}; u}- \frac{1}{3} R_{k\ov{k}u\ov{m}; m}
\Big)\frac{\partial f_{\,x_0}}{\partial \ov{z}_{u}}(0) \\
 &+ \frac{\partial f_{\,x_0}}{\partial z_{u}}(0)
\Big( \frac{1}{6} R^E_{k\ov{k} ;  \ov{u}}
-\frac{1}{2} R^E_{q\ov{u};\ov{q}}\Big)
+ \frac{1}{3} R^E_{m\ov{u}; \ov{m}}
\frac{\partial f_{\,x_0}}{\partial z_{u}}(0) \\
&+ \Big( \frac{1}{6} R^E_{k\ov{k}; u}
- \frac{1}{2} R^E_{u\ov{q};q}\Big)
\frac{\partial f_{\,x_0}}{\partial \ov{z}_{u}}(0)\
+ \frac{1}{3}\frac{\partial f_{\,x_0}}{\partial \ov{z}_{u}}(0)
R^E_{u\ov{m} ; m}.
\end{split}
\end{equation}
\begin{lemma} \label{bkt4.3} The following identity holds
with $K_{3f}$ defined in \eqref{bk4.17}:
\begin{align}\label{bk4.18}
 \cK\Big[ J_{3,\,x_0}\,, \frac{\partial f_{\,x_0}}{\partial Z_{u}}(0) Z_{u}
 \Big](0,0)
+  \cK\Big[1,  \frac{\partial f_{\,x_0}}{\partial Z_{u}}(0) Z_{u}
 J_{3,\,x_0}\Big](0,0)
=\frac{1}{\pi^2} K_{3f}.
\end{align}
\end{lemma}
\begin{proof}[Proof of Lemma \ref{bkt4.3}] Set
\begin{subequations}
\begin{align} \label{bk4.19}
B_1(b,Z):=&\Big\{ \frac{1}{6}  R_{k\ov{s}\ell\ov{q}; \ov{m}} z_{k} z_{\ell}
b_s b_q + \frac{2\pi}{3}  R_{q\ov{s}k\ov{q}; \ov{m}}  \ov{z}_{s} z_{k}\\
&-  \Big[\frac{2\pi}{15} R_{k\ov{q}\ell\ov{s}; \ov{m}}
 z_{\ell} \ov{z}_q +\frac{1}{3} R_{\ell\ov{\ell}k\ov{s} ; \ov{m}}
- \frac{1}{3} R_{k\ov{s}q\ov{m}; \ov{q}}
- \frac{2}{3} R^E_{k\ov{s}; \ov{m}} \Big] z_{k} b_s \nonumber\\
&- \frac{2\pi}{15} R_{k\ov{s}q\ov{m}; \ov{q}}  z_{k}  \ov{z}_s
-\frac{2}{3} R_{\ell\ov{\ell}q\ov{m};\ov{q}}
 -\frac{2}{3} R_{\ell\ov{\ell}q\ov{q};  \ov{m}}
+ \frac{2}{3} R^E_{q\ov{m};\ov{q}}
- 2R^E_{q\ov{q} ; \ov{m}}  \Big\} \ov{z}_m ,\nonumber\\
\label{bk4.21}
\mB_1(Z): =&\frac{2\pi^2}{5} R_{k\ov{s}\ell\ov{q}; \ov{m}} z_{k} z_{\ell}
\ov{z}_s \ov{z}_q \ov{z}_m
%-\frac{4\pi}{3}R_{s\ov{q}l\ov{s}; \ov{m}}\ov{z}_q z_{\ell}\ov{z}_m
+ \frac{8\pi}{15} R_{k\ov{s}q\ov{m}; \ov{q}}  z_{k}  \ov{z}_s
\ov{z}_m\\
&+  \frac{4\pi}{3} R^E_{k\ov{s} ;  \ov{m}}  z_k \ov{z}_{s}\ov{z}_m
-\frac{2}{3} \Big[ R_{\ell\ov{\ell}q\ov{m};\ov{q}}
 + R_{\ell\ov{\ell}q\ov{q};  \ov{m}}
-  R^E_{q\ov{m};\ov{q}}
+ 3 R^E_{q\ov{q} ; \ov{m}}  \Big] \ov{z}_m. \nonumber
\end{align}
\end{subequations}
Then by Lemma \ref{lmt1.6}, (\ref{bk3.1}), we have
\begin{equation}\label{bk4.22}
\begin{split}
B_1(b,Z) & \cP(Z,0)=\Big\{
\frac{4\pi^2}{6}  R_{k\ov{s}\ell\ov{q}; \ov{m}} z_{k} z_{\ell}\ov{z}_s \ov{z}_q
+ \frac{2\pi}{3}  R_{q\ov{s}k\ov{q}; \ov{m}}  \ov{z}_{s}z_{k}\\
&-  \Big[\frac{2\pi}{15} R_{k\ov{q}\ell\ov{s}; \ov{m}}
  z_{\ell} \ov{z}_q +\frac{1}{3} R_{\ell\ov{\ell}k\ov{s} ; \ov{m}}
- \frac{1}{3} R_{k\ov{s}q\ov{m}; \ov{q}}
- \frac{2}{3} R^E_{k\ov{s} ;  \ov{m}}\Big]  z_k \cdot  2\pi\ov{z}_s \\
&- \frac{2 \pi}{15} R_{k\ov{s}q\ov{m}; \ov{q}}  z_{k}  \ov{z}_s
- \frac{2}{3} R_{\ell\ov{\ell}q\ov{m};\ov{q}}
 -\frac{2}{3} R_{\ell\ov{\ell}q\ov{q};  \ov{m}}
+ \frac{2}{3} R^E_{q\ov{m};\ov{q}}
- 2R^E_{q\ov{q} ; \ov{m}}  \Big\} \ov{z}_{m}\cP(Z,0)\\
&= \mB_1(Z)\cP(Z,0).
\end{split}
\end{equation}
Observe that the commutation relations \eqref{bk2.66} imply that
\begin{equation*}
\begin{split}
R_{k\ov{s}\ell\ov{q}; \ov{m}} z_{k} z_{\ell} b_{s}b_{q}b_{m}
=  b_{s}b_{q}b_{m}R_{k\ov{s}\ell\ov{q}; \ov{m}} z_{k} z_{\ell}
&+ b_{s}b_{m} ( 8 R_{q\ov{s}k\ov{q}; \ov{m}}  + 4R_{k\ov{s}q\ov{m};
\ov{q}} )z_{k} \\
 &+ b_{m}(8 R_{q\ov{s}s\ov{q}; \ov{m}}  + 16R_{s\ov{s}q\ov{m};
\ov{q}} )\,.
\end{split}
\end{equation*}
By (\ref{bk2.66}), (\ref{bk3.1}) and (\ref{bk4.21}), we have
\begin{equation}\label{bk4.24}
\begin{split}
\mB_1(Z)\cP(Z,0)
=\frac{1}{\pi} \Big\{ \Big(\frac{1}{20}
R_{k\ov{s}l\ov{q}; \ov{m}} z_{k} z_{\ell}  b_q
+ \frac{2}{15} R_{k\ov{s}q\ov{m}; \ov{q}}  z_{k}
+  \frac{1}{3} R^E_{k\ov{s} ;  \ov{m}}  z_k \Big) b_{s}b_m\\
- \frac{1}{3}  \Big[ R_{\ell\ov{\ell}q\ov{m};\ov{q}}
 + R_{\ell\ov{\ell}q\ov{q};  \ov{m}}
-  R^E_{q\ov{m};\ov{q}}
+ 3 R^E_{q\ov{q} ; \ov{m}}  \Big] b_m\Big\}\cP(Z,0)\\
=\frac{1}{\pi} \Big\{ \frac{b_s b_q  b_m}{20} 
R_{k\ov{s}\ell\ov{q}\,;\,\ov{m}} z_{k} z_{\ell}
+ b_{s}b_m\Big(\frac{2}{5} R_{q\ov{s}k\ov{q}\,;\,\ov{m}}
+ \frac{1}{3}R_{k\ov{s}q\ov{m}; \ov{q}}
+  \frac{1}{3} R^E_{k\ov{s} ;  \ov{m}} \Big) z_{k} \\
+   b_m\Big( \frac{1}{15}   % (\frac{2}{5} -\frac{1}{3})
R_{s\ov{s}q\ov{q}; \ov{m}}+ %(\frac{4}{5} + \frac{8}{15} -\frac{1}{3})
R_{s\ov{s}q\ov{m}; \ov{q}}
- \frac{1}{3} R^E_{s\ov{s} ;  \ov{m}} + R^E_{q\ov{m};\ov{q}}
\Big) \Big\}\cP(Z,0).
\end{split}
\end{equation}
Note that by Theorem \ref{bkt2.17} and \eqref{bk3.1}, we have
$(\cP^\bot z_{u} \cP)(Z,0)=(\cP \ov{z}_{u} \cP)(Z,0)=0$.
Taking into account that $\cP(0,0)=1$ and 
relations \eqref{toe1.6}, \eqref{bk2.24} and \eqref{bk2.31}, we get
\begin{equation}\label{bk4.25}
\begin{split}
 \cK\Big[ J_{3,x_0}, \frac{\partial f_{\,x_0}}{\partial Z_{u}}(0) Z_{u}
 \Big](0,0)
= -( \cL^{-1} \cP^\bot\mO_3\cP \frac{\partial f_{\,x_0}}{\partial z_{u}}(0)
z_{u}\cP )(0,0) \\
-(  \cP \mO_3\cL^{-1}\cP^\bot
\tfrac{\partial f_{\,x_0}}{\partial \ov{z}_{u}}(0) \ov{z}_{u}  \cP )(0,0) .
\end{split}
\end{equation}
By Remark \ref{toet2.7}, $(\cP \mO_3\cL^{-1} \ov{z}_{u} \cP )(0,0)$ is the
adjoint of $(\cP  z_{u}\cL^{-1}\cP^\bot\mO_3\cP  )(0,0)$,
thus we will compute only the latter.
By Lemma \ref{lmt1.6}, \eqref{bk2.66}, \eqref{lm01.11} and \eqref{bk3.2},
we get as in  \eqref{bk3.3},
\begin{equation}\label{bk4.26}
\begin{split}
\mO_3=&\, \frac{1}{6}  R_{k\ov{s}\ell\ov{q}; Z} z_{k} z_{\ell}
b_s b_q + \frac{2\pi}{3}  R_{q\ov{s}k\ov{q}; Z}  \ov{z}_{s}z_{k}
+ c_s(b,b^+, Z) b^+_s\\
&-  \Big(\frac{2\pi}{15} R_{\ell\ov{q}k\ov{s}; Z}  z_{\ell} \ov{z}_q
 +\frac{1}{3} R_{\ell\ov{\ell}k\ov{s} ; Z}
+\frac{1}{3}  \big(R_{\ell\ov{q}k\ov{s}; q} z_{\ell}
-R_{q\ov{m}k\ov{s}; \ov{q}} \ov{z}_{m}\big)
- \frac{2}{3} R^E_{k\ov{s} ; Z} \Big) z_{k} b_s \\
&+  \frac{2\pi}{15} (R_{k\ov{s}\ell\ov{q}; q}z_\ell
-R_{k\ov{s}q\ov{m}\,;\, \ov{q}} \ov{z}_{m} ) z_{k}  \ov{z}_s
-\frac{2}{3} \big( R_{\ell\ov{\ell}k\ov{q} ;q} z_k 
+ R_{\ell\ov{\ell}q\ov{m};\ov{q}} \ov{z}_m\big)\\
& -\frac{2}{3} R_{\ell\ov{\ell}q\ov{q}; Z}
-\frac{2}{3} \big( R^E_{k\ov{q} ;q} z_k -R^E_{q\ov{m};\ov{q}}\ov{z}_m  \big)
- 2R^E_{q\ov{q} ; Z}, 
\end{split}
\end{equation}
and $c_s(b,b^+, Z)$ are polynomials in $b, b^+, Z$, whose precise
formula will not be used,
and $R_{k\ov{q}\ell\ov{s}; Z} $, $R^E_{k\ov{s} ; Z}$ are defined by
replacing $\frac{\partial}{\partial z_{s}}$ by $Z$ in \eqref{lm01.4}.

From Lemma \ref{lmt1.6}, \eqref{bk2.66}, \eqref{bk4.19} 
and \eqref{bk4.26} we deduce that
the only term in $\mO_3$ not containing
$b^+$ and
having total degree in $b,\ov{z}$ bigger than its degree in $z$, is $B_1(b,Z)$.
Now, Theorem \ref{bkt2.17}, Remark \ref{bkt3.5}, \eqref{toe1.4},
\eqref{bk2.66}, \eqref{bk4.22}, \eqref{bk4.24} and \eqref{bk4.26} imply
\begin{equation}\label{bk4.27}
\begin{split}
 -   (&\cP z_{u} \cL^{-1}  \cP^\bot\mO_3\cP  )(0,0)
= -  (\cP z_{u}\cL^{-1} \cP^\bot   B_{1}(b,Z)\cP  )(0,0)\\
&=- \Big(\cP z_{u}   \frac{b_m}{4 \pi^2}\Big(\frac{1}{15}
R_{s\ov{s}q\ov{q}; \ov{m}}+ R_{s\ov{s}q\ov{m}; \ov{q}}
- \frac{1}{3} R^E_{s\ov{s} ;  \ov{m}}
+ R^E_{q\ov{m};\ov{q}}\Big) \cP  \Big)(0,0)\\
&= -  \frac{1}{2 \pi^2}\Big(\frac{1}{15}
R_{s\ov{s}q\ov{q}; \ov{u}}+ R_{s\ov{s}q\ov{u}; \ov{q}}
- \frac{1}{3} R^E_{s\ov{s} ;  \ov{u}} + R^E_{q\ov{u};\ov{q}}\Big).
\end{split}
\end{equation}
By Remark \ref{bkt3.5}, (\ref{toe1.4}), (\ref{bk2.66}) and (\ref{bk4.26}), 
we also have
\begin{equation}\label{bk4.28}
\begin{split}
- ( \cL^{-1} \cP^\bot\mO_3\cP z_{u} \cP )(0,0)
= - ( \cL^{-1} \cP^\bot  B_{1}(b,Z) z_{u}\cP )(0,0)\\
= - \Big\{\cL^{-1} \cP^\bot  \big[z_{u}B_{1}(b,Z)
-2 \tfrac{\partial}{\partial b_{u}}B_{1}(b,Z)\big]\cP \Big\}(0,0) .
\end{split}
\end{equation}
Lemma \ref{lmt1.6}, (\ref{bk3.16d}), (\ref{bk3.16g}),  
(\ref{lm3.38}), (\ref{bk4.21}) and (\ref{bk4.22}) yield
\begin{equation}\label{bk4.29}
\begin{split}
- &(\cL^{-1} \cP^\bot  z_{u}B_{1}(b,Z)  \cP )(0,0)
= - (\cL^{-1} \cP^\bot  z_{u}\mB_{1}(Z)  \cP )(0,0) \\
&= \frac{11}{24\pi^2} \cdot
\frac{2}{5} \big(4 R_{s\ov{s}m\ov{u}; \ov{m}} +2 R_{s\ov{s}q\ov{q}; \ov{u}}\big)
+\frac{3}{8\pi^2} \cdot  \Big(\frac{16}{15} R_{s\ov{s}q\ov{u}; \ov{q}}
+\frac{4}{3} R^E_{s\ov{s} ;  \ov{u}}
+\frac{4}{3} R^E_{m\ov{u} ;  \ov{m}} \Big)\\
&\hspace{4cm}- \frac{1}{6 \pi^2}\Big( R_{\ell\ov{\ell}q\ov{u};\ov{q}}
 + R_{\ell\ov{\ell}q\ov{q};  \ov{u}}  -  R^E_{q\ov{u};\ov{q}}
+ 3 R^E_{q\ov{q} ; \ov{u}}  \Big) \\
&= \frac{1}{\pi^2} \Big( %(\frac{17}{15}- \frac{1}{6})
\frac{29}{30} R_{s\ov{s}m\ov{u}; \ov{m}}
%+ (\frac{11}{30}-\frac{1}{6})
+ \frac{1}{5} R_{s\ov{s}q\ov{q}; \ov{u}}
%+( \frac{1}{2}-\frac{1}{2})  R^E_{s\ov{s} ;  \ov{u}}
+\frac{2}{3} R^E_{m\ov{u} ; \ov{m}}\Big)\,.
\end{split}
\end{equation}
Moreover, by (\ref{bk4.19}),
\begin{equation}\label{bk4.30}
\begin{split}
\frac{\partial}{\partial
b_{u}}B_{1}(b,Z)  = &\,\frac{1}{3}  R_{k\ov{s}\ell\ov{u}; \ov{m}}
z_{k} z_{\ell} b_{s} \ov{z}_{m}\\
&-  \Big(\frac{2\pi}{15} R_{k\ov{q}\ell\ov{u}; \ov{m}}
  z_{\ell} \ov{z}_q +\frac{1}{3} R_{\ell\ov{\ell}k\ov{u} ; \ov{m}}
- \frac{1}{3} R_{k\ov{u}q\ov{m}; \ov{q}}
- \frac{2}{3} R^E_{k\ov{u} ;  \ov{m}}   \Big)
z_{k} \ov{z}_{m}  .
\end{split}
\end{equation}
Lemma \ref{lmt1.6}, (\ref{bk3.1}), (\ref{bk3.16d}), 
(\ref{bk3.16g}) and (\ref{bk4.30}) yield
\begin{equation}\label{bk4.31}
\begin{split}
( &\cL^{-1} \cP^\bot (\tfrac{\partial}{\partial b_{u}}B_{1}(b,Z))\cP)(0,0)  \\
&=  -\frac{3}{8\pi^2} \cdot  %(\frac{2}{3}- \frac{2}{15})
\frac{8}{15} \cdot 2
R_{s\ov{s}m\ov{u}; \ov{m}}
+ \frac{1}{4\pi^2} \Big(\frac{1}{3} R_{m\ov{q}q\ov{u} ; \ov{m}}
- \frac{1}{3} R_{m\ov{u}q\ov{m}; \ov{q}}
- \frac{2}{3} R^E_{m\ov{u} ;  \ov{m}} \Big) \\
&= -\frac{2}{5\pi^2} R_{s\ov{s}m\ov{u}; \ov{m}}
- \frac{1}{6\pi^2} R^E_{m\ov{u} ;  \ov{m}}.
\end{split}
\end{equation}
Formulas \eqref{bk4.28}--\eqref{bk4.31} entail altogether
\begin{equation}\label{bk4.32}
-\pi^2 ( \cL^{-1} \cP^\bot\mO_3\cP z_{u} \cP )(0,0)
%= (\frac{17}{15}- \frac{1}{6}-\frac{4}{5})R_{s\ov{s}m\ov{u}; \ov{m}}
%+ \frac{1}{5} R_{s\ov{s}q\ov{q}; \ov{u}}
%+\frac{1}{3}R^E_{m\ov{u} ; \ov{m}}\\
=\frac{1}{6}R_{s\ov{s}m\ov{u}; \ov{m}}
+ \frac{1}{5} R_{s\ov{s}q\ov{q}; \ov{u}}
+\frac{1}{3}R^E_{m\ov{u} ; \ov{m}} .
\end{equation}
Combining \eqref{bk4.25}, \eqref{bk4.27} with \eqref{bk4.32}, we get
\begin{equation}\label{bk4.33}
\begin{split}
\pi^2 \cK\Big[ J_{3,x_0}, \frac{\partial f_{\,x_0}}{\partial Z_{u}}(0) Z_{u}
 \Big](0,0)= \Big(\frac{1}{6}R_{s\ov{s}m\ov{u}; \ov{m}}
+ \frac{1}{5} R_{s\ov{s}q\ov{q}; \ov{u}}
+\frac{1}{3}R^E_{m\ov{u} ; \ov{m}} \Big)
\frac{\partial f_{\,x_0}}{\partial z_{u}}(0) \\
+  \Big(-\frac{1}{30} R_{s\ov{s}q\ov{q}; u}
- \frac{1}{2} R_{s\ov{s}u\ov{q}; q}
+ \frac{1}{6} R^E_{s\ov{s}; u} -\frac{1}{2}R^E_{u\ov{q}; q} \Big)
\frac{\partial f_{\,x_0}}{\partial \ov{z}_{u}}(0).
\end{split}
\end{equation}
Since $\cK\Big[1,\frac{\partial f_{\,x_0}}{\partial Z_{u}}(0) Z_{u}
J_{3,x_0} \Big]$ is the adjoint of
$\cK\Big[ J_{3,x_0}, \tfrac{\partial f_{\,x_0}}{\partial Z_{u}}(0) Z_{u}
 \Big]$, Lemma \ref{lmt1.6} yield
\begin{equation}\label{bk4.34}
\begin{split}
\pi^2 \cK\Big[1, \frac{\partial f_{\,x_0}}{\partial Z_{u}}(0) Z_{u}
J_{3,x_0} \Big](0,0)
=\frac{\partial f_{\,x_0}}{\partial \ov{z}_{u}}(0)
 \Big(\frac{1}{6}R_{s\ov{s}u\ov{m}; m}
+ \frac{1}{5} R_{s\ov{s}q\ov{q}; u}
+\frac{1}{3}R^E_{u\ov{m}; m} \Big)\\
+\frac{\partial f_{\,x_0}}{\partial z_{u}}(0)
\Big(-\frac{1}{30} R_{s\ov{s}q\ov{q}; \ov{u}}
- \frac{1}{2} R_{s\ov{s}q\ov{u}; \ov{q}}
+ \frac{1}{6} R^E_{s\ov{s}; \ov{u}} -\frac{1}{2}R^E_{q\ov{u};\ov{q}} \Big).
\end{split}
\end{equation}
Finally, \eqref{bk4.33} and \eqref{bk4.34} deliver \eqref{bk4.18}.
The proof of Lemma \ref{bkt4.3} is completed.
\end{proof}
We continue with the proof of Theorem \ref{toet4.1a}. We'll write now 
the formulas in terms of connections.
%Recall that  for $f\in \cC^{\infty}(X,\End(E))$, on $\R^{2n}$,
%by  (\ref{toe2.5}), as in (\ref{bk2.87}), we have
For $f\in \cC^{\infty}(X,\End(E))$, we obtain  
by \eqref{toe2.5} (as in \eqref{bk2.87}) the following formula 
in normal coordinates: 
\begin{align}\label{bk4.41}
(\Delta^{E} f)(Z) = -g^{ij}(\nabla^{E}_{e_{i}}\nabla^{E}_{e_{j}}
- \Gamma^{l}_{ij}\nabla^{E}_{e_{l}}) f, \quad
\nabla^{E} f = df + [\Gamma^{E}(\cdot), f].
\end{align}
By (\ref{bk2.85}),
\begin{equation}\label{bk4.42}
\begin{split}
    \nabla^{E}_{e_{i}}\nabla^{E}_{e_{i}}
    &= \frac{\partial^{2}}{\partial Z_{i}^{2}}
    + R^E(\mR,e_{i}) \frac{\partial}{\partial Z_{i}}
  +   \frac{1}{4} R^E (\mR,e_{i}) R^E (\mR,e_i)
+ \frac{2}{3} R^E_{\, ; Z} (\mR,e_i) \frac{\partial}{\partial
Z_{i}}\\
&+  \frac{1}{3} R^E_{\, ; e_{i}} (\mR,e_i)
+\frac{1}{8} \Big(2 R^E_{\, ; (Z,e_{i})} (\mR,e_i)
+ \frac{1}{3}\left \langle  R^{TX} (\mR,e_i)e_{i},
e_k\right \rangle R^E (\mR,e_k) \Big)\\& +  \cO(|Z|^3).
\end{split}
\end{equation}
Note that by \eqref{lm01.2},
$\dfrac{\partial^{2}}{\partial Z_{m}^{2}}R^E_{\, ; (Z,e_{i})} (\mR,e_i)
= 2 R^E_{\, ; (e_{m},e_{i})} (e_{m},e_i) =0$.
Thus from (\ref{alm01.0}), (\ref{bk3.2}), (\ref{bk4.42}),
 and taking into account that $R^E(e_{m}, e_{i})$ is anti-symmetric in
$m,i$, and $\ric(e_{m}, e_{i})$ is symmetric in $m,i$, we infer
\begin{multline}\label{bk4.43}
   \Big( \frac{\partial^{2}}{\partial Z_{m}^{2}}
   \nabla^{E}_{e_{i}}\nabla^{E}_{e_{i}} f\Big)(0)
= \frac{\partial^{4} f_{x_0}}{\partial Z_{m}^{2}\partial Z_{i}^{2}}(0)
+ 2 \bigg[R^E(e_{m}, e_{i}), \frac{\partial^{2}f_{x_0}}{\partial Z_{i}\partial
Z_{m}}(0)\bigg]  \\
+ \frac{1}{2} \Bigl[R^E(e_{m},e_{i}), \Bigl[R^E(e_{m},e_{i}),
f(x_{0})\Bigr] \Bigr]
+ \frac{2}{3} \bigg[R^E_{\, ; e_{m}} (e_{m},e_i) , \frac{\partial
f_{x_0}}{\partial Z_{i}}(0)\bigg]\\
= 16\frac{\partial ^4
f_{x_0}  }{\partial z_i \partial z_q\partial\ov{z}_i \partial\ov{z}_q}(0)
- 4 \Big[R^{E}_{k\ov{\ell}}, \Big[R^{E}_{\ell\ov{k}}, f(x_{0})\Big]\Big]
+ \frac{8}{3} \bigg[R^E_{m\ov{q} ; \ov{m}} , \frac{\partial
f_{x_0}}{\partial z_{q}}(0)\bigg]
-\frac{8}{3} \bigg[R^E_{q\ov{m} ; m} , \frac{\partial
f_{x_0}}{\partial \ov{z}_{q}}(0)\bigg] .
\end{multline}
By \eqref{alm01.1}, \eqref{abk2.83}, \eqref{bk2.85} and \eqref{bk4.41},
%\begin{multline}\label{bk4.44}
%-    \frac{\partial^{2}}{\partial Z_{m}^{2}}
%   (\Gamma^{\ell}_{ii}\nabla^{E}_{e_{\ell}} f)(0)
%   = - \frac{4}{3} \ric(e_{m}, e_{\ell})
%   \Big( \frac{\partial^{2}f_{x_0}}{\partial Z_{\ell}\partial
%Z_{m}}(0) + \frac{1}{2} \left[R^E(e_{m},e_{\ell}),
%f(x_{0})\right]\Big)\\
%- \frac{1}{6}\Big(4\ric_{\, ; e_{m}}(e_{m}, e_{\ell}) 
%-2 \ric_{\, ; e_{q}}(e_{q}, e_{\ell})
%+ \ric_{\, ; e_{\ell}}(e_{m}, e_{m}) \Big)  \frac{\partial f_{x_0}}{\partial
%Z_{\ell}}(0)\\
%=  - \frac{32}{3} \ric_{m\ov{\ell}}\frac{\partial ^2
%f_{x_0}  }{\partial z_{\ell}  \partial\ov{z}_{m}}(0)
%- \frac{4}{3} \Big( \big(\ric_{m\ov{\ell}; \ov{m}} 
%+ \ric_{m\ov{m}; \ov{\ell}}\big) \frac{\partial
%f_{x_0} }{\partial z_\ell}(0)  + \big(\ric_{\ell\ov{m}; m} + 
%\ric_{m\ov{m}; \ell}\big)
%\frac{\partial f_{x_0} }{\partial\ov{z}_\ell}(0)\Big)\,.
%\end{multline}
\begin{equation}\label{bk4.44}
\begin{split}
-    \frac{\partial^{2}}{\partial Z_{m}^{2}}&
   (\Gamma^{\ell}_{ii}\nabla^{E}_{e_{\ell}} f)(0)\\
   &= - \frac{4}{3} \ric(e_{m}, e_{\ell})
   \Big( \frac{\partial^{2}f_{x_0}}{\partial Z_{\ell}\partial
Z_{m}}(0) + \frac{1}{2} \left[R^E(e_{m},e_{\ell}),
f(x_{0})\right]\Big)\\
&\quad - \frac{1}{6}\Big(4\ric_{\, ; e_{m}}(e_{m}, e_{\ell}) 
-2 \ric_{\, ; e_{q}}(e_{q}, e_{\ell})
+ \ric_{\, ; e_{\ell}}(e_{m}, e_{m}) \Big)  \frac{\partial f_{x_0}}{\partial
Z_{\ell}}(0)\\
&=  - \frac{32}{3} \ric_{m\ov{\ell}}\frac{\partial ^2
f_{x_0}  }{\partial z_{\ell}  \partial\ov{z}_{m}}(0)\\
&\quad - \frac{4}{3} \Big( \big(\ric_{m\ov{\ell}; \ov{m}} 
+ \ric_{m\ov{m}; \ov{\ell}}\big) \frac{\partial
f_{x_0} }{\partial z_\ell}(0)  + \big(\ric_{\ell\ov{m}; m} + \ric_{m\ov{m}; \ell}\big)
\frac{\partial f_{x_0} }{\partial\ov{z}_\ell}(0)\Big)\,.
\end{split}
\end{equation}

Formulas \eqref{alm01.1}, \eqref{lm01.38}, \eqref{bk2.85}, \eqref{bk4.3a}
and \eqref{bk4.41}--\eqref{bk4.44} yield
\begin{equation}\label{bk4.45}
\begin{split}
((\Delta^E)^2 f&)(x_{0}) =
     - \Big( \frac{\partial^{2}}{\partial Z_{m}^{2}} \Delta^E f\Big)(0) \\
     = &\,\frac{2}{3} \ric(e_{i}, e_{q}) (
     \nabla^{E}_{e_{i}}\nabla^{E}_{e_{q}}f)(0)
     +  \Big( \frac{\partial^{2}}{\partial Z_{m}^{2}}
   \nabla^{E}_{e_{i}}\nabla^{E}_{e_{i}} f\Big)(0)
   -    \frac{\partial^{2}}{\partial Z_{m}^{2}}
   (\Gamma^{l}_{ii}\nabla^{E}_{e_{l}} f)(0)\\
     = &\,16\frac{\partial ^4
f_{x_0}  }{\partial z_i \partial z_q\partial\ov{z}_i \partial\ov{z}_q}(0)
- \frac{16}{3} \ric_{m\ov{\ell}}\frac{\partial ^2
f_{x_0} }{\partial z_{\ell}  \partial\ov{z}_{m}}(0)
- 4 \left[R^{E}_{k\ov{\ell}}\,, [R^{E}_{\ell\ov{k}}\,, f(x_{0})]\right]\\
& + \frac{8}{3} \bigg[R^E_{m\ov{q} ; \ov{m}} \,, \frac{\partial
f_{x_0}}{\partial z_{q}}(0)\bigg]
-\frac{8}{3} \bigg[R^E_{q\ov{m} ; m} \,, \frac{\partial
f_{x_0}}{\partial \ov{z}_{q}}(0)\bigg]\\
 & - \frac{4}{3} \Big( \big(\ric_{m\ov{\ell}; \ov{m}} 
+ \ric_{m\ov{m}; \ov{\ell}}\big) \frac{\partial
f_{x_0}  }{\partial z_\ell}(0)  + \big(\ric_{\ell\ov{m}; m} + \ric_{m\ov{m}; \ell}\big)
\frac{\partial f_{x_0}  }{\partial\ov{z}_\ell}(0)\Big)\,.
\end{split}
\end{equation}
By Lemma \ref{bkt4.3}, the discussion after \eqref{bk4.9},
formulas \eqref{bk4.9}, \eqref{abk4.11}, \eqref{abk4.13},
\eqref{bk4.15} and \eqref{bk4.16}, we have
\begin{equation}\label{bk4.46}
\begin{split}
  \pi^2  Q_{4,\,x_0}(f) (0,0)= &\, 2 K_{41} f(x_{0})
    +  K_{42} f(x_{0}) +f(x_{0})  K_{42}+K_{2f}+ K_{3f}\\
    &+\bb_{\C f}+ \bb_{Ef2}
     - \frac{1}{3} R_{\ell\ov{k}k\ov{q}}
     \frac{\partial ^2 f_{\,x_0}}{\partial z_{q}\partial \ov{z}_{\ell}} (0)
+ \frac{1}{2}\frac{\partial ^4
f_{\,x_0}  }{\partial z_i \partial z_q\partial\ov{z}_i
\partial\ov{z}_q}(0).
\end{split}
\end{equation}
%By Lemma \ref{lmt1.6} and \eqref{lm3.56}, we have
%\begin{align}\label{bk4.47}
%R_{i\ov{i}m\ov{m}; \ov{u}}=  R_{i\ov{i}m\ov{u}; \ov{m}},
%\quad R_{i\ov{i}m\ov{m}; u}= R_{i\ov{i}u\ov{m}; m}
%\end{align}
Note that
$\left[R^{E}_{k\ov{\ell}}\,, [R^{E}_{\ell\ov{k}}\,, f(x_{0})]\right]
= R^{E}_{k\ov{\ell}}R^{E}_{\ell\ov{k}}f(x_{0})
- 2 R^{E}_{k\ov{\ell}}f(x_{0}) R^{E}_{\ell\ov{k}}
+ f(x_{0}) R^{E}_{k\ov{\ell}}R^{E}_{\ell\ov{k}}$\,,
so by \eqref{alm01.5}, \eqref{bk4.5}, (\ref{bk4.45}) and \eqref{bk4.46},
we get \eqref{bk4.4}.
The proof of Theorem \ref{toet4.1a} is completed.

%%%%%%%%%%%%%%%%%%%%%%%%%%%%%%%%%%%%%%%%%%
\subsection{Composition of Berezin-Toeplitz operators: proofs of Theorems
\ref{toet4.6}, \ref{toet4.5}}\label{toes4.3}
%%%%%%%%%%%%%%%%%%%%%%%%%%%%%%%%%%

%\newpage

\begin{proof}[\textbf{Proof of Theorem \ref{toet4.6}}]
   By Lemma \ref{toet2.3}, we deduce as in the proof of
Lemma \ref{toet2.3}, that for $Z,Z^\prime \in T_{x_0}X$,
$\abs{Z},\abs{Z^{\prime}}<\var/4$, we have 
(cf.\,\cite[(4.79)]{MM08b}, \cite[(7.4.6)]{MM07})
\begin{align} \label{toe4.6}
p^{-n}(T_{f,\,p}\circ  T_{g,\,p})_{x_0}(Z,Z^\prime)\cong
\sum^k_{r=0}(Q_{r,\,x_0}(f,g)\cP_{x_0})(\sqrt{p}\,Z,\sqrt{p}\,Z^{\prime})
p^{-\frac{r}{2}} + \mO(p^{-\frac{k+1}{2}}),
\end{align}
%and with the notation \eqref{toe1.6},
where
\begin{align} \label{toe4.7}
Q_{r,\,x_0}(f,g)= \sum_{r_1+r_2=r}  \cK[Q_{r_1,\,x_0}(f),
Q_{r_2,\,x_0}(g)]\in\End(E)_{\,x_{0}}[Z,Z^{\prime}],
\end{align}
is a polynomial in $Z,Z^{\prime}$ with the same parity as $r$. 
%with values in $\End(E)_{\,x_{0}}$.

The existence of the expansion \eqref{toe4.30} and the expression of 
$\bb_{0,\,f,\,g}$ follow from
\eqref{toe2.15}, \eqref{toe4.6} and \eqref{toe4.7}; we get also
    \begin{align}\label{toe4.33}
	\bb_{r,\,f,\,g}(x_{0}) = Q_{2r,\,x_{0}}(f,g)(0,0).
\end{align}
By \eqref{toe4.7}, %\cite[(7.4.7)]{MM07},
\begin{align}\label{toe4.10}
Q_{2,\,x_{0}}(f,g) =   \cK[f(x_{0}), Q_{2,\,x_{0}}(g)]
+  \cK[Q_{1,\,x_{0}}(f), Q_{1,\,x_{0}}(g)]
+ \cK[Q_{2,\,x_{0}}(f), g(x_{0})].
\end{align}
Formulas (\ref{toe1.6}), (\ref{bk4.6})--(\ref{bk4.7}) yield 
\begin{align}\label{toe4.11}
\begin{split}
 \cK[Q_{2,\,x_{0}}(f), g(x_{0})] \cP
&=  (Q_{2,\,x_{0}}(f)\cP)\cP  g(x_{0})  \\
&=  \Big(  \cP \sum_{|\alpha|=2}\frac{\partial ^\alpha f_{x_{0}}}
{\partial Z^\alpha}(0) \frac{Z^\alpha}{\alpha !} \cP
-  \cL^{-1}\mO_2\cP f(x_{0}) \Big) g(x_{0}),\\
 \cK[f(x_{0}), Q_{2,\,x_{0}}(g)] \cP
&= f(x_{0})  \cP (Q_{2,\,x_{0}}(g)\cP)\\
&=  f(x_{0}) \Big(  \cP \sum_{|\alpha|=2}
\frac{\partial ^\alpha g_{x_{0}}}{\partial Z^\alpha}(0)
\frac{Z^\alpha}{\alpha !} \cP   - g(x_{0}) \cP
\mO_2\cL^{-1}\cP^\bot\Big).
\end{split}
\end{align}
Using \eqref{bk2.6}, \eqref{bk3.7}, \eqref{bk3.8}, \eqref{abk3.8}, 
\eqref{bk4.3a} and \eqref{toe4.11}, we obtain
\begin{align}\label{toe4.12} \begin{split}
    \cK[Q_{2,\,x_{0}}(f), g(x_{0})] (0,0)
    =&  \left(- \frac{1}{4\pi} (\Delta^{E}f)(x_{0})
    + \frac{1}{2} (\bb_{1} f)(x_{0})\right)g(x_{0}),\\
    \cK[f(x_{0}), Q_{2,\,x_{0}}(g)] (0,0)
	=&f(x_{0}) \left( - \frac{1}{4\pi} (\Delta^{E} g)(x_{0})
	+ \frac{1}{2} (g\bb_{1})(x_{0})\right).
\end{split}\end{align}
By  (\ref{toe1.6}), (\ref{toe2.14}) and (\ref{bk2.33}),
we have (cf.\,\cite[(7.4.14)]{MM07}),
%\cite[(7.2.21)]{MM07},
\begin{align}\label{toe4.13}
 Q_{1,\,x_{0}}(f) (Z,Z^{\prime})
 =   \cK\Big[1, \frac{\partial f_x}{\partial Z_q}(0) Z_q\Big](Z,Z^{\prime})
 = \frac{\partial f_{x_{0}}}{\partial z_i}(0)
 z_i+\frac{\partial f_{x_{0}}}{\partial \ov{z}_i}(0)
 \ov{z}^{\,\prime}_{i}.
\end{align}
Thus from (\ref{abk3.8}),  (\ref{toe4.12}),
we get as in \cite[(7.4.15)]{MM07} (cf.\,also \cite[(7.1.11)]{MM07}),
\begin{align}\label{toe4.16}
    \cK[Q_{1,\,x_{0}}(f), Q_{1,\,x_{0}}(g)] (0,0)=
    \sum_{i=1}^n \frac{1}{\pi}
    \frac{\partial f_{x_{0}}}{\partial \ov{z}_i}(0)\,
    \frac{\partial g_{x_{0}}}{\partial z_i}(0).
\end{align}
Now, \eqref{toe4.33}, \eqref{toe4.10}, \eqref{toe4.12} and \eqref{toe4.16}
imply the formula for $\bb_{1,\,f,\,g}$ given in (\ref{toe4.31}).

We prove next \eqref{toe4.32}. It suffices to consider 
$f,g\in \cC^\infty(X,\R)$, which we henceforth assume.
By \eqref{toe2.9}, we have $Q_{r,\,x_0}(1)=J_{r,\,x_0}$.
Hence, taking $f=1$ in \eqref{toe4.7}, we get
\begin{equation}\label{toe4.36}
    Q_{4,\,x}(g)=  \cK[1, Q_{4,\,x}(g)]
 + \cK[J_{2,\,x}, Q_{2,\,x}(g)]
+ \cK[J_{3,\,x}, Q_{1,\,x}(g)]+ g(x) \cK[J_{4,\,x}, 1].
\end{equation}
Taking $g=1$ in \eqref{toe4.7} yields
\begin{equation}\label{toe4.37}
    Q_{4,\,x}(f)= f(x) \cK[1 , J_{4,\,x}]
+ \cK[Q_{1,\,x}(f), J_{3,\,x}]\\
+ \cK[Q_{2,\,x}(f), J_{2,\,x}]
+ \cK[Q_{4,\,x}(f), 1].
\end{equation}
By (\ref{toe2.15}) and  (\ref{toe4.7}),  we get
\begin{equation}\label{toe4.34}
\begin{split}
Q_{4,\,x_0}&(f,g)=  \cK[ f(x_0), Q_{4,\,x_0}(g)]
+ \cK[Q_{1,\,x_0}(f), Q_{3,\,x_0}(g)]\\
&+ \cK[Q_{2,\,x_0}(f), Q_{2,\,x_0}(g)]
+ \cK[Q_{3,\,x_0}(f), Q_{1,\,x_0}(g)]+ \cK[Q_{4,\,x_0}(f), g(x_0)].
\end{split}
\end{equation}
Set
  \begin{align}\label{toe4.39}\begin{split}
 &\wi{Q}_{3,\,x_0}(g)  =  \cK\Big[1,  \frac{\partial  g_{x_{0}}}
{\partial Z_{q}}(0)  Z_{q} J_{2,\,x_0}\Big]
+\cK\Big[ J_{2,\,x_0}, \frac{\partial  g_{x_{0}}}
{\partial Z_{q}}(0)  Z_{q} \Big]\\
&\hspace{20mm}+ \cK\Big[1, \sum_{|\alpha|=3}\frac{\partial ^\alpha g_{x_{0}}}
{\partial Z^\alpha}(0)  \frac{Z^\alpha}{\alpha !}\Big],\\
&\bI_{4,f,g}= \cK\bigg[\cK\Big[1, \sum_{|\alpha|=2}
\frac{\partial ^\alpha g_{x_{0}}}
{\partial Z^\alpha}(0)  \frac{Z^\alpha}{\alpha !}\Big],\,
\cK\Big[1, \sum_{|\alpha|=2}\frac{\partial ^\alpha g_{x_{0}}}
{\partial Z^\alpha}(0)  \frac{Z^\alpha}{\alpha !}\Big]\bigg].
\end{split}\end{align}
Note that by (\ref{bk2.24}) and (\ref{bk2.31}),
we have
  \begin{align}\label{toe4.38}
 J_{3,\,x_0} =  \cK[1, J_{3,\,x_0}]+ \cK[J_{3,\,x_0}, 1].
\end{align}
and \eqref{toe2.14}, \eqref{bk2.33} together with \eqref{toe4.38} imply
\begin{align}\label{toe4.40}
    Q_{3,\,x_0}(g)  = g(x_{0}) J_{3,\,x_0} +  \wi{Q}_{3,\,x_0}(g) .
\end{align}
Since $g\in\cC^{\infty}(X,\R)$, \eqref{bk4.6} entails
\begin{align}\label{toe4.41}
Q_{2,\,x_0}(g)  = g(x_{0}) J_{2,\,x_0}
+  \cK \Big[1, \sum_{|\alpha|=2}\frac{\partial ^\alpha g_{x_{0}}}
{\partial Z^\alpha}(0)  \frac{Z^\alpha}{\alpha !} \Big].
\end{align}
By Remark \ref{bkt3.5}, \eqref{bk4.10}, \eqref{toe4.36}--\eqref{toe4.41}, we get
\begin{equation}\label{toe4.42}
\begin{split}
    Q_{4,\,x_0}(f,g)= &\,f(x_0) Q_{4,\,x_0}(g) + g(x_0)Q_{4,\,x_0}(f)
  -  f(x_0)  g(x_0) J_{4,\,x_0}\\
  &+ \cK[Q_{1,\,x_0}(f),   \wi{Q}_{3,\,x_0}(g)]
  + \cK[  \wi{Q}_{3,\,x_0}(f), Q_{1,\,x_0}(g)]
  + \bI_{4,f,g}.
\end{split}
\end{equation}
By (\ref{abk3.8}) and (\ref{bk4.8}), we get
\begin{equation}\label{toe4.43}
 \begin{split}
    (\bI_{4,f,g}&\cP)(0,0)\\
    & =
\Big\{ \cP   \sum_{|\beta|=2}\frac{\partial ^\beta f_{x_{0}}}
{\partial Z^\beta}(0)  \frac{Z^\beta}{\beta !}
\Big(\sum_{|\alpha|=2}\frac{\partial^\alpha g_{x_0}}{\partial {z}^\alpha} (0)
\frac{z^\alpha}{\alpha !}
+ \frac{1}{\pi} \frac{\partial ^2 g_{x_0}}{\partial z_{i} \partial \ov{z}_{i}} (0)
\Big) \cP\Big\}(0,0)\\
&= \frac{1}{\pi^2} \Big( \frac{1}{2}
\frac{\partial ^2 f_{x_0}}{\partial \ov{z}_{i} \partial \ov{z}_{q}} (0)
\frac{\partial ^2 g_{x_0}}{\partial z_{i} \partial z_{q}} (0)
+ \frac{\partial ^2 f_{x_0}}{\partial z_{q} \partial \ov{z}_{q}} (0)
\frac{\partial ^2 g_{x_0}}{\partial z_{i} \partial \ov{z}_{i}}
(0)\Big).
\end{split}
\end{equation}
By Remark \ref{toet2.7},
$Q_{1,\,x_0}(f)$,   $Q_{3,\,x_0}(g)$ are self-adjoint
for $f,g \in \cC^{\infty}(X,\R)$, thus by (\ref{toe4.40}),
\begin{align}\label{toe4.44}
    \cK[  \wi{Q}_{3,\,x_0}(f), Q_{1,\,x_0}(g)]
    = \cK[ Q_{1,\,x_0}(g), \wi{Q}_{3,\,x_0}(f)]^{*}.
\end{align}
Thus we only need to compute the fourth term in \eqref{toe4.42}.

An examination of \eqref{bk3.3} shows that in each term of the sum 
giving $\mO_2$, the total degree in
$b^+,z$ equals the total degree in $b,\ov{z}$.
Hence Remark \ref{bkt3.5}, \eqref{toe1.6}, \eqref{bk2.31},
\eqref{bk3.7}, \eqref{bk3.8}, \eqref{abk3.8}, \eqref{bk4.13} and \eqref{toe4.13} yield
\begin{equation}\label{toe4.48}
\begin{split}
    \cK\Big[Q_{1,\,x_0}(f)&,  \cK\Big[1,  \frac{\partial  g_{x_{0}}}
{\partial Z_{q}}(0)  Z_{q} J_{2,\,x_0}\Big]\Big] \cP(0,0)\\
&= \Big(\cP \frac{\partial f_{x_0}}{\partial \ov{z}_{u}} (0)\ov{z}_{u}
\frac{\partial g_{x_0}}{\partial z_{v}} (0) z_{v}
\cF_{2,\,x_{0}}\Big)(0,0)\\
&= \frac{1}{2\pi^{2}} \frac{\partial f_{x_0}}{\partial \ov{z}_{u}} (0)
\frac{\partial g_{x_0}}{\partial z_{v}}(0)
\Big[  \delta_{uv} (R_{s\ov{s}q\ov{q}} + R^E_{q\ov{q}}) 
-\frac{4}{3} R_{u\ov{k}k\ov{v}}- R^E_{u\ov{v}} \Big].
\end{split}
\end{equation}
Using \eqref{bk2.31}, and the formula 
$(\cP^{\bot}z_{i}\cP)(Z,0)=(\cP\ov{z}_{i}\cP)(Z,0)=0$ (cf. \eqref{toe1.4}),
 we get
\begin{equation}\label{toe4.49}
\begin{split}
   &( \cK[ J_{2,\,x_0}, \frac{\partial  g_{x_{0}}}
    {\partial Z_{q}}(0)  Z_{q} ] \cP)(Z,0) =
%(-\cL^{-1} \mO_2 \cP-\cP \mO_2 \cL^{-1}\cP^{\bot})
(\cF_{2,\,x_{0}}\frac{\partial  g_{x_{0}}}
    {\partial Z_{q}}(0)  Z_{q} \cP)(Z,0) \\
&=(-\cL^{-1} \mO_2 \cP  \frac{\partial  g_{x_{0}}}
    {\partial z_{q}}(0)  z_{q} \cP
-\cP \mO_2 \cL^{-1}\cP^{\bot} \frac{\partial  g_{x_{0}}}
    {\partial \ov{z}_{q}}(0)  \ov{z}_{q} \cP)(Z,0 ).
    \end{split}
\end{equation}
%Note that by \eqref{bk3.3}, in each term of $\mO_2$, the degree on
%$b^+,z$ and the degree on $b,\ov{z}$ are same. 
By Remark \ref{bkt3.5}, \eqref{bk2.66},
\eqref{toe1.6}, \eqref{toe4.13} and
\eqref{toe4.49}, we obtain as in \eqref{toe4.48}
\begin{equation}\label{toe4.50}
\begin{split}
    &\cK\bigg[Q_{1,\,x_0}(f),  \cK\Big[ J_{2,\,x_0}, \frac{\partial  g_{x_{0}}}
    {\partial Z_{q}}(0)  Z_{q} \Big] \bigg]\cP (0,0)\\
    &= \bigg(\cP \frac{\partial f_{x_{0}}}{\partial \ov{z}_{u}}(0) \ov{z}_{u}
\cK\Big[ J_{2,\,x_0}, \frac{\partial g_{x_{0}}}
    {\partial Z_{q}}(0)  Z_{q} \Big] \cP\bigg) (0,0)\\
&=- \Big(\cP \frac{\partial f_{x_{0}}}{\partial \ov{z}_{u}}(0)
\ov{z}_{u} \cL^{-1} \mO_2 \cP  \frac{\partial  g_{x_{0}}}
    {\partial z_{q}}(0)  z_{q} \cP\Big) (0,0)=0,
    \end{split}
\end{equation}
since by \eqref{toe1.4} and \eqref{bk3.10} we have
$(\cP \ov{z}_{u} b^{\alpha}z^{\beta}\cP)(Z,0)=0$
for $|\alpha|\geqslant1$.

Note that for homogeneous degree $3$ polynomials $H$ in $Z$ 
the analogue of formula \eqref{bk3.15} holds for $(H\cP)(Z,0)$.   
Using this analogue together with
\eqref{toe1.4} and \eqref{bk3.10} we obtain
\begin{align}\label{abk4.8b}
&(\cP H \cP)(Z,0)=
\Big(\sum_{|\alpha|=3}\frac{\partial^3 H}{\partial {z}^\alpha}
\frac{z^\alpha}{\alpha !}
+\frac{z_q}{\pi} \frac{\partial ^3 H}
{\partial z_{q}\partial z_{i} \partial \ov{z}_{i}}
\Big) \cP(Z,0).
\end{align}
%we get \eqref{abk4.8b} by using the similar equation 
%\eqref{bk3.15} for $(H\cP)(Z,0)$.
Finally, \eqref{toe1.6}, \eqref{abk3.8},
 \eqref{toe4.13}, \eqref{abk4.8b} and the equality $\cP(0,0)=1$ imply
\begin{align}\label{toe4.51}
\cK\Big[Q_{1,\,x_0}(f),
\cK\Big[1, \sum_{|\alpha|=3}\frac{\partial ^\alpha g_{x_{0}}}
{\partial Z^\alpha}(0)  \frac{Z^\alpha}{\alpha !}\Big]\Big](0,0)
= \frac{1}{\pi^2} \frac{\partial  f_{x_{0}}}
    {\partial \ov{z}_{u}}(0) \frac{\partial^3 g_{x_{0}}}
    {\partial z_{u}\partial z_{i}\partial \ov{z}_{i}}(0) .
\end{align}
By Lemma \ref{lmt1.6}, \eqref{lm01.38},
\eqref{abk2.83} and \eqref{bk4.41} for $E=\C$, we get
\begin{equation}\label{toe4.53}
\frac{\partial}{\partial z_{q}} (\Delta g) (0)
= - 4  \frac{\partial^3 g_{x_{0}}}
    {\partial z_{q}\partial z_{i}\partial \ov{z}_{i}}(0)
+ \frac{4}{3} \ric_{q\ov{\ell}} \frac{\partial g_{x_{0}}}{\partial z_{\ell}}(0).
\end{equation}
Lemma \ref{lmt1.6}, \eqref{toe4.39} and
\eqref{toe4.48}--\eqref{toe4.53} entail
\begin{equation}\label{toe4.54}
\begin{split}
    &\pi^2\cK[Q_{1,\,x_0}(f),   \wi{Q}_{3,\,x_0}(g)](0,0)\\
&= \frac{\partial  f_{x_{0}}}
    {\partial \ov{z}_{u}}(0) \frac{\partial^3 g_{x_{0}}}
    {\partial z_{u}\partial z_{i}\partial \ov{z}_{i}}(0)
 +\frac{1}{2} \frac{\partial f_{x_0}}{\partial \ov{z}_{u}} (0)
\frac{\partial g_{x_0}}{\partial z_{v}}(0)
\Big[ \delta_{uv} (R_{s\ov{s}q\ov{q}} + R^E_{q\ov{q}})
- \frac{4}{3} R_{u\ov{k}k\ov{v}}-   R^E_{u\ov{v}} \Big]\\
&= - \frac{1}{8} \langle \ov{\partial} f, \partial \Delta g\rangle
+  \frac{1}{4} \langle \ov{\partial} f, \partial g\rangle
( R_{s\ov{s}q\ov{q}} + R^E_{q\ov{q}})
- \frac{1}{8}  \langle  \ov{\partial} f\wedge \partial g, R^E\rangle_{\om}.
\end{split}
\end{equation}
From  \eqref{toe4.44} and \eqref{toe4.54}, we have
\begin{equation}\label{toe4.55}
\begin{split}
    \pi^2\cK[ & \wi{Q}_{3,\,x_0}(f),   Q_{1,\,x_0}(g)](0,0)\\
&= - \frac{1}{8}\langle \ov{\partial} \Delta f, \partial g\rangle
+ \frac{1}{4} \langle \ov{\partial} f, \partial g\rangle
( R_{s\ov{s}q\ov{q}} + R^E_{q\ov{q}})
- \frac{1}{8}  \langle  \ov{\partial} f\wedge \partial g, R^E\rangle_{\om}.
\end{split}
\end{equation}
By \eqref{bk2.6}, \eqref{bk4.3a},  \eqref{toe4.42},
\eqref{toe4.43}, \eqref{toe4.54} and \eqref{toe4.55},
we get (\ref{toe4.32}). The proof of Theorem \ref{toet4.6} is completed.
\end{proof}

\begin{proof}[\textbf{Proof of Theorem \ref{toet4.5}}]
The existence of the expansion \eqref{toe4.2} and formula
$C_0(f,g)=fg$  were established in \cite[Th.\,1.1]{MM08b} (cf.\,also \cite[Th.\,7.4.1]{MM07})
in general symplectic settings.

By \eqref{toe4.2}, \eqref{toe2.13}, \eqref{bk4.5},
\eqref{toe4.6} and \eqref{toe4.33}, we obtain (cf.\;also \cite[(7.4.9)]{MM07}),
\begin{equation}\label{toe4.8}
C_1(f,g)= (Q_{2,x}(f,g)- Q_{2,x}(fg))(0,0) = \bb_{1,\,f,\,g}-\bb_{1,\,fg} .
\end{equation}
Hence \eqref{toe2.3}, \eqref{bk4.3}, \eqref{toe4.31} and \eqref{toe4.8}
yield the formula for $C_1(f,g)$ given in \eqref{toe4.3}.
Moreover, \eqref{toe4.2}, \eqref{bk4.5} and \eqref{toe4.33} imply 
the formula for $C_2(f,g)$ from \eqref{toe4.3}.

We will prove \eqref{toe4.3a} now. 
%We still denote by $\nabla^{TX}$ the connection on 
%$\Lambda(T^*X)$ induced by the Levi-Civita connection $\nabla^{TX}$. 
%Let $T^{(1,0)}X$ be the holomorphic tangent bundle on $X$ 
%and $T^{*(1,0)}X$ its dual bundle.
Let $\{e_i\}$ be an orthonormal frame of $(TX, g^ {TX})$, 
and $\{w_i\}$ be an orthonormal frame of 
$T^{(1,0)}X$.
Let $\square= \ov{\partial}^*\ov{\partial}+ \ov{\partial}\,\ov{\partial}^*$ 
be the Kodaira Laplacian on  $\Lambda(T^*X)\otimes _{\R}\C$, 
and let $\Delta$ be the Bochner Laplacian on  
$\Lambda(T^*X)\otimes _{\R}\C$
%= \Lambda(T^{*(0,1)}X)\widehat{\otimes}\Lambda(T^{*(1,0)}X)$
associated with the connection $\nabla^{\Lambda(T^{*}X)}$ 
on $\Lambda(T^*X)$
induced by $\nabla^{TX}$ (cf.\,\eqref{toe2.3}). 
Let $R^{\Lambda(T^{*(1,0)}X)}$
be the curvature of the holomorphic Hermitian connection 
on $\Lambda(T^{*(1,0)}X)$.
By Lichnerowicz formula \cite[Remark 1.4.8]{MM07}, 
and \eqref{bk4.2a}, we have
\begin{align}\label{toe4.61}\begin{split}
&R^{\Lambda(T^{*(1,0)}X)}= 
- \langle R^{TX}w_l,\ov{w}_k\rangle w^l\wedge i_{w_k},\\
&2\square =\Delta - R^{\Lambda(T^{*(1,0)}X)}(w_q,\ov{w}_q)
+ ( 2 R^{\Lambda(T^{*(1,0)}X)} + \ric) (w_l,\ov{w}_k) 
\ov{w}^k \wedge i_{\ov{w}_l}.
\end{split}\end{align}
Since $(X,\om,J)$ is K{\"a}hler, $\square$ commutes 
with the operators $\partial,\ov{\partial},d$ (cf.\,\cite[Cor.\,1.4.13]{MM07}),
 and \eqref{toe4.61} shows that $2\square f=\Delta f$ 
 for any $f\in\cC^\infty(X)$. From Lemma \ref{lmt1.6}, 
 \eqref{bk4.2a} and \eqref{toe4.61}, we have for any $f\in\cC^\infty(X)$:
\begin{align}\label{toe4.62}\begin{split}
\Delta \partial f =&  \,\partial \Delta f - \ric(\cdot, \ov{w}_k) \, w_k(f),\\
\Delta \ov{\partial} f=& \, \ov{\partial} \Delta f 
- \ric(\cdot, w_k) \, \ov{w}_k(f),\\
\Delta d f=&  \, d \Delta f 
- \ric(\cdot, e_j)  \, e_j(f).
\end{split}\end{align}
Thus \eqref{toe2.3}, \eqref{toe4.62} yield for any $f,g\in \cC^\infty(X)$:
\begin{align}\label{toe4.63}\begin{split}
&\Delta(fg)= g \Delta f + f\Delta g -2 \langle d f,dg\rangle,\\ 
&\Delta\langle \partial f,\ov{\partial}g\rangle
=\langle \Delta\partial f,\ov{\partial}g\rangle 
+\langle \partial f,\Delta\ov{\partial}g\rangle
-2\langle \nabla^{T^{*}X}\partial f,\nabla^{T^{*}X}\ov{\partial}g\rangle\\
&= \langle \partial \Delta f,\ov{\partial}g\rangle 
+\langle \partial f,\ov{\partial}\Delta g\rangle
-2\langle \nabla^{T^{*}X}\partial f,\nabla^{T^{*}X}\ov{\partial}g\rangle
-2 \ric(w_m,\ov{w}_q)\, \ov{w}_m(g) w_q(f). 
\end{split}\end{align}
Using \eqref{toe2.3}, \eqref{toe4.62} and \eqref{toe4.63}, we infer
\begin{align}\label{toe4.64}\begin{split}
\Delta^2 (fg)=& f\Delta^2g + g\Delta^2 f + 2 (\Delta f)\Delta g 
- 4 \langle d\Delta  f,d g\rangle
- 4 \langle d f,d\Delta g\rangle\\
&+ 4 \langle \nabla^{T^{*}X}d f,\nabla^{T^{*}X}d g\rangle
+ 4 \ric(e_i, e_j)  \, e_i(g)  e_j(f).
\end{split}\end{align}

We examine now closely the expression of $\pi^2 C_2(f,g)$ given 
by \eqref{toe4.3}. Using 
\eqref{abk4.4}, \eqref{toe4.32}, \eqref{toe4.3},
we see that the term of differential order $0$ in $f,g$ 
from the expression of $\pi^2 C_2(f,g)$ is zero, and 
the term of total differential order $2$ in $f,g$,
disregarding the term involving $R^E$, in $\pi^2  C_2(f,g)$ is 
\begin{align}\label{toe4.67}\begin{split}
C_{22}= \frac{\sqrt{-1}}{8} \left\langle \ric_\om, 
- f\partial\ov{\partial}g -g \partial\ov{\partial}f 
+\partial\ov{\partial}(fg)\right\rangle.
%= \frac{1}{2}\ric_{m\ov{q}} 
%\Big( \frac{\partial f_{x_{0}}}{\partial z_q} 
%\frac{\partial g_{x_{0}}}{\partial \ov{z}_m} 
%+ \frac{\partial f_{x_{0}}}{\partial \ov{z}_m} 
%\frac{\partial g_{x_{0}}}{\partial z_q} \Big),
\end{split}\end{align}
The term of total differential order $4$ in $f,g$ 
in the expression of $\pi^2  C_2(f,g)$ is 
\begin{align}\label{toe4.65}\begin{split}
C_{24}=& \frac{1}{32}f \Delta^2 g + \frac{1}{32}g \Delta^2 f
 -\frac{1}{8}\langle \ov{\partial} f, \partial \Delta g\rangle
- \frac{1}{8}\langle \ov{\partial} \Delta f, \partial g\rangle
+ \frac{1}{16} (\Delta f) \Delta g\\
&+ \frac{1}{8} \langle D^{0,1}\ov{\partial} f, D^{1,0}\partial g\rangle
-  \frac{1}{32} \Delta^2 (fg) 
-  \frac{1}{8} \Delta \langle\partial f,\ov{\partial} g\rangle.
\end{split}\end{align}
By \eqref{toe4.63}, \eqref{toe4.64}, \eqref{toe4.65} and by the formula 
$\langle D^{1,0}\ov{\partial} f, D^{0,1}\partial g\rangle
=\langle D^{0,1}\partial f, D^{1,0}\ov{\partial} g\rangle$, we get 
\begin{align}\label{toe4.66}\begin{split}
C_{24}=\frac{1}{8} \langle D^{1,0}\partial f, D^{0,1}\ov{\partial} g\rangle
+\frac{\sqrt{-1}}{8} \left\langle \ric_\om, 
\partial f\wedge \ov{\partial}g 
-\partial g \wedge\ov{\partial}f \right\rangle.
%+  \frac{1}{2}\ric_{m\ov{q}} 
%\Big( \frac{\partial f_{x_{0}}}{\partial z_q} 
%\frac{\partial g_{x_{0}}}{\partial \ov{z}_m} 
%- \frac{\partial f_{x_{0}}}{\partial \ov{z}_m} 
%\frac{\partial g_{x_{0}}}{\partial z_q} \Big).
\end{split}\end{align}

Finally, by inspecting \eqref{abk4.4}, \eqref{toe4.32}, \eqref{toe4.3}, 
we see that the term involving $R^E$ in the expression of 
$\pi^{2} C_2(f,g)$ is 
\begin{equation}\label{toe4.68}
\begin{split} 
\frac{\sqrt{-1}}{8} \Big(-f\Delta g - g \Delta f +
\Delta(fg)\Big)R^E_{\Lambda}
+ \frac{1}{4}\langle g\partial  \ov{\partial} f +
f\partial\ov{\partial}g - \partial\ov{\partial}(fg), R^E\rangle_{\om}\\
+ \frac{\sqrt{-1}}{4}\langle \ov{\partial} f,
\partial g\rangle R^E_{\Lambda}
-  \frac{1}{4} \langle  \ov{\partial} f\wedge \partial g,
R^E\rangle_{\om}
+\frac{\sqrt{-1}}{4} \langle  \partial f,\ov{\partial} g\rangle  R^E_{\Lambda}.
\end{split}
\end{equation}
Combining \eqref{toe4.67}, \eqref{toe4.66} and  \eqref{toe4.68}, 
 we get \eqref{toe4.3a}. 
The proof of Theorem \ref{toet4.5} is completed.
\end{proof}

%\newpage
%%%%%%%%%%%%%%%%%%%%%%%%%%%%%%%%%%%%%%%%
%%%%%%%%%%%%%%%%%%%%%%%%%%%%%%%%%%%
\section{Donaldson's $Q$-operator}\label{toes5}
%%%%%%%%%%%%%%%%%%%%%%%%%%%%%%%%%%%%%%%%%%

In this section
we study the asymptotics of the sequence of operators introduced by 
Donaldson \cite{D09}.  We suppose henceforth that $E=\C$.
%We to study the asymptotics when $p\to\infty$ of the sequence of operators 
%$Q_p$ introduced by Donaldson \cite{D09}.
Set ${\rm Vol}(X,dv_X):= \int_X dv_X$. Following \cite[\S 4]{D09}, set
\begin{align}\label{1n1}
K_p(x,x'):= |P_{p}(x,x')|^2_{h^{L^p}_x\otimes h^{L^{p*}}_{x'}}, \quad
R_p:= (\dim H^0(X, L^p))/ {\rm Vol}(X, dv_X).
\end{align}
Let $K_p$, $Q_{K_p}$ be the integral operators associated to $K_p$, 
defined for $f\in \cC^\infty (X)$ by
\begin{equation}\label{1n2}
(K_p f)(x):=   \int_X K_p(x,y) f(y)dv_X(y),
\quad Q_{K_p} (f) =  \frac{1}{R_p} K_p f.
\end{equation}
%By \cite[Th.\,26]{LM09}, the operator $Q_{K_p}$ has 
%the following asymptotics
%\footnote{The right hand side of the second equation of \cite[(33)]{LM09},
%and \cite[(34)]{LM09}, should be modified as follows:
%$C p^{-1/2} \left|f\right|_{\cC^{m+1}(X)}$
%  or  $\frac{C}{p} \left|f\right|_{\cC^{m+2}(X)}$.
%}:
%\begin{align*}
%\begin{split}
%& \left| Q_{K_p}f- \frac{{\rm Vol}(X, \nu)}{{\rm Vol}(X, dv_X)}\eta f
%\right|_{\cC^m(X)}
%\leqslant C p^{-1/2} \left|f\right|_{\cC^{m+1}(X)}
% \text{ or } \frac{C}{p}    \left|f\right|_{\cC^{m+2}(X)} \,,
%\end{split}
%\end{align*}
%where $d\nu$ is a volume form on $X$ and $\eta$ is the positive function on
%$X$ defined by $dv_X =\eta d\nu$.
Recall that, just as the Bergman kernel appears when comparing 
a K{\"a}hler metric $\omega$ to its algebraic approximations $\omega_p$ 
(i.e.\ pull-backs of the Fubini-Study metrics by the Kodaira embeddings), 
the operators $Q_{K_p}$ appear when one
relates infinitesimal deformations of the metric $\omega$ 
to the corresponding deformations of the approximations $\omega_p$\,.
The asymptotics of the operator $K_p$ were obtained 
in \cite[Th.\,26]{LM09}\footnote{Note that in the present context 
\cite[Th.\,26]{LM09} should be modified as follows:
$\lim_{p\to \infty}Q_{K_p}f=\frac{{\rm Vol}(X, \nu)}{{\rm Vol}(X, 
dv_X)}\,\eta f$ in $\cC^{m}(X)$, or
\begin{align*}
\begin{split}
& \left| Q_{K_p}f- \frac{{\rm Vol}(X, \nu)}{{\rm Vol}(X, dv_X)}\,\eta f
\right|_{\cC^m(X)}
\leqslant C p^{-1/2} \left|f\right|_{\cC^{m+1}(X)}
 \text{ or } Cp^{-1}    \left|f\right|_{\cC^{m+2}(X)},
\end{split}
\end{align*}
since the right hand side of the second equation of \cite[(33)]{LM09}
and \cite[(34)]{LM09} should read
as convergence  in  $\cC^{m}(X)$ without the speed, or 
$C p^{-1/2} \left|f\right|_{\cC^{m+1}(X)}$
or  $Cp^{-1} \left|f\right|_{\cC^{m+2}(X)}$.}
and used in \cite{Fine08}. The following result refines \cite[Th.\,26]{LM09} 
and is applied in the recent paper \cite{Fine10}.
%We consider now the asymptotics of the operator $K_p$:
\begin{thm} \label{nt1} For every $m\in \N$,
 there exists  $C>0$ such that for any
$f\in \cC^\infty (X)$, $p\in \N^{*}$,
\begin{align}\label{1n4}
\begin{split}
& \left|\frac{1}{p^n} K_p f - f + \frac{1}{8\pi  p} (-\br  f + 2\Delta f)
\right|_{\cC^m(X)}
\leqslant C p^{-3/2}   \left|f\right|_{\cC^{m+3}(X)}
\text{ or } C p^{-2}   \left|f\right|_{\cC^{m+4}(X)}.
\end{split}
\end{align}
\end{thm}

\begin{proof}
By \eqref{toe2.9} with $Z=0$, \eqref{bk2.24}, \eqref{bk2.31} and \eqref{bk2.33}, we get
\begin{equation}\label{n50}
\begin{split}
\left |\Big(\frac{1}{p^{2n}} K_{p,x_0}(0,Z')\kappa_{x_0}(Z')
- \Big(1+  \sum_{r=2}^k p^{-r/2}  J^\prime_{r} (0,\sqrt{p} Z')\Big)
e ^{-\pi p |Z^\prime|^2} \Big)\right|_{\cC^m(X)}\\
\leqslant C p^{-(k+1)/2}
(1+|\sqrt{p} Z'|)^N
\exp (- C_0 \sqrt{p} |Z'|)+ \cO(p^{-\infty}),
\end{split}
\end{equation}
with
\begin{align}\label{n51}
J'_2(0,Z') =(J_2 + \ov{J_2})(0,Z') .
\end{align}
Now we have the analogue of \cite[(32)]{LM09},
\begin{equation}\label{n52}
\begin{split}
\left|p^{-n}  K_{p} f
- p^n \int_{|Z^\prime|\leqslant \var}
\Big(1+  \sum_{r=2}^k p^{-r/2}  J^\prime_{r} (0,\sqrt{p} Z')\Big)
e^{-\pi p |Z^\prime|^2}f_{x_0}(Z^\prime)
dv_{T_{x_0}X} (Z^\prime)\right|_{\cC^m(X)}\\
\leqslant C p^{-(k+1)/2} \left|f\right|_{\cC^m(X)}.
\end{split}
\end{equation}
But as in the proof of \cite[Th\,.2.29\,(2)]{BeGeVe}, we get
\begin{align}\label{n53}
\begin{split}
&\left| p^n \int_{|Z^\prime|\leqslant \var} J^\prime_{r} (0,\sqrt{p} Z')
e^{-\pi p |Z^\prime|^2}f_{x_0}(Z^\prime)
dv_{T_{x_0}X}(Z^\prime)\right|_{\cC^m(X)}
\leqslant C \left|f\right|_{\cC^m(X)}.
\end{split}\end{align}
Moreover
\begin{equation}\label{n54}
\begin{split}
\left| p^n \int_{|Z^\prime|\leqslant \var}
e^{-\pi p |Z^\prime|^2}f_{x_0}(Z^\prime)
dv_{T_{x_0}X}(Z^\prime)- f(x_{0})
+ \frac{1}{4\pi p}(\Delta f) (x_{0})  \right|_{\cC^m(X)}\\
\leqslant C p^{-3/2} \left|f\right|_{\cC^{m+3}(X)}
\text{  or } C p^{-2} \left|f\right|_{\cC^{m+4}(X)}.
\end{split}
\end{equation}
Finally, %we need to compute
%$\int_{Z'\in \C^n} (J_2 + \ov{J_2})(0,Z') |\cP|^2 (0, Z') dZ'
%= 2 \Re \int_{Z'\in \C^n} \ov{J_2}(0,Z') |\cP|^2 (0, Z') dZ'$.
 by  \eqref{bk3.8} (cf.\,also \cite[(4.1.110)]{MM07}), we have
\begin{equation}\label{n55}
\begin{split}
 \int_{Z'\in \C^n} \ov{J_2}(0,Z') |\cP|^2(0, Z') dZ'
=  \int_{Z'\in \C^n} \cP (0, Z')J_2(Z', 0) \cP (Z', 0)dZ'\\
=  (\cP J_2 \cP)(0,0)
= - (\cP   \mO_2\cL^{-1}\cP^\bot)(0,0)
%=  - (\cL^{-1} \cP^\bot\mO_2 \cP)(0,0)^*
= \frac{1}{16\pi} \br. %= \frac{1}{2} b_1(x_0).
\end{split}
\end{equation}
Thus
\begin{equation}\label{n56}
\begin{split}
\left| p^n \int_{|Z^\prime|\leqslant \var} J'_2(0,\sqrt{p} Z')
e^{-\pi p |Z^\prime|^2}f_{x_0}(Z^\prime)dv_{T_{x_0}X}(Z^\prime)
- \big(\frac{ \br}{8\pi} f\big) (x_{0})  \right|_{\cC^m(X)}\\
\leqslant C p^{-1/2} \left|f\right|_{\cC^{m+1}(X)}
\text{  or } C p^{-1} \left|f\right|_{\cC^{m+2}(X)}.
\end{split}
\end{equation}
The proof of Theorem \ref{nt1} is completed.
\end{proof}

\end{document}